\documentclass[12pt]{amsart}

\usepackage{hyperref}
\usepackage[dvips]{graphicx}

\newtheorem{theorem}{Theorem}[section]
\newtheorem{lemma}[theorem]{Lemma}
\newtheorem{proposition}[theorem]{Proposition}

\newtheorem{corollary}[theorem]{Corollary}

\theoremstyle{definition}
\newtheorem{definition}[theorem]{Definition}
\newtheorem{notation}[theorem]{Notation}
\newtheorem{example}[theorem]{Example}

\theoremstyle{remark}
\newtheorem{remark}[theorem]{Remark}

\newcommand{\NN}{ {\mathbb N} }

\newcommand{\CC}{{\mathbb C}}

\newcommand{\cP}{ {\mathcal P} }
\newcommand{\cPP}{{\mathcal{PS}}}

\newcommand{\GG}{{\mathcal G}}
\newcommand{\cV}{{\mathcal V}}
\newcommand{\cW}{{\mathcal W}}
\newcommand{\cM}{{\mathcal M}}
\newcommand{\KK}{{\kappa}}
\newcommand{\cU}{{\mathcal U}}
\newcommand{\UU}{\mathcal{U}}

\newcommand{\cX}{{\mathcal X}}

\newcommand{\ff}{\varphi}

\newcommand{\tr}{\mathrm{tr}}
\newcommand{\Tr}{\mathrm{Tr}}
\newcommand{\EE}{\mathrm{E}}

\newcommand{\cB}{\mathcal{B}}
\newcommand{\cA}{\mathcal{A}}

\newcommand{\SNC}{S_{NC}}

\newcommand{\NC}{NC}

\newcommand{\Wg}{\mathrm{Wg}}

\newcommand{\cS}{\mathcal{S}}

 \newcommand{\Sur}{\mathcal{SS}}

\newcommand{\cc}{k}

\newcommand{\cR}{\mathcal{R}}

\newcommand{\cov}{\mathrm{cov}}
\newcommand{\kk}{\kappa}
\newcommand{\moeb}{\mathrm{\text{M\"ob}}}
\newcommand{\cWg}{C}
\newcommand{\tovert}{\to}
\newcommand{\IZ}{\mathrm{IZ}}
\newcommand{\Young}{\mathbf{Y}}

\newcommand{\ab}{\allowbreak}
\newcommand{\ol}{\overline}
\newcommand{\ie}{\textit{i.e.}\,}
\newcommand{\wh}{\widehat}

\begin{document}

\title[Fluctuations of Random Matrices]{Second Order Freeness \\and Fluctuations of Random
Matrices:\\
III. Higher order freeness and free cumulants}

\author[B. Collins]{Beno\^{\i}t Collins $^{\dagger}$}
\address{%
CNRS, Institut Camille Jordan,
Universit\'e Claude Bernard Lyon 1,
43 boulevard du 11 novembre 1918,
69622 Villeurbanne Cedex,
France}
\email{collins@math.univ-lyon1.fr}
\thanks{$\dagger$ Research supported by JSPS and COE 
postdoctoral fellowships}


\author[J. A. Mingo]{James A. Mingo $^{(*)}$}
\address{Queen's University, Department of Mathematics and Statistics,
Jeffery Hall, Kingston, ON, K7L 3N6, Canada}
\email{mingo@mast.queensu.ca}
\thanks{$^*$ Research supported by Discovery Grants and a Leadership
Support Initiative Award from the Natural Sciences and Engineering
Research Council of Canada}

\author[P. \'Sniady]{Piotr \'Sniady $^{(\ddagger)}$}
\thanks{$\ddagger$ 
Research supported by MNiSW (project 1 P03A 013 30), EU Research Training
Network ``QP-Applications", (HPRN-CT-2002-00279) and by
European Commission Marie Curie Host Fellowship for the
Transfer of Knowledge ``Harmonic Analysis, Nonlinear Analysis
and Probability" (MTKD-CT-2004-013389).}
\address{Instytut Matematyczny, Uniwersytet Wroclawski,
 pl~Grunwald\-zki~2/4, 50-384 Wroclaw, Poland}
\email{Piotr.Sniady@math.uni.wroc.pl}

\author[R. Speicher]{\hbox{Roland Speicher
$^{(*)(\P)}$}}
\thanks{$^\P\,$Research supported by a Premier's
  Research Excellence Award from the Province of Ontario and a Killam Fellowship
  from the Canada Council for the Arts}
\address{Queen's University, Department of Mathematics and Statistics,
Jeffery Hall, Kingston, ON, K7L 3N6, Canada}
\email{speicher@mast.queensu.ca}

\begin{abstract}
We extend the relation between random matrices and free
probability theory from the level of expectations to the level
of all correlation functions (which are classical cumulants of
traces of products of the matrices). We introduce the notion of
``higher order freeness" and develop a theory of corresponding
free cumulants. We show that two independent random matrix
ensembles are free of arbitrary order if one of them is
unitarily invariant. We prove R-transform formulas for second
order freeness. Much of the presented theory relies on a
detailed study of the properties of ``partitioned permutations".
\end{abstract}

\maketitle

\section{Introduction}

Random matrix models and their large dimension behavior have been an important subject of study
in Mathematical Physics and Statistics since Wishart and Wigner. Global fluctuations of the
eigenvalues (that is, linear functionals of the eigenvalues) of random matrices have been
widely investigated in the last decade; see, e.g., \cite{MR1487983, Diaconis, Radulescu,
MingoNica2004annular, HigherOrderFreeness2}. Roughly speaking,
the trend of these investigations is that for a wide class of
converging random matrix models, the non-normalized trace
asymptotically behaves like a Gaussian variable whose variance
only depends on macroscopic parameters such as moments. The
philosophy of these results, together with the freeness results
of Voiculescu served as a motivation for our series of papers on
second order freeness.

One of the main achievements of the free probability theory of
Voiculescu \cite{Voiculescu1991,VoiculescuDykemaNica} was an
abstract description via the notion of ``freeness" 
of the expectation of these Gaussian variables
for a large class of
non-commuting tuples of random matrices. 

In the previous articles of this series
\cite{HigherOrderFreeness1,HigherOrderFreeness2} we showed that for
many interesting ensembles of random matrices an analogue of the
results of Voiculescu for
expectations holds also true on the level of variances as well; thus pointing in the
direction that the structure of random matrices and the fine structure of their
eigenvalues can be studied in much more detail
by using the new concept of ``second order freeness".
One of the main obstacles for such a detailed study
was the absence of an effective
machinery for doing concrete calculations in this framework. Within
free probability theory of first order, such a machinery was provided
by Voiculescu with the concept of the $R$-transform, and by Speicher
with the concept of free cumulants; see, e.g., \cite{VoiculescuDykemaNica, TheBook}.

One of the main achievements of the present article is to develop a theory of
second order cumulants (and show that
the original
definition of second order freeness from Part I of
this series \cite{HigherOrderFreeness1} is equivalent to the vanishing
of mixed second order cumulants)
and provide the corresponding $R$-transform
machinery.

In Section \ref{section2} we will give a more detailed (but still quite condensed) survey of
the connection between Voiculescu's free probability theory and random matrix
theory. We will there also provide the main motivation, notions and concepts for our extension
of this theory to the level of fluctuations (second order), as well as the statement of
our main results concerning second order cumulants and $R$-transforms.

Having first and second order freeness it is, of course, a natural question whether
this theory can be generalized to higher orders. It turns out that this is the case,
most of the general theory is the same for all orders. So we will in this paper
consider freeness of all orders from the very beginning and develop a general
theory of higher order freeness and higher order cumulants. Let us, however,
emphasize that first and second order freeness seem to be more important than the higher
order ones. Actually,  we can prove some of
the most important results (e.g. the $R$-transform
machinery) only for first and second order,
mainly because of the complexity of the underlying combinatorial objects.

The basic combinatorial notion behind the (usual) free cumulants
are non-crossing partitions. Basically, passage to higher order
free cumulants corresponds to a change to multi-annular non-crossing
permutations \cite{MingoNica2004annular}, or more general
objects which we call ``partitioned permutations".
For much of the conceptual
framework there is no difference between different levels of
freeness, however for many concrete questions it seems that
increasing the order makes some calculations much harder. This
relates to the fact that $n$-th order freeness is described in
terms of planar permutations which connect points on $n$ different
circles.
Whereas enumeration of all non-crossing permutations in the case of
one circle is quite easy, the case of two circles gets more
complicated, but is still feasible; for the case of three or more circles,
however, the answer does not seem to be of a nice compact form.

In the present paper we develop the notion and combinatorial machinery for freeness of all
orders by a careful analysis of the main example: unitarily invariant random matrices. We start
with the calculation of mixed correlation functions for random matrices and use the structure
which we observe there as a motivation for our combinatorial setup. In this way the concept of
partitioned permutations and the moment--cumulant relations appear quite canonically.

We want to point out that even though our notion of second and
higher order freeness is modeled on the situation found
for correlation functions of random matrices, this notion and
theory also have some far-reaching applications. Let us
mention in this respect two points.

Firstly, recently one of us \cite{Sniady2005GaussuanFluctuationsofYoungdiagrams} developed a
quite general theory for fluctuations of characters and shapes of random Young diagrams
contributing to many natural representations of symmetric groups. The results presented there
are closely (though, not explicitly) related to combinatorics of higher order cumulants. This
connection will be studied in detail in the part IV of this series \cite{HigherOrderFreeness4}
where we prove that under some mild technical conditions Jucys-Murphy elements, which arise
naturally in the study of symmetric groups, are examples of free random
variables of higher order.

In another direction, the description of subfactors in von Neumann
algebras via planar algebras
\cite{Jonesplanar}
relies very much on the
notions of annular non-crossing partitions and thus resembles the
combinatorial objects lying at the basis of our theory of second
order freeness. This indicates that our results could have some
relevance for subfactors.

\subsection*{Overview of the article}

In Section \ref{section2} we will
give a compact survey of
the connection between Voiculescu's free probability theory and random matrix
theory, provide the main motivation, notions and concepts for our extension
of this theory to the level of fluctuations (second order), as well as the statement of
our main results concerning second order cumulants and $R$-transforms.
We will also make a few general remarks about higher order freeness.

In Section \ref{section3} we will introduce the basic notions and relevant results on
permutations, partitions, classical cumulants, Haar unitary random 
matrices, and the Weingarten function.

In Section \ref{sec:correlation} we study the correlation functions (classical cumulants of
traces) of random matrix models. We will see how those are related to cumulants of entries of
the matrices for unitarily invariant random matrices and we will in particular look on the
correlation functions for products of two independent ensembles of random matrices, one of
which is unitarily invariant. The limit of those formulas if the size $N$ of the matrices goes
to infinity will be the essence of what we are going to call ``higher order freeness". Also our
main combinatorial objects, ``partitioned permutations", will arise very naturally in these
calculations.

In Section \ref{sec:multiplicative} we will forget for a while random variables and just look
on the combinatorial essence of our formulas, thus dealing with multiplicative functions on
partitioned permutations and their convolution. The Zeta and M\"obius functions on partitioned
permutations will play an important role in these considerations.

In Section \ref{sec:R-transform} we will derive, for the case of second order, the analogue of
the R-transform formulas.

In Section \ref{sec:higher} we will finally come back to a (non-commutative) probabilistic
context, give the definition and work out the basic properties of ``higher order freeness".

In Section \ref{sec:RMT-IZ} we introduce the notion of ``asymptotic higher order freeness".
We show that the Itzykson-Zuber integral encodes all information
about higher order freeness. We also indicate how our techniques can give some insight into
the computation of some limits of matrix integrals.

In an appendix, Section \ref{sec:surfaced}, we provide a graphical interpretation of
partitioned permutations as a special case of ``surfaced permutations".

\section{Motivation and Statement of our Main Results Concerning Second Order
Freeness and Cumulants}\label{section2} In this section we will first recall in a quite compact
form the main connection between Voiculescu's free probability theory and questions about
random matrices. Then we want to motivate our notion of second order freeness by extending
these questions from the level of expectations to the level of fluctuations. We will recall the
relevant results from the papers
\cite{HigherOrderFreeness1,HigherOrderFreeness2} and state the main new results of the present
paper. Even though in the later parts of the paper our treatment will include freeness of
arbitrarily high order, we restrict ourselves in this section mainly to the second order.
The reason for this is that (apart from first order) second order freeness seems to be the most
important order for applications, so that it seems worthwhile to spell out our general results for
this case more explicitly. Furthermore, it is only there that we have an analogue of $R$-transform
formulas. We will make a few general remarks about higher order freeness at the end of this
section.

\subsection{Moments of random matrices and asymptotic freeness}

Assume we know the eigenvalue distribution of two matrices $A$ and $B$. What can we say about
the eigenvalue distribution of the sum $A+B$ of the matrices? Of course, the latter is not just
determined by the eigenvalues of $A$ and the eigenvalues of $B$, but also by the relation
between the eigenspaces of $A$ and of $B$. Actually, it is a quite hard problem (Horn's
conjecture) --- which was only solved recently --- to
characterize all possible eigenvalue distributions of $A+B$.
However, if one is asking this question in the context of
$N\times N$-random matrices, then in many situations the answer
becomes deterministic in the limit
$N\to\infty$.

\begin{definition}
Let $A=(A_N)_{N\in\NN}$ be a sequence of $N\times N$-random matrices. We say that $A$ has a
limit eigenvalue distribution if the limit of all moments
$$\alpha_n:=\lim_{N\to\infty} E[\tr(A_N^n)]\qquad (n\in\NN)$$
exists, where $E$ denotes the expectation and $\tr$ the normalized trace.
\end{definition}

In this language, our question becomes: Given two random matrix ensembles of $N\times N$-random
matrices, $A=(A_N)_{N\in\NN}$ and $B=(B_N)_{N\in\NN}$, with limit eigenvalue distribution, does
also their sum $C=(C_N)_{N\in\NN}$, with $C_N=A_N+B_N$, have a limit eigenvalue distribution,
and furthermore, can we calculate the limit moments $\alpha_n^C$ of $C$ out of the limit
moments $(\alpha_k^A)_{k\geq 1}$ of $A$ and the limit moments $(\alpha_k^B)_{k\geq 1}$ of $B$
in a deterministic way. It turns out that this is the case if the two ensembles are
in generic position, and then the rule for calculating the limit moments of $C$ are given
by Voiculescu's concept of ``freeness". Let us recall this fundamental result of Voiculescu.

\begin{theorem}[Voiculescu \cite{Voiculescu1991}]
Let $A$ and $B$ be two random matrix ensembles of $N\times N$-random matrices,
$A=(A_N)_{N\in\NN}$ and $B=(B_N)_{N\in\NN}$, each of them with a limit eigenvalue distribution.
Assume that $A$ and $B$ are independent (i.e., for each $N\in\NN$, all entries of $A_N$ are
independent from all entries of $B_N$), and that at least one of them is unitarily invariant
(i.e., for each $N$, the joint distribution of the entries does not change if we conjugate the random matrix
with an arbitrary unitary $N\times N$ matrix).
Then $A$ and $B$ are asymptotically free in the sense of the following definition.
\end{theorem}

\begin{definition}[Voiculescu \cite{Voiculescu-first}]\label{def:freeness-first}
Two random matrix ensembles $A=(A_N)_{N\in\NN}$ and $B=(B_N)_{N\in\NN}$ with limit eigenvalue
distributions are \emph{asymptotically free} if we have for all $p\geq 1$ and all
$n(1),m(1),\dots,n(p)$, $m(p)\geq 1$ that
\begin{multline*}
\lim_{N\to\infty} E\Bigl[\tr \bigl\{ (A_N^{n(1)}-\alpha_{n(1)}^A 1)
\cdot (B_N^{m(1)}-\alpha_{m(1)}^B 1)\cdots \\ 
\cdots
(A^{n(p)}-\alpha_{n(p)}^A 1)\cdot (B^{m(p)}-\alpha_{m(p)}^B 1) \bigr\}\Bigr]=0
\end{multline*}
\end{definition}

One should realize that asymptotic freeness is actually a rule which allows to calculate all
mixed moments in $A$ and $B$, i.e. all expressions
$$\lim_{N\to\infty} E[\tr(A^{n(1)}B^{m(1)}A^{n(2)}B^{m(2)}\cdots A^{n(p)}B^{m(p)})]$$
out of the limit moments of $A$ and the limit moments of $B$. In particular, this means that
all limit moments of $A+B$ (which are sums of mixed moments) exist, thus $A+B$ has a limit
distribution, and are actually determined in terms of the limit moments of $A$ and the limit
moments of $B$. The actual calculation rule is not directly clear from the above definition but
a basic result of Voiculescu shows how this can be achieved by going over from the moments
$\alpha_n$ to new quantities $\kk_n$. In \cite{Speicher1994}, the combinatorial structure behind
these $\kk_n$ was revealed and the name ``free cumulants" was coined for them. Whereas in
the later parts of this paper we will have to rely crucially on the combinatorial description
and their extensions to higher orders, as well as on the definition of more general ``mixed"
cumulants, we will here state the results in the simplest possible form in terms of
generating power series, which
avoids the use of combinatorial objects.

\begin{definition}[Voiculescu \cite{Voiculescu-R}, Speicher \cite{Speicher1994}]
Given the moments $(\alpha_n)_{n\geq 1}$ of some distribution (or limit moments of some
random matrix ensemble), we define
the corresponding \emph{free
cumulants}
$(\kk_n)_{n\geq 1}$ by the following relation between their generating power series: If we put
$$M(x):=1+\sum_{n\geq 1} \alpha_n x^n\qquad\text{and}\qquad
C(x):=1+\sum_{n\geq 1} \kk_n x^n,$$ then we require as a relation between these formal power
series that
$$C(xM(x))=M(x).$$
\end{definition}

Voiculescu actually formulated the relation above in a slightly different way using the
so-called
$R$-transform
$\cR(x)$, which is related to $C(x)$ by the relation
$$C(x)=1+z\cR(x)$$
and in terms of the Cauchy transform $G(x)$ corresponding to a measure with moments $\alpha_n$,
which is related to $M(x)$ by
$$G(x)=\frac { M(\frac 1x)}x.$$ In these terms
the equation $C(xM(x))=M(x)$ says that
\begin{equation}\label{eq:r-transform-first}
\frac 1{G(x)}+\cR(G(x))=x,
\end{equation}
i.e., that $G(x)$ and $K(x):=\frac 1x+\cR(x)$ are inverses of each other under composition.

One should also note that the relation $C(xM(x))=M(x)$ determines the moments uniquely in terms of
the cumulants and the other way around. The relevance of the $\kk_n$ and
the $R$-transform for our problem comes from the following result of Voiculescu,
which provides, together with (\ref{eq:r-transform-first}), a very efficient way
for calculating eigenvalue distributions of the sum of asymptotically free random matrices.

\begin{theorem}[Voiculescu \cite{Voiculescu-R}]
Let $A$ and $B$ be two random matrix ensembles which are asymptotically free. Denote by
$\kk_n^{A}$, $\kk_n^{B}$, $\kk_n^{A+B}$ the free cumulants of $A$, $B$, $A+B$, respectively.
Then one has for all $n\geq 1$ that
$$\kk_n^{A+B}=\kk_n^{A}+\kk_n^{B}.$$
Alternatively,
$$\cR^{A+B}(x)=\cR^A(x)+\cR^B(x).$$
\end{theorem}

This theorem is one reason for calling the $\kk_n$ cumulants, but there is also another
justification for this, namely they are also the limit of classical cumulants of the entries of
our random matrix, in the case that this is unitarily invariant. This description will follow
from our formulas \eqref{eq:KK-k} and \eqref{eq:KK-limit}. We denote the classical cumulants by
$\cc_n$, considered as multi-linear functionals in $n$ arguments.

\begin{theorem}\label{firstlimitcumulant}
Let $A=(A_N)_{N\in\NN}$ be a unitarily invariant random matrix ensemble
of $N\times N$ random matrices $A_N$ whose limit eigenvalue
distribution exists. Then the free cumulants of this matrix ensemble can also be expressed as
the limit of special classical cumulants of the entries of the random matrices: If
$A_N=(a_{ij}^{(N)})_{i,j=1}^N$, then
$$\kk_n^A=\lim_{N\to\infty}N^{n-1} \cc_n(a^{(N)}_{i(1)i(2)},a^{(N)}_{i(2)i(3)},\dots,
a^{(N)}_{i(n),i(1)})$$ for any choice of distinct $i(1),\dots,i(n)$.
\end{theorem}

\subsection{Fluctuations of random matrices and asymptotic second order
freeness} There are many more refined questions about the limiting eigenvalue distribution of
random matrices. In particular, questions around fluctuations have received a lot of interest
in the last decade or so. The main motivation for introducing the concept of ``second order
freeness" was to understand the global fluctuations of the eigenvalues, which means that we
look at the probabilistic behavior of traces of powers of our matrices. The limiting eigenvalue
distribution, as considered in the last section, gives us the limit of the average of this
traces. However, one can make more refined statements about their distributions. Consider a
random matrix $A=(A_N)_{N\in\NN}$ and look on the normalized traces $\tr(A_N^k)$. Our
assumption of a limit eigenvalue distribution means that the limits
$\alpha_k:=\lim_{N\to\infty} E[\tr(A_N^k)]$ exist. It turned out that in many cases the
fluctuation around this limit,
$$\tr(A_N^k)-\alpha_k$$
is asymptotically Gaussian of order $1/N$; i.e., the random variable
$$N\cdot( \tr(A_N^k)-\alpha_k)=\Tr(A_N^k)-N\alpha_k=\Tr(A_N^k-\alpha_k1)$$
(where $\Tr$ denotes the unnormalized trace)
converges for $N\to\infty$ to a normal variable. Actually, the whole family
of centered unnormalized traces $(\Tr(A^k_N)-N\alpha_k)_{k\geq 1}$ converges to a centered
Gaussian family.
(One should note that we restrict all our considerations to complex random matrices;
in the case of real random matrices there are additional complications, which will
be addressed in some future investigations.) Thus the main information about fluctuations of our considered ensemble is
contained in the covariance matrix of the limiting Gaussian family, i.e., in the
quantities
$$\alpha_{m,n}:=\lim_{N\to\infty} \cov(\Tr(A_N^m),\Tr(A_N^n)).$$
Let us emphasize that the $\alpha_n$ and the $\alpha_{m,n}$ are actually limits of classical
cumulants of traces; for the first and second order, with expectation as first and variance
as second cumulant, this
might not be so visible, but it will become evident when we go over to higher orders.
Nevertheless, the $\alpha$'s will
behave and will also be treated like moments; accordingly we will call
the $\alpha_{m,n}$ `fluctuation moments'. We will later define some other quantities
$\kk_{m,n}$, which take the role of cumulants in this context.

This kind of convergence to a Gaussian family was formalized in \cite{HigherOrderFreeness1} as
follows. Note that convergence to Gaussian means that all higher order classical cumulants
converge to zero. As before, we denote the classical cumulants by $\cc_n$; so $\cc_1$ is just
the expectation, and $\cc_2$ the covariance.

\begin{definition}
Let $A=(A_N)_{N\in\NN}$ be an ensemble of $N\times N$ random matrices $A_N$.
We say that it has a \emph{second order
limit distribution} if for all $m,n\geq 1$ the limits
$$\alpha_n:=\lim_{N\to\infty} \cc_1(\tr(A_N^n))$$
and
$$\alpha_{m,n}:=\lim_{N\to\infty} \cc_2(\Tr(A_N^m),\Tr(A_N^n))$$
exist and if
$$\lim_{N\to\infty}\cc_r\bigl(\Tr(A_N^{n(1)}),\dots,\Tr(A_N^{n(r)})\bigr)=0$$
for all $r\geq 3$ and all $n(1),\dots,n(r)\geq 1$.
\end{definition}

We can now ask the same kind of question for the limit fluctuations as for the limit moments;
namely, if we have two random matrix ensembles $A$ and $B$ and we know the second order limit
distribution of $A$ and the second order limit distribution of $B$, does this imply that we have a second order limit distribution
for $A+B$, and, if so, is there an effective way for calculating it. Again, we can only hope
for a positive solution to this if $A$ and $B$ are in a kind of generic position. As it turned
out, the same requirements as before are sufficient for this. The rule for calculating mixed
fluctuations constitutes the essence of the definition of the concept of second order freeness.

\begin{theorem}[Mingo, \'Sniady, Speicher
\cite{HigherOrderFreeness2}] Let $A$ and $B$ be two random
matrix ensembles of $N\times N$-random matrices,
$A=(A_N)_{N\in\NN}$ and $B=(B_N)_{N\in\NN}$, each of them having a second order limit
distribution. Assume that $A$ and $B$ are independent and that at least one of them is unitarily
invariant. Then $A$ and $B$ are asymptotically free of second order in the sense of the
following definition.
\end{theorem}

\begin{definition}
[Mingo, Speicher \cite{HigherOrderFreeness1}]\label{def:freeness-second}
Consider two random
matrix ensembles $A=(A_N)_{N\in\NN}$ and $B=(B_N)_{N\in\NN}$, each of them with a second order
limit distribution. Denote by
$$Y_N\bigl(n(1),m(1),\dots,n(p),m(p)\bigr)$$
the random variable
$$\Tr\bigl((A_N^{n(1)}-\alpha_{n(1)}^A1)(B_N^{m(1)}-\alpha_{m(1)}^B1)\cdots
(A_N^{n(p)}-\alpha^A_{n(p)}1)(B_N^{m(p)}-\alpha_{m(p)}^B1)\bigr).$$
The random matrices $A=(A_N)_{N\in\NN}$ and $B=(B_N)_{N\in\NN}$ are
\emph{asymptotically free of second order} if for all $n,m\geq 1$
$$\lim_{N\to\infty}\cc_2\bigl(\Tr(A_N^n-\alpha^A_n1),\Tr(B_N^m-\alpha^B_m1)\bigr)=0$$
and for all $p,q\geq 1$ and $n(1),\dots,n(p)$,$m(1),\dots,m(p)$,$\tilde n(1),\dots,\tilde
n(q)$, $\tilde m(1),\dots,\tilde m(q)\geq 1$ we have
$$
\lim_{N\to\infty}\cc_2\Bigl(Y_N\bigl(n(1),m(1),\dots,n(p),m(p)\bigr),
Y_N\bigl(\tilde n(1),\tilde m(2),\dots,\tilde n(q),
\tilde m(q)\bigr)\Bigr)=0$$
if $p\not=q$, and otherwise (where we count modulo $p$ for the arguments of the
indices, i.e., $n(i+p)=n(i)$)
\begin{multline*}
\lim_{N\to\infty}\cc_2\Bigl(Y_N\bigl(n(1),m(1),\dots,n(p),m(p)\bigr),Y_N\bigl(
\tilde n(p),\tilde m(p),\dots,\tilde n(1),
\tilde m(1)\bigr)\Bigr)\\=
\sum_{k=1}^p \prod_{i=1}^p\bigl(\alpha^A_{n(i+k)
+\tilde n(i)}-\alpha^A_{n(i+k)}\alpha^A_{\tilde n(i)}\bigr)
\bigl(\alpha^B_{m(i+k)+\tilde m(i+1)}-\alpha^B_{m(i+k)}\alpha^B_{\tilde m(i+1)}\bigr).
\end{multline*}
 \end{definition}

Again, it is crucial to realize that this definition allows one (albeit in a complicated way) to
express every second order mixed moment, i.e., a limit of the form
$$\lim_{N\to\infty} \cc_2\bigl(\Tr(A_N^{n(1)}B_N^{m(1)}\cdots A_N^{n(p)}B_N^{m(p)}),
\Tr (A_N^{\tilde n(1)}
B_N^{\tilde m(1)}\cdots A_N^{\tilde n(q)}B_N^{\tilde m(q)})\bigr)$$
in terms of the second order limits of $A$ and the second order limits of $B$. In particular,
asymptotic freeness of second order also implies that the sum $A+B$ of our random matrix
ensembles has a second order limit distribution and allows one to express them in principle in
terms of the second order limit distribution of $A$ and the second order limit
distribution of $B$. As in the case of first
order freeness, it is not clear at
all how this calculation of the fluctuations of $A+B$ out of the fluctuations of $A$ and the
fluctuations of $B$ can be performed effectively.
It is one of the main results of the present paper to achieve such an effective description.
We are able to solve this problem by providing a second order
cumulant machinery, similar to the first order case.  Again, the idea is to go over to
quantities which behave like cumulants in this setting. The actual description of those relies
on combinatorial objects (annular non-crossing permutations),
but as before this can be reformulated in terms of formal power
series. Let us spell out the definition here in this form. (That this is equivalent to
our actual definition of the cumulants will follow from Theorem \ref{second-R}.)

\begin{definition}\label{def:cumulants-second}
Let $(\alpha_n)_{n\geq 1}$ and $(\alpha_{m,n})_{m,n\geq 1}$ describe the first and second
order limit moments of a random matrix ensemble.
We
define the corresponding \emph{first and second order free cumulants}
$(\kk_n)_{n\geq 1}$ and $(\kk_{m,n})_{m,n\geq 1}$
by the
following requirement in terms of the corresponding generating power series. Put
$$C (x):=1+\sum_{n\geq 1} \kk_n x^n,\qquad
C(x,y):=\sum_{m,n\geq 1}\kk_{m,n}x^my^n$$ and
$$M(x):=1+\sum_{n\geq 1} \alpha_nx^n,\qquad
M(x,y):=\sum_{m,n\geq 1}\alpha_{m,n}x^my^n.$$
Then we require as relations between these
formal power series that
\begin{equation}
C(xM(x))=M(x)
\end{equation}
and for the second order
\begin{equation}
M(x,y)=H\bigl(xM(x),yM(y)\bigr)\cdot \frac{\frac d{dx}(xM(x))}{M(x)}\cdot \frac{
\frac d{dy}(yM(y))}{M(y)},
\end{equation}
where
\begin{equation}
H(x,y):=C(x,y)-xy\frac{\partial^2}{\partial x\partial y}\log\Bigl(\frac{xC(y)-yC(x)}{x-y}\Bigr),
\end{equation}
or equivalently,
\begin{multline}\label{eq:M-R-motivation}
M(x,y)=C\bigl(xM(x),yM(y)\bigr)\cdot \frac{\frac d{dx}(xM(x))}{M(x)}\cdot \frac{
\frac d{dy}(yM(y))}{M(y)}\\
+xy\Bigl(
\frac{\frac d{dx}(xM(x))\cdot \frac d{dy}(yM(y))}{(xM(x)-yM(y))^2}-\frac 1{(x-y)^2}
\Bigr).\end{multline}
\end{definition}

As in the first order case, instead of the moment power series $M(x,y)$ one
can consider a kind of second order Cauchy transform, defined by
$$G(x,y):=\frac {M(\frac 1x,\frac 1y)}{xy}.$$
If we also define a kind of second order $R$ transform $\cR(x,y)$ by
$$\cR(x,y):=\frac 1{xy}C(x,y),$$
then the formula (\ref{eq:M-R-motivation}) takes on a particularly nice form:
\begin{equation}\label{eq:G-R-motivation}
G(x,y)=G'(x)G'(y)\Bigl\{
\cR(G(x),G(y))+\frac 1{(G(x)-G(y))^2}\Bigl\}-\frac 1{(x-y)^2}.
\end{equation}
$G(x)$ is here, as before, the first order Cauchy transform,
$G(x)=\frac 1x{M(1/x)}$.

The $\kk_{m,n}$ defined above deserve the name ``cumulants" as they linearize the problem of
adding random matrices which are asymptotically free of second order. Namely,
as will follow from our Theorem \ref{thm:additivity-of-cumulants}, we have the
following theorem, which provides, together with (\ref{eq:G-R-motivation}), an effective
machinery for calculating the fluctuations of the sum of asymptotically free random matrices.

\begin{theorem}
Let $A$ and $B$ be two random matrix ensembles which are asymptotically free.
Then one has for
all $m,n\geq 1$ that
$$\kk_n^{A+B}=\kk_n^{A}+\kk_n^{B}\qquad\text{and}\qquad
\kk_{m,n}^{A+B}=\kk_{m,n}^A+\kk_{m,n}^B.$$
Alternatively,
$$\cR^{A+B}(x)=\cR^A(x)+\cR^B(x)$$
and
$$\cR^{A+B}(x,y)=\cR^A(x,y)+\cR^B(x,y).$$
\end{theorem}

Again, one can express the second order cumulants as limits of classical cumulants of entries
of a unitarily invariant matrix. In contrast to the first order case, we have now to run over
two disjoint cycles in the indices of the matrix entries. This theorem will follow from
our formulas \eqref{eq:KK-k} and \eqref{eq:KK-limit}.

\begin{theorem}\label{secondlimitcumulant}
Let $A=(A_N)_{N\in\NN}$ be a unitarily invariant random matrix ensemble which has a second
order limit distribution. Then the second order free cumulants of this matrix ensemble can also
be expressed as the limit of classical cumulants of the entries of the random matrices:
If $A_N=(a_{ij}^{(N)})_{i,j=1}^N$, then
\begin{multline*}
$$\kk_{m,n}^A=\lim_{N\to\infty}N^{m+n} \cc_{m+n}(a^{(N)}_{i(1)i(2)},a^{(N)}_{i(2)i(3)},\dots,
a^{(N)}_{i(m),i(1)},\\a^{(N)}_{j(1)j(2)},a^{(N)}_{j(2)j(3)},\dots, a^{(N)}_{j(n),j(1)})
\end{multline*}
for any
choice of distinct $i(1),\dots,i(m),j(1),\dots,j(n)$.
\end{theorem}

This latter theorem makes it quite obvious that the second order cumulants for Gaussian
as well as for Wishart matrices vanish identically, i.e., $\cR(x,y)=0$ and thus
we obtain in these cases that the second order Cauchy transform
is totally determined in terms of the first order Cauchy transform (i.e., in terms of
the limiting eigenvalue distribution) via
\begin{equation}
G(x,y)=\frac {G'(x)G'(y)}{(G(x)-G(y))^2}-\frac 1{(x-y)^2}.
\end{equation}
This formula for fluctuations of Wishart matrices
was also derived by Bai and Silverstein in \cite{BS}.

\subsection{Higher order freeness}
The idea for higher order freeness is the same as for second order one.
For a random matrix ensemble $A=(A_N)_{N\in\NN}$ we define $r$-th order limit moments as
the scaled limit of classical cumulants of $r$ traces of powers of our matrices,
$$\alpha_{n_1,\dots,n_r}:=
\lim_{N\to\infty}N^{2-r}\cc_r\bigl(\Tr(A_N^{n(1)}),\dots,\Tr(A_N^{n(r)})\bigr).$$ (The choice
of $N^{2-r}$ is motivated by the fact that this is the leading order for many interesting
random matrix ensembles, e.g. Gaussian or Wishart. Thus our theory of higher order freeness
captures the features of random matrix ensembles whose cumulants of traces decay in the same
way as for Gaussian random matrices.) Then we look at two random matrix ensembles $A$ and $B$
which are independent, and one of them unitarily invariant. The mixed moments in $A$ and $B$ of
order $r$ are, in leading order in the limit $N\to\infty$, determined by the limit moments of
$A$ up to order $r$ and the limit moments of $B$ up to order $r$.  The structure of these
formulas motivates directly the definition of cumulants of the considered order. The definition
of those is in terms of a moment-cumulant formula, which gives a moment in terms of cumulants
by summing over special combinatorial objects, which we call ``partitioned permutations". Most
of the theory we develop relies on an in depth analysis of properties of these partitioned
permutations and the corresponding convolution of multiplicative functions on partitioned
permutations. Our definition of ``higher order freeness" is then in terms of the vanishing of
mixed cumulants. It follows quite easily that in the first and second order case this gives the
same as the relations in Definitions \ref{def:freeness-first} and \ref{def:freeness-second},
respectively. For higher orders, however, we are not able to find an explicit relation of that
type.

This reflects somehow the observation that our general formulas
in terms of sums over partitioned permutations are the same for
all orders, but that evaluating or simplifying these sums (by
doing partial summations) is beyond our abilities for orders
greater than 2. Reformulating the combinatorial
relation between moments and cumulants in terms of generating
power series is one prominent example for this. Whereas this is
quite easy for first order, the complexity of the arguments and
the solution (given in Definition \ref{def:cumulants-second}) is
much higher for second order, and out of reach for higher order.

One should note that an effective (analytic or symbolic)
calculation of higher order moments of a sum $A+B$
for $A$ and $B$ free of higher order relies usually on the presence
of such generating power series formulas. In this sense, we have succeeded in
providing an effective machinery for dealing with fluctuations (second order), but we
were not able to do so for higher order.

Our results for higher orders are more of a theoretical nature. One of the main problems we
have to address there is the associativity of the notion of higher order freeness. Namely, in
order to be an interesting concept, our definition that $A$ and $B$ are free of higher order
should of course imply that any function of $A$ is also free of higher order from any function
of $B$. Whereas for first and second order this follows quite easily from the equivalent
characterization of freeness in terms of moments as in Definitions \ref{def:freeness-first} and
\ref{def:freeness-second}, the absence of such a characterization for higher orders makes this
a more complicated matter. Namely, what we have to see is that the vanishing of mixed cumulants
in random variables implies also the vanishing of mixed cumulants in elements from the
generated algebras. This is a quite non-trivial fact and requires a careful analysis, see
section \ref{sec:higher}.

\section{Preliminaries}\label{section3}

\subsection{Some general notation}
For natural numbers $m,n\in\NN$ with $m<n$, we denote by $[m,n]$
the interval of natural numbers between $m$ and $n$, i.e.,
$$[m,n]:=\{m,m+1,m+2,\dots,n-1,n\}.$$
For a matrix $A=(a_{ij})_{i,j=1}^N$, we denote by $\Tr$ the
unnormalized and by $\tr$ the normalized trace,
$$\Tr(A):=\sum_{i=1}^N a_{ii},\qquad \tr(A):=\frac 1N \Tr(A).$$

\subsection{Permutations}

We will denote the set of permutations on $n$ elements by
$S_n$.  We will quite often use the cycle notation for such
permutations, i.e., $\pi=(i_1,i_2,\dots,i_r)$ is a cycle
which sends $i_k$ to $i_{k+1}$ ($k=1,\dots,r$), where
$i_{r+1}=i_1$.

\subsubsection{Length function}

For a permutation $\pi\in S_n$ we denote by $\#\pi$ the number of cycles of $\pi$ and by
$\vert\pi\vert$ the minimal number of transpositions needed to write $\pi$ as a product of
transpositions. Note that one has
$$\vert\pi\vert+\#\pi=n\qquad\text{for all $\pi\in S_n$}.$$

\subsubsection{Non-crossing permutations}

Let us denote by $\gamma_n\in S_n$ the cycle
$$\gamma_n=(1, 2 ,\dots, n).$$
For all $\pi\in S_n$ one has that
$$\vert\pi\vert+\vert\gamma_n\pi^{-1}\vert\leq n-1.$$
If we have equality then we call $\pi$ \emph{non-crossing}. Note
that this is equivalent to
$$\#\pi+\#(\gamma_n\pi^{-1})=n+1.$$ If $\pi$ is
non-crossing, then so are $\gamma_n\pi^{-1}$ and
$\pi^{-1}\gamma_n$; the latter is called the
\emph{(Kreweras) complement} of $\pi$.

We will denote the set of non-crossing permutations in $S_n$
by $NC(n)$. Note that such a non-crossing permutation can be
identified with a non-crossing partition, by forgetting the
order on the cycles. There is exactly one cyclic order on
the blocks of a non-crossing partition which makes it into a
non-crossing permutation.

\subsubsection{Annular non-crossing permutations}

Fix $m,n\in\NN$ and denote by $\gamma_{m,n}$ the product of the
two cycles
$$\gamma_{m,n}=(1, 2,\dots, m)(m+1, m+2,\dots, m+n).$$
More generally, we shall denote by $\gamma_{m_1,\dots,m_k}$ the
product of the corresponding $k$ cycles.

We call a $\pi\in S_{m+n}$ \emph{connected} if the pair
$\pi$ and $\gamma_{m,n}$ generates a transitive subgroup in
$S_{m+n}$. A connected permutation $\pi\in S_{m+n}$ always
satisfies
\begin{equation}\label{annular-NC}
\vert\pi\vert+\vert\gamma_{m,n}\pi^{-1}\vert\leq m+n.
\end{equation}
If $\pi$ is connected and if we have equality in that
equation then we call $\pi$ \emph{annular
non-crossing}. Note that if $\pi$ is annular non-crossing then
$\gamma_{m,n}\pi^{-1}$ is also annular non-crossing. Again, we
call the latter the \emph{complement} of $\pi$. Of course,
all the above notations depend on the pair $(m,n)$; if we
want to emphasize this dependency we will also speak about
$(m,n)$-connected permutations and $(m,n)$-annular
non-crossing permutations.

We will denote the set of $(m,n)$-annular non-crossing permutations by $\SNC(m,n)$. A cycle of
a $\pi\in \SNC(m,n)$ is called a \emph{through-cycle} if it contains points on both cycles.
Each $\pi\in\SNC(m,n)$ is connected and must thus have at least one through-cycle. The subset
of $\SNC(m,n)$ where all cycles are through-cycles will be denoted by $\SNC^{all}(m,n)$.

Again one can go over from annular non-crossing permutations to annular non-crossing partitions
by forgetting the cyclic orders on cycles; however, in the annular case, the relation between
non-crossing permutation and non-crossing partition is not one-to-one. Since we will not use
the language of annular partitions in the present paper, this is of no relevance here.

Annular non-crossing permutations and partitions were
introduced in
\cite{MingoNica2004annular}; there, many different
characterizations---in particular, the one
(\ref{annular-NC}) above in terms of the length
function---were given.

\subsection{Partitions}

We say that $\cV=\{V_1,\dots,V_k\}$ is a partition of a set $[1,n]$ if the sets $V_i$ are
disjoint and non--empty and their union is equal to $[1,n]$. We call $V_1,\dots,V_k$ the blocks
of partition $\cV$.

If $\cV=\{V_1,\dots,V_k\}$ and $\cW=\{W_1,\dots,W_l\}$ are partitions of the same set, we say
that $\cV\leq \cW$ if for every block $V_i$ there exists some block $W_j$ such that
$V_i\subseteq W_j$. For a pair of partitions $\cV,\cW$ we denote by $\cV\vee \cW$ the smallest
partition $\cU$ such that $\cV\leq \cU$ and $\cW\leq \cU$. We denote by $1_n=\big\{ [1,n]
\big\}$ the biggest partition of the set $[1,n]$.

If $\pi\in S_n$ is a permutation, then we can associate to $\pi$ in a natural way
a partition whose blocks consist exactly of the cycles of $\pi$; we will denote this
partition either by $0_\pi\in\cP(n)$ or, if the context makes the meaning clear,
just by $\pi\in\cP(n)$.

For a permutation $\pi\in S_n$ we say that a partition $\cV$ is $\pi$-invariant if $\pi$
preserves each block of $\cV$. This means that $0_\pi\leq \cV$ (which we
will usually write just as $\pi\leq \cV$).

If $\cV=\{V_1,\dots,V_k\}$ is a partition of the set $[1,n]$ and if, for $1\leq i\leq k$,
$\pi_i$ is a permutation of the set $V_i$ we denote by $\pi_1\times \cdots\times \pi_k\in S_n$
the concatenation of these permutations. We say that $\pi=\pi_1\times\cdots\times \pi_k$ is a
cycle decomposition if additionally every factor $\pi_i$ is a cycle.

\subsection{Classical cumulants}
Given some classical probability space $(\Omega,\ab P)$ we 
denote by
$\EE$ the expectation with respect to the corresponding probability measure,
$$\EE(a):=\int_\Omega a(\omega)dP(\omega)$$
and by $L^{\infty-}(\Omega,P)$ the algebra of random variables for which all moments
exist. 
Let us for the following put $\cA:=L^{\infty-}(\Omega,P)$.

We extend the linear functional $\EE:\cA\to\CC$ 
to a corresponding multiplicative functional on all
partitions by ($\cV\in \cP(n)$, $a_1,\dots,a_n\in\cA$)
\begin{equation}
\EE_\cV[a_1,\dots,a_n]:=\prod_{V\in\cV} \EE[a_1,\dots,a_n\vert_V],
\end{equation}
where we use the notation
$$\EE[a_1,\dots,a_n\vert_V]:=\EE(a_{i_1}\cdots a_{i_s})\qquad\text{for}\qquad
V=(i_1<\dots<i_s)\in\cV.$$

Then, for $\cV\in \cP(n)$, we define the \emph{classical cumulants} $\cc_\cV$ as multilinear
functionals on $\cA$ by
\begin{equation}
\cc_\cV[a_1,\dots,a_n]=\sum_{\cW\in \cP(n)\atop \cW\leq \cV} \EE_\cW[a_1,\dots,a_n]\cdot
\moeb_{\cP(n)}(\cW,\cV),
\end{equation}
where $\moeb_{\cP(n)}$ denotes the M\"obius function on $\cP(n)$
(see \cite{rota}).

The above definition is, by M\"obius inversion on $\cP(n)$, equivalent to
$$\EE(a_1\cdots a_n)=\sum_{\pi\in\cP(n)}\cc_\pi[a_1,\dots,a_n].$$
The $\cc_\pi$ are also multiplicative with respect to the blocks
of $\cV$ and thus determined by the values of
$$\cc_n(a_1,\dots,a_n):=\cc_{1_n}[a_1,\dots,a_n].$$

Note that we have in particular
$$\cc_1(a)=\EE(a)\qquad\text{and}\qquad
\cc_2(a_1,a_2)=\EE(a_1a_2)-\EE(a_1)\EE(a_2).$$

An important property of classical cumulants is the following formula
of Leonov and Shiryaev \cite{LS} for cumulants with products as arguments.

Let $m,n\in\NN$ and $1 \leq i(1) < i(2) < \cdots < i(m)=n$. Define
$\cU\in\cP(n)$ by
\begin{equation*}
\cU=\bigl\{\bigl(1,\dots,i(1)\bigr),\bigl(i(i)+1,\dots,i(2)\bigr),
\dots, \bigl({i({m-1})+1},\dots,i(m)\bigr)\bigr\}.
\end{equation*}
Consider now random variables $a_1,\dots,a_n\in\cA$ and define
\begin{align*}
A_1:&=a_1 \cdots a_{i(1)}\\
A_2:&=a_{i(1)+1}\cdots a_{i(2)}\\
&\vdots\\
A_m:&= a_{i({m-1})+1} \cdots a_{i(m)}.
\end{align*}
Then we have
\begin{equation}\label{eq:product-formula-classical-cumulants}
  \cc_m(A_1,A_2,\dots,A_m)
            = \sum_{\cV \in \cP(n) \atop \cV \vee \cU = 1_n}
    \cc_{\cV} [a_1,\dots,a_n].
\end{equation}
The sum on the right-hand side is running over those partitions of $n$ elements
which satisfy $\cV\vee\cU=1_n$, which are, informally speaking, those
partitions which connect
all the arguments of the cumulant $\cc_m$, when written in terms of the $a_i$.

Here is an example for this formula; for $\cc_2(a_1a_2,a_3a_4)$.
In order to reduce the number
of involved terms we will restrict to the special case where $\EE(a_i)=0$ (and thus also
$\cc_1(a_i)=0$) for all $i=1,2,3,4$.
There are three partitions $\pi\in\cP(4)$ without singletons which satisfy
$$\pi\vee \{(1,2),(3,4)\}=1_4,$$ namely
$$
\setlength{\unitlength}{0.3cm}
\begin{picture}(10,10)
  \thicklines
  \put(0,3){\line(0,1){2}}
  \put(0,3){\line(1,0){9}}
  \put(6,4){\line(0,1){1}}
  \put(3,4){\line(0,1){1}}
  \put(9,3){\line(0,1){2}}
\put(3,4){\line(1,0){3}} \thicklines
  \put(0,6){\line(0,1){2}}
  \put(0,6){\line(1,0){9}}
  \put(3,6){\line(0,1){2}}
  \put(6,6){\line(0,1){2}}
\put(9,6){\line(0,1){2}}
  \thicklines
  \put(0,1){\line(0,1){1}}
\put(0,1){\line(1,0){6}} \put(6,1){\line(0,1){1}}
  \put(3,0){\line(0,1){2}}
  \put(3,0){\line(1,0){6}}
  \put(9,0){\line(0,1){2}}
 \put(-.25,9.5){\dashbox{0.2}(4.2,1.5){}}
 \put(5.5,9.5){\dashbox{0.2}(4,1.5){}}
\put(0,10){$a_1$} \put(2.5,10){$a_2$} \put(5.75,10){$a_3$} \put(8,10){$a_4$}
\end{picture}
$$
and thus formula (\ref{eq:product-formula-classical-cumulants}) gives in this case
\begin{multline*}
\cc_2(a_1a_2,a_3a_4)=\cc_4(a_1,a_2,a_3,a_4)\\+
\cc_2(a_1,a_4)\cc_2(a_2,a_3)+\cc_2(a_1,a_3)\cc_2(a_2,a_4).
\end{multline*}

As a consequence of (\ref{eq:product-formula-classical-cumulants}) one has the following
important corollary: If $\{a_1,\dots,a_n\}$ and $\{b_1,\dots,b_n\}$ are independent
then
\begin{equation}\label{eq:product-special-case}
\cc_\cW[a_1b_1,\dots,a_nb_n]=\sum_{\cV,\cV'\in\cP(n)\atop
\cV\vee\cV'=\cW}
\cc_{\cV}[a_1,\dots,a_n]\cdot \cc_{\cV'}[b_1,\dots,b_n].
\end{equation}

\subsection{Haar distributed unitary random matrices and
            the Weingarten function}

In the following we will be interested in the asymptotics of
special matrix integrals over the group $\UU(N)$ of unitary
$N\times N$-matrices. We always equip the compact group
$\UU(N)$ with its Haar probability measure.
A random matrix whose distribution is this measure
will be called a \emph{Haar distributed unitary random
matrix}. Thus the expectation $\EE$ over this ensemble is
given by integrating with respect to the Haar measure.

The expectation of products of entries of Haar distributed
unitary random matrices can be described in terms of a
special function on the permutation group. Since such
considerations go back to Weingarten
\cite{Weingarten1976}, Collins
\cite{Collins2002} calls this function the \emph{Weingarten function} and
denotes it by $\Wg$. We will follow his notation. In the
following we just recall the relevant information about this
Weingarten function, for more details we refer to
\cite{Collins2002,CollinsSniady2004}.

We use the following definition of the Weingarten function. For
$\pi\in S_n$ and $N\geq n$ we put
$$\Wg(N,\pi)=
\EE[u_{11}\cdots
u_{nn}\overline{u_{1\pi(1)}}\cdots\overline{u_{n\pi(n)}}],$$
where
$U=(u_{ij})_{i,j=1}^N$ is an $N\times N$ Haar distributed unitary random
matrix. Sometimes we will suppress the dependency on $N$ and just
write $\Wg(\pi)$. This $\Wg(N,\pi)$ depends only on the conjugacy
class of $\pi$. General matrix integrals over the unitary group
can be calculated as follows:
\begin{multline}\label{eq:Wick-unitary}
\EE[u_{i'_1j'_1}\cdots u_{i'_nj'_n}\overline {u_{i_1j_1}}
\cdots\overline{ u_{i_nj_n}}]\\= \sum_{\alpha,\beta\in S_n}
\delta_{i_1i'_{\alpha(1)}}\cdots \delta_{i_ni'_{\alpha(n)}}
\delta_{j_1j'_{\beta(1)}}\cdots\delta_{j_nj'_{\beta(n)}}
\Wg(\beta\alpha^{-1}).
\end{multline}
This formula for the calculation of moments of the entries
of a Haar unitary random matrix bears some resemblance to the
Wick formula for the joint moments of the entries
of Gaussian random matrices; thus we will call
(\ref{eq:Wick-unitary}) the \emph{Wick formula for Haar unitary
matrices}.

The Weingarten function is quite a complicated object, and its
full understanding is at the basis of questions around
Itzykson-Zuber integrals. One knows (see, e.g.,
\cite{Collins2002,CollinsSniady2004}) that the leading order in $1/N$ is given by
$\vert\pi\vert+n$ and increases in steps of 2.

\subsection{Cumulants of the Weingarten function}
We will also need some (classical) relative cumulants of the Weingarten function, which were
introduced in \cite[\S 2.3]{Collins2002}. As before, let $\moeb_{\cP(n)}$ be the M\"obius
function on the partially ordered set of partitions of $[1,n]$ ordered by inclusion.

Let us first extend the Weingarten function by multiplicative extension, for $\cV\geq\pi$, by
$$\Wg(\cV,\pi):=\prod_{V\in\cV}\Wg(\pi|_V),$$
where $\pi|_V$ denotes the restriction of $\pi$ to the block $V\in\cV$ (which is invariant
under $\pi$ since $\pi\leq\cV$).

The \emph{relative cumulant of the Weingarten function} is now, for $\sigma\leq\cV\leq\cW$,
defined by
\begin{equation}\label{eq:cWg-Wg}
\cWg_{\cV,\cW}(\sigma)=
\sum_{\cU\in\cP(n) \atop \cV \leq \cU \leq \cW}
\moeb(\cU,\cW) \cdot \Wg(\cU, \sigma)
\end{equation}
Note that, by M\"obius inversion, this is, for any
$\sigma\leq\cV\leq\cW$, equivalent to
\begin{equation}\label{eq:Wg-cWg}
\Wg(\cW,\sigma)=\sum_{\cU\in\cP(n)\atop \cV\leq\cU\leq\cW}\cWg_{\cV,\cU}(\sigma).
\end{equation}
In \cite[Cor. 2.9]{Collins2002} it was shown that the order of $\cWg_{\cV,\cW}(\sigma)$ is at
most
\begin{equation}\label{eq:max-order-Wg}
N^{-2n+\#\sigma+2\#\cW-2\#\cV}.
\end{equation}

\section{Correlation functions for random matrices}
\label{sec:correlation}

\subsection{Correlation functions and partitioned permutations}
Let us consider $N\times N$-random matrices $B_1,\dots,B_n:\Omega\to M_N(\CC)$.
The
main information we are interested in are the ``correlation
functions" $\ff^{}_n$ of these matrices, given by classical
cumulants of their traces, i.e.,
$$\ff^{}_n(B_1,\dots,B_n):=\cc_n(\Tr(B_1),\dots,\Tr(B_n)).$$
Even though these correlation functions are cumulants, it is more adequate
to consider them as a kind of moments for our random matrices. Thus, we will
also call them sometimes \emph{correlation moments}.

We will also need to consider traces of products which are best
encoded via permutations. Thus, for $\pi\in S_n$,
$\ff^{}(\pi)[B_1,\dots,B_n]$ shall mean that we take cumulants of
traces of products along the cycles of $\pi$. For an $n$-tuple
$B=(B_1,\dots,B_n)$ of random matrices and a cycle
$c=(i_1,i_2,\dots,i_k)$ with $k\leq n$ we denote
$$B\vert_c:=B_{i_1}B_{i_2}\cdots B_{i_k}.$$
(We do not distinguish between products which differ by a cyclic rotation of the factors;
however, in order to make this definition well-defined we could normalize our cycle
$c=(i_1,i_2,\dots,i_k)$ by the requirement that $i_1$ is the smallest among the appearing
numbers.) For any $\pi\in S(n)$ and any $n$-tuple $B=(B_1,\dots,B_n)$ of random matrices we put
$$\ff(\pi)[B_1,\dots,B_n]:=\ff_r(B\vert_{c_1},\dots,
B\vert_{c_r}),$$ where $\pi$ consists of the cycles $c_1,\dots,c_r$.

Example:
\begin{align*}
\ff^{}((1,3)(2,5,4))[B_1,B_2,B_3,B_4,B_5]&=
\ff^{}_2(B_1B_3,B_2B_5B_4)\\
&=\cc_2(\Tr(B_1B_3),\Tr(B_2B_5B_4))
\end{align*}

Furthermore, we also need to consider more general products of such
$\ff(\pi)$'s. In order to index such products we will use pairs
$(\cV,\pi)$ where $\pi$ is, as above, an element in $S_n$, and
$\cV\in\cP(n)$ is a partition which is compatible with the cycle
structure of $\pi$, i.e., each block of $\cV$ is fixed under $\pi$,
or to put it another way, $\cV\geq\pi$. In the latter inequality we
use the convention that we identify a permutation with the partition
corresponding to its cycles if this identification is obvious from
the structure of the formula; we will write this partition 
$0_\pi$ or just 0 if no confusion will result.

\begin{notation}
A \emph{partitioned permutation} is a pair $(\cV,\pi)$ consisting of
$\pi\in S_n$ and $\cV\in\cP(n)$ with $\cV\geq\pi$. We will denote
the set of partitioned permutations of $n$ elements by $\cPP(n)$. We
will also put
$$\cPP:=\bigcup_{n\in\NN}\cPP(n).$$
\end{notation}

For such a $(\cV,\pi)\in\cPP$ we denote finally
$$\ff^{}(\cV,\pi)[B_1,\dots,B_n]:=\prod_{V\in\cV}\ff(\pi\vert_V)
[B_1,\dots,B_n\vert_V].$$

Example:
\begin{align*}
\ff\bigl(\{1,3,4\}\{2\},(1,3)(2)(4)\bigr)&[B_1,B_2,B_3,B_4]\\&=
\ff_2(B_1B_3,B_4)\cdot \ff_1(B_2)\\
&=\cc_2(\Tr(B_1B_3),\Tr(B_4))\cdot \cc_1(\Tr(B_2))
\end{align*}

Let us denote by $\Tr_\sigma$ as usual a product of traces along the cycles of $\sigma$. Then
we have the relation
$$\EE\{\Tr_\sigma[A_1,\dots,A_n]\}=
\sum_{\cW\in\cP(n)\atop \cW\geq\sigma} \ff(\cW,\sigma)[A_1,\dots,A_n].$$

By using the formula (\ref{eq:product-formula-classical-cumulants}) of Leonov and Shiryaev one
sees that in terms of the entries of our matrices $B_k=(b_{ij}^{(k)})_{i,j=1}^N$  our
$\ff(\cU,\gamma)$ can also be written as
\begin{equation}\label{eq:ff-entries}
\ff(\cU,\gamma)[B_1,\dots,B_n]=\sum_{\cV\leq\cU\atop \cV\vee\gamma=\cU} \sum_{i(1),
\dots,i(n)=1}^N
\cc_\cV[b_{i(1)i(\gamma(1))}^{(1)},\dots,
b_{i(n)i(\gamma(n))}^{(n)}].
\end{equation}

\subsection{Moments of unitarily invariant random matrices}

For unitarily invariant random matrices there exists a definite relation between
cumulants of traces and cumulants of entries.
We want to work out this connection
in this section.

\begin{definition}
Random matrices $A_1,\dots,A_n$ are called \emph{unitarily
invariant} if the joint distribution of all their entries
does not change by global conjugation with any unitary
matrix, i.e., if, for any unitary matrix $U$, the
matrix-valued random variables $A_1,\dots, A_n:\Omega\to
M_N(\CC)$ have the same joint distribution as the
matrix-valued random variables $UA_1U^*, \dots, UA_nU^* :
\Omega \to M_N(\CC)$.
\end{definition}

Let $A_1,\dots,A_n$ be unitarily invariant random matrices. We will now
try expressing the microscopic quantities ``cumulants of entries
of the $A_i$" in terms of the macroscopic quantities ``cumulants
of traces of products of the $A_i$".

In order to make this connection we have to use the unitary invariance of our
ensemble. By definition,
this
means that $A_1,\dots,A_n$ has the same distribution as
$\tilde A_1,\dots,\tilde A_n$ where $\tilde A_i:=UA_iU^*$.
Since this holds for any unitary $U$, the same is true after averaging over such $U$, i.e.,
we can take in the definition of the $\tilde A_i$ the $U$ as
Haar distributed
unitary random matrices, independent from $A_1,\dots,A_n$.
This reduces calculations for unitarily invariant ensembles
essentially to properties of Haar unitary random matrices; in particular,
the Wick
formula for the $U$'s implies that we have an analogous Wick formula for joint moments in the
entries of the  $A_i$. Let us write
$A_k=(a_{ij}^{(k)})_{i,j=1}^N$ and $\tilde A_k=(\tilde a_{ij}^{(k)})_{i,j=1}^N$. Then we
can calculate:
\begin{align*}
\EE\bigl\{a_{p_1 r_{1}}^{(1)} \cdots\
         a_{p_{n} r_n}^{(n)}
\bigr\}&=\EE
\bigl\{\tilde a_{p_1 r_{1}}^{(1)} \cdots\
         \tilde a_{p_{n} r_n}^{(n)}
\bigr\}\\
&=\sum_{i,j} \EE\{u_{p_1 i_1}a^{(1)}_{i_1j_1}\ol
{u_{r_1j_1}}\cdots u_{p_ni_n}a^{(n)}_{i_nj_n}\ol {u_{r_nj_n}}\}\\
&=\sum_{i,j}\EE\{u_{p_1 i_1}\ol {u_{r_1j_1}}\cdots
u_{p_ni_n}\ol {u_{r_nj_n}}\}\cdot \EE\{a^{(1)}_{i_1j_1}\cdots
a^{(n)}_{i_nj_n}\}
\\
&= \sum_{i,j}\sum_{\pi,\sigma\in S_n} \delta_{r,p\circ
\pi}\delta_{j,i\circ \sigma}\Wg(\sigma\pi^{-1}) \cdot
\EE\{a^{(1)}_{i_1j_1}\cdots a^{(n)}_{i_nj_n}\}\\
&=\sum_{\pi\in S_n}\delta_{r,p\circ\pi}\cdot \GG(\pi)[A_1,\dots,A_n],
\end{align*}
where
\begin{align}
\label{eq:definitionofGG}
\GG(\pi)[A_1,\dots,A_n]:&= \sum_{\sigma\in S_n}
\Wg(\sigma\pi^{-1}) \cdot\sum_{i}
\EE\{a^{(1)}_{i_1i_{\sigma(1)}}\cdots a^{(n)}_{i_ni_{\sigma(n)}}\}\\
&=\sum_{\sigma\in S_n} \Wg(\sigma\pi^{-1}) \cdot \EE\{\Tr_\sigma[A_1,\dots,A_n]\}.\notag\\
&=\sum_{\sigma\in S_n} \Wg(\sigma\pi^{-1})\cdot\sum_{\cW\in\cP(n)\atop \cW\geq\sigma}
\ff(\cW,\sigma)[A_1,\dots,A_n]\notag\\
&=\sum_{(\cW,\sigma)\in\cPP(n)}\Wg(\sigma\pi^{-1})\cdot \ff(\cW,\sigma)[A_1,\dots,A_n].\notag
\end{align}
The important point here is that $\GG(\pi)[A_1,\dots,A_n]$ depends only on the macroscopic
correlation moments of $A$.

We can extend the above to products of expectations by
\begin{align*}
\EE_\cV [a_{p_1 r_{1}},\dots
        a_{p_{n} r_n}]=\sum_{\pi\in S_n\atop \pi\leq \cV}
\delta_{r,p\circ\pi}\cdot \GG(\cV,\pi)[A_1,\dots,A_n],
\end{align*}
where $\GG(\cV,\pi)$ is given by
multiplicative extension:
\begin{align}
\label{eq:definitionofGGmulti}\GG(\cV,\pi)[A_1,\dots,A_n]:&=\prod_{V\in
\cV}\GG(\pi\vert_V)[A_1,\dots,A_n\vert_V]\notag\\
&=\sum_{(\cW,\sigma)\in\cPP(n)\atop \cW\leq\cV} \Wg(\cV,\sigma\pi^{-1})\cdot
\ff(\cW,\sigma)[A_1,\dots,A_n]
\end{align}

Now we can look on the cumulants of the entries of our unitarily invariant matrices $A_i$; they
are given by
\begin{align*}k_\cV\bigl\{a^{(1)}_{p_1 r_{1}},\dots,a^{(n)}_{p_{n} r_n}\}&=
\sum_{{\cU}\in\cP(n)\atop {\cU}\leq \cV} \moeb_{\cP(n)}({\cU},\cV)\cdot\EE_{\cU} [a^{(1)}_{p_1
r_{1}},\dots
        a^{(n)}_{p_{n} r_n}]\\
&=\sum_{{\cU}\leq \cV}\sum_{\pi\in S_n\atop \pi\leq {\cU}}
\delta_{r,p\circ\pi}\cdot \moeb_{\cP(n)}({\cU},\cV)\cdot\GG(\cU,\pi)[A_1,\dots,A_n]\\
&=\sum_{\pi\in S_n\atop \pi\leq \cV} \delta_{r,p\circ\pi}\sum_{\cU \in\cP(n)\atop \cV\geq
{\cU}\geq\pi}\moeb_{\cP(n)}({\cU},\cV)\cdot \GG(\cU,\pi)[A_1,\dots,A_n].
\end{align*}
With the definition
\begin{equation}\label{eq:def-KK}
\KK(\cV,\pi)[A_1,\dots,A_n]:=\sum_{{\cU\in\cP(n)}\atop \cV\geq {\cU}\geq\pi}
\moeb_{\cP(n)}({\cU},\cV)\cdot\GG(\cU,\pi)[A_1,\dots,A_n].
\end{equation}
we thereby get 
\begin{equation}\label{eq:k-KK}
k_\cV\bigl\{a^{(1)}_{p_1 r_{1}},\dots,a^{(n)}_{p_{n} r_n}\}=
\sum_{\pi\in S_n\atop \pi\leq \cV}
\delta_{r,p\circ\pi}\cdot\KK(\cV,\pi)[A_1,\dots,A_n].
\end{equation}
It
follows that
\begin{align*}
\ff^{}(\cU,\gamma)[A_1,\dots,A_n]&=\sum_{\cV\leq\cU\atop
\cV\vee\gamma=\cU} \sum_{i(1),\dots,i(n)=1}^N k_\cV[a^{(1)}_{i(1)i(\gamma(1))},\dots,
a^{(n)}_{i(n)i(\gamma(n))}]\\
&=\sum_{\cV\leq\cU\atop \cV\vee\gamma=\cU} \sum_{i(1),\dots,i(n)=1}^N \sum_{\pi\in
S_n\atop \pi\leq\cV}
\delta_{i\circ \gamma,i\circ \pi}\cdot \KK(\cV,\pi)[A_1,\dots,A_n]\\
&=\sum_{\cV\leq\cU\atop \cV\vee\gamma=\cU}\sum_{\pi\in S_n\atop
\pi\leq\cV} \KK(\cV,\pi)[A_1,\dots,A_n]\cdot N^{\#(\gamma\pi^{-1})}.
\end{align*}
Since $\cV\vee\gamma=\cU$ is, under the assumption $\pi\leq\cV$,
equivalent to $\cV\vee\gamma\pi^{-1}=\cU$ we can write this also as
\begin{align}\label{eq:moment-cumulant-N}
\ff(\cU,\gamma)[A_1,\dots,A_n] &=\sum_{(\cV,\pi)\in\cPP(n)\atop 
\cV\vee\gamma\pi^{-1}=\cU}\KK(\cV,\pi)[A_1,\dots,A_n] \cdot N^{\#(\gamma\pi^{-1})}.
\end{align}

\begin{remark}

1) Note that although the quantity $\KK$ is defined by \eqref{eq:def-KK} in terms of the
macroscopic moments of the $A_i$, they have also a very concrete meaning in terms of cumulants
of entries of the $A_i$. Namely, if we choose $\pi \in S_n$ and
distinct $1\leq i(1),\dots, i(n)\leq N$ then equation
\eqref{eq:k-KK} becomes, when we set $\cV = 1_n$,
\begin{equation}
\KK(1_n,\pi)[A_1,\dots,A_n]=
\cc_n\bigl(a^{(1)}_{i(1) i(\pi(1))},\dots,a^{(n)}_{i(n)
i(\pi(n))}\bigr)
\end{equation}
as the the only term in the sum that survives is the one for
$\pi$.


2) Equation \eqref{eq:moment-cumulant-N} should be considered as
a kind of moment-cumulant formula in our context, thus it should
contain all information for defining the ``cumulants"
$\KK$ in terms of the moments $\ff$. Actually, we can solve this
linear system of equations for
$\KK$ in terms of $\ff$, by using equation
\eqref{eq:def-KK} to define $\KK$ and equation
\eqref{eq:definitionofGGmulti} for $\GG$.
\begin{align*}
&\KK(\cV,\pi)[A_1,\dots,A_n]\\
&=\sum_{{\cU\in\cP(n)}\atop \cV\geq {\cU}\geq\pi} \moeb_{\cP(n)}({\cU},\cV)\cdot
\sum_{(\cW,\sigma)\in\cPP(n)\atop \cW\leq\cU} \Wg(\cU,\sigma\pi^{-1})
\cdot\ff(\cW,\sigma)[A_1,\dots,A_n]\\
&=\sum_{(\cW,\sigma)\in\cPP(n)} \ff(\cW,\sigma)[A_1,\dots,A_n]\cdot
 \sum_{{\cU\in\cP(n)}\atop \cV\geq {\cU}\geq\pi\vee\cW}
\moeb_{\cP(n)}({\cU},\cV)\cdot
 \Wg(\cU,\sigma\pi^{-1}).
\end{align*}
Thus, by using the relative cumulants of the Weingarten function 
from \eqref{eq:cWg-Wg}, we get finally
\begin{equation}\label{eq:KK-ff-cWg}
\KK(\cV,\pi)[A_1,\dots,A_n]=\sum_{(\cW,\sigma)\in\cPP(n)\atop
W\leq\cV}\ff(\cW,\sigma)[A_1,\dots,A_n]\cdot
 \cWg_{\pi\vee\cW,\cV}(\sigma\pi^{-1}).
\end{equation}

3) One should also note that we have defined the Weingarten function only for $N\geq n$; thus
in the above formulas we should always consider sufficiently large $N$. This is consistent with
the observation that the system of equations \eqref{eq:moment-cumulant-N} might not be
invertible for $N$ too small; the matrix $\bigl(N^{\#(\sigma\pi^{-1})}\bigr)_{\sigma,\pi\in
S_n}$ is invertible for $N\geq n$, however, in general not for all $N<n$ (e.g, clearly not for
$N=1$). One can make sense of some formulas involving the Weingarten function also for $N<n$
(see \cite{CollinsSniady2004}). However, since we are mainly interested in the asymptotic
behavior of our formulas for $N\to\infty$,  we will not elaborate on this.
\end{remark}

\subsection{Product of two independent ensembles}
Let us now calculate the correlation functions for a product of two
independent ensembles $A_1,\dots,A_n$ and $B_1,\dots,B_n$
of random matrices, where we assume that one
of them, let's say the $B_i$'s, is unitarily invariant. We have,
by using \eqref{eq:ff-entries} and
the special version \eqref{eq:product-special-case} of the formula
of Leonov and Shiryaev, the following:

\begin{align*}
&\ff(\cU,\gamma)[A_1B_1,\dots,A_nB_n]\\
& = \mathop{\sum_{i(1), \cdots, i(n)}}_{j(1), \cdots, j(n)}
\sum_{\cV,\cV'\leq
\cU\atop \cV\vee\cV'\vee\gamma=\cU}
\cc_{\cV}[a^{(1)}_{j(1)i(1)},\dots,a^{(n)}_{j(n)i(n)}]\cdot
\cc_{\cV'}[b^{(1)}_{i(1) j(\gamma(1))},
\dots,b^{(n)}_{i(n)j(\gamma(n))}] \\
&\stackrel{\eqref{eq:def-KK}}{=}
\sum_{i,j}\sum_{\cV,\cV'\leq\cU\atop
\cV\vee\cV'\vee\gamma=\cU}
\sum_{\pi\in S_n\atop \pi\leq\cV}\delta_{i, j\circ\pi}\cdot
\KK(\cV,\pi)[A_1,\dots,A_n]\cdot
\cc_{\cV'}[b^{(1)}_{i(1) j(\gamma(1))},
\dots,b^{(n)}_{i(n)j(\gamma(n))}]\\
&=\sum_{\pi\in S_n}\sum_{\cV\in\cP(n)\atop\cU\geq\cV\geq \pi}
\KK(\cV,\pi)[A_1,\dots,A_n]\cdot\\
&\qquad\qquad\qquad\qquad\cdot
\Big( \sum_{\cV'\leq\cU\atop \cV'\vee\cV\vee\gamma=\cU}\sum_{i}
\cc_{\cV'}[b^{(1)}_{i(1)i(\pi^{-1}\gamma(1))},\dots,b^{(n)}_{i(n)
i(\pi^{-1}\gamma(n))}]\Big)\\
\end{align*}
In order to evaluate the second factor we note first that, under
the assumption
$\pi\leq\cV$, the condition $\cV'\vee\cV\vee\gamma=\cU$ is equivalent to $\cV'\vee\cV
\vee\pi^{-1}\gamma=\cU$. Next, we rewrite the sum over all $\cV'\in\cP(n)$ with $\cV'\leq\cU$
and $\cV'\vee\cV\vee\pi^{-1}\gamma=\cU$ as a double sum over all $\cW\in\cP(n)$ with
$\cV\vee\cW=\cU$ and all $\cV'\in\cP(n)$ with $\cV'\leq \cW$ and $\cV'\vee\pi^{-1}\gamma=\cW$.
\begin{align*}
\sum_{\cV'\in\cP(n) \atop \cV'\leq \cU, \cV' \vee
\cV \vee \gamma = \cU} 
&\sum_{i}
\cc_{\cV'} [b^{(1)}_{i(1)i(\pi^{-1}\gamma(1))}, \dots,
b^{(n)}_{i(n)i(\pi^{-1}\gamma(n))}]\\ 
&=\sum_{\cW \in \cP(n) \atop \cV \vee \cW = \cU}
\sum_{\cV' \leq \cW \atop \cV'\vee \pi^{-1} \gamma = \cW} 
\sum_{i}
\cc_{\cV'}[b^{(1)}_{i(1)i(\pi^{-1}\gamma(1))},
\dots, b^{(n)}_{i(n)i(\pi^{-1}\gamma(n))}]\\
&= \sum_{\cW\in\cP(n)\atop \cW\geq \pi^{-1}\gamma,\cV\vee\cW=\cU}
\ff(\cW,\pi^{-1}\gamma)[B_1,\dots,B_n].
\end{align*}
Thus we finally get
\begin{align*}
&\ff(\cU,\gamma)[A_1B_1,\dots,A_nB_n]\\
&=\sum_{\pi\in S_n}\sum_{\cV\in\cP(n)\atop\cU\geq\cV\geq \pi}
 \sum_{\cW\in\cP(n)\atop \cW\geq \pi^{-1}\gamma, \cV\vee\cW=\cU}
\KK(\cV,\pi)[A_1,\dots,A_n]\cdot
\ff(\cW,\pi^{-1}\gamma)[B_1,\dots,B_n]\\
&=\sum_{(\cV,\pi),(\cW,\sigma)\in\cPP(n)\atop
\cV\vee\cW=\cU, \pi\sigma=\gamma}
\KK(\cV,\pi)[A_1,\dots,A_n]\cdot
\ff(\cW,\pi^{-1}\gamma)[B_1,\dots,B_n].
\end{align*}

Let us summarize the result of our calculations in the following theorem. In order to indicate
that our main formulas are valid for any fixed $N$, we will decorate the relevant quantities
with a superscript $^{(N)}$. Note that up to now we have not made any asymptotic consideration.

\begin{theorem}\label{distrib-product}
Let $\cM_N:=M_N\otimes L^\infty(\Omega)$ be an ensemble of $N\times
N$-random matrices. Define correlation functions $\ff_n^{(N)}$ on
$\cM_N$ by ($n\in\NN$, $D_1,\dots,D_n\in\cM_N$)
\begin{equation}
\ff^{(N)}_n(D_1,\dots,D_n):= \cc_n(\Tr(D_1),\dots,\Tr(D_n))
\end{equation}
and corresponding ``cumulant functions" $\KK^{(N)}$ (for $n\leq N$) by
\begin{equation}\label{kkff}
\KK^{(N)}(\cV,\pi)
[A_1,\dots,A_n]  = 
\kern-1em\mathop{\sum_{\cW \in \cP(n),\, 
\pi \in S_n}}_{W \leq \cV} \kern-1em
\ff^{(N)}(\cW,\sigma)[A_1,\dots,A_n]\cdot
 \cWg^{(N)}_{\pi\vee\cW,\cV}(\sigma\pi^{-1}).
\end{equation}
or equivalently by the implicit system of equations
\begin{align}\label{ffkk}
\ff^{(N)}(\cU,\gamma)[D_1,\dots,D_n] 
&=\sum_{\cV,\pi } \KK^{(N)} (\cV,\pi) [D_1,\dots,D_n]
\cdot N^{\#(\gamma\pi^{-1})}.
\end{align}
where the sum is over all $\cV \in \cP(n)$ all $\pi \in S_n$
such that $\pi \leq \cV$ and $\cV \vee \gamma\pi^{-1} = \cU$.

1) Let $\cA_N$ be an algebra of unitarily invariant random
matrices in
$\cM_N$. Then we have for all $n\leq N$, all distinct
$i(1), \dots , i(n)$, all
$A_k=\bigl(a_{ij}^{(k)}\bigr)_{i,j=1}^N\in\cA$, and all $\pi \in
S_n$ that
\begin{equation}\label{eq:KK-k}
\KK^{(N)}(1_n,\pi)[A_1,\dots,A_n]= \cc_n\bigl(a^{(1)}_{i(1) i(\pi(1))},\dots,a^{(n)}_{i(n)
i(\pi(n))}\bigr)
\end{equation}

2) Assume that we have two subalgebras $\cA_N$ and $\cB_N$ of
$\cM_N$ such that \par
{\leftskip1cm
\noindent
\llap{$\diamond$}
$\cA_N$ is a unitarily invariant ensemble,\\
\llap{$\diamond$}
$\cA_N$ and $\cB_N$ are independent.
\par}
Then we have for all $n\in\NN$ with $n\leq N$ and all
$A_1,\dots,A_n\in\cA_N$ and  $B_1,\dots,B_n\in\cB_M$:
\begin{multline}\label{convolution}
\ff^{(N)}(\cU,\gamma)[A_1B_1,\dots,A_nB_n]\\
=\sum_{\cV,\pi,\cW,\sigma}
\KK^{(N)}(\cV,\pi)[A_1,\dots,A_n]\cdot
\ff^{(N)}(\cW,\sigma)[B_1,\dots,B_n].
\end{multline}
where the sum is over all
$\cV, \cW \in \cP(n)$ and all $\pi, \sigma \in S_n$ such that
$\pi
\leq \cV$, $\sigma \leq \cW$, $\cV \vee \cW = \cU$, and $\gamma
= \pi \sigma$. 

\end{theorem}

\subsection{Large $N$ asymptotics for moments and cumulants}
Our main interest in this paper will be the large $N$ limit of
formula \eqref{convolution}. This structure in leading order between
independent ensembles of random matrices which are randomly rotated
against each other will be captured in our abstract notion of higher
order freeness.

Of course, now we must make an assumption about the asymptotic
behavior in $N$ of our correlation functions. We will require that
the cumulants of traces of our random matrices decays in
$N$ with the same order as in the case of Gaussian or Wishart
random matrices. In these cases one has very detailed ``genus
expansions" for those cumulants; see, e.g. \cite{Okounkov,
MingoNica2004annular} and one knows that the $n$-th cumulant of
unnormalized traces in polynomials of those random matrices decays
like $N^{2-n}$ (see e.g. \cite[Thm. 3.1 and Thm.
3.5]{HigherOrderFreeness1}).

\begin{definition}
%
%
Let, for each $N\in\NN$, $B^{(N)}_1, \dots, B^{(N)}_r \subset M_N
\otimes L^{\infty-} (\Omega)$ be $N\times N$-random matrices. 
Suppose that the leading term of the correlation
moments of $B^{(N)}_1, \dots, B^{(N)}_r$ are of order $2 - n$,
i.e., that for all $n \in \NN$ and all polynomials $p_1, \dots,
p_t$ in $r$ non-commuting variables the limits
$$
\lim_{N\to\infty}\ff^{(N)}_n(
p_1(B^{(N)}_1, \dots, B^{(N)}_r), \dots,
p_t(B^{(N)}_1, \dots, B^{(N)}_r )) \cdot N^{n-2}
$$ 
exist. Then we will say that $\{B^{(N)}_1, \dots ,
B^{(N)}_r\}$ has {\em limit distributions of all orders}. Let
$\cB$ be the free algebra generated by generators $b_1, \dots,
b_r$. Then we define the {\em limit correlation functions} of
$\cB$ by  
\begin{align*}
\ff_n(& p_1(b_1, \dots, b_r), \dots ,
p_t(b_1, \dots, b_r) ) \\ 
&=
\lim_{N\to\infty}\ff^{(N)}_n(
p_1(B^{(N)}_1, \dots, B^{(N)}_r), \dots,
p_t(B^{(N)}_1, \dots, B^{(N)}_r )) \cdot N^{n-2} 
\end{align*}
\end{definition}

Note
that this assumption implies that the leading term for the quantities
$\ff^{(N)}(\cV,\pi)$ is of order $2\#(\cV)-\#(\pi)$. Indeed, if
$\cV$ has $k$ blocks and the $i^{th}$ block of $\cV$ contains $r_i$
cycles of $\pi$ then $\ff^{(N)}(\cV, \pi) = \ff_{r_1} \cdots
\ff_{r_k}$ and each $\ff_{r_i}$ has order ${2 - r_i}$. Then the
order of $\ff^{(N)}(\cV, \pi)$ is $(2 - r_1) + \cdots + (2 - r_k)
= 2 k - (r_1 + \cdots + r_k) = 2\, \#(\cV) - \#(\pi)$. Thus
\begin{align*}
\ff(& \cV,\pi)( p_1(b_1, \dots, b_r), \dots ,
p_t(b_1, \dots, b_r) ) \\ 
&=
\lim_{N\to\infty}\ff^{(N)}(\cV,\pi)(
p_1(B^{(N)}_1, \dots, B^{(N)}_r), \dots, 
p_t(B^{(N)}_1, \dots, B^{(N)}_r )) \\
& \hskip0.75in \cdot N^{-2\#(\cV)+\#(\pi)} 
\end{align*}

{}From formula (\ref{ffkk}) one can deduce that the leading order
of $\KK^{(N)}(\cV,\pi)$ is given by the term
$(\cU,\gamma)=(\cV,\pi)$ and thus must be of order
$$N^{-n+2\#\cV-\#\pi}.$$
(Indeed, this also follows from equation
\eqref{eq:KK-ff-cWg} and the leading order
of the relative cumulant of the Weingarten
function given in equation \eqref{eq:max-order-Wg}.)

Thus we can define the \emph{limiting cumulant functions}
to be the limit of the leading order of the cumulants by the
equation
\begin{equation}\label{eq:KK-limit}
\KK(\cV,\pi)[b_1,\dots,b_n]:=\lim_{N\to\infty}
N^{n-2\#\cV+\#\pi} \cdot \KK^{(N)} (\cV,\pi) [B^{(N)}_1,
\dots, B^{(N)}_n]
\end{equation}

When $(\cV, \pi) = (1_n, \gamma_n)$ and $B_1 = B_2 = \cdots =
B_n = B$ equation \eqref{eq:KK-k} becomes
$$
\KK^{(N)}(1_n, \gamma_n)[B, \dots , B] =
\cc_n\bigl(b^{(1)}_{i(1) i(2)},\dots,b^{(n)}_{i(n)
i(1)})
$$
Thus to prove Theorem \ref{firstlimitcumulant} we must show
that $\KK^{(N)}(1_n, \gamma_n)[B, \dots , B]\ab \cdot N^{n-1}$
converges to $\KK^b_n$ the $n^{th}$ free cumulant of the limiting
eigenvalue distribution of $B^{(N)}$. 

When $(\cV, \pi) = (1_{m+n}, \gamma_{m,n})$ equation
\eqref{eq:KK-k} becomes
$$
\KK^{(N)}(1_{m+n}, \gamma_{m,n})[B, \dots , B] =
\cc_{m+n}\bigl(b^{(1)}_{i(1) i(2)},\dots,b^{(n)}_{i(n)
i(1)})
$$
Thus to prove Theorem \ref{secondlimitcumulant} we must show
that $\KK^{(N)}(1_{m+n}, \gamma_{m,n})[B, \ab \dots , B]
\cdot N^{m+n}$ converges to $\KK^b_{m,n}$ the $(m,n)^{th}$ free
cumulant of second order of the limiting second order
distribution of
$B^{(N)}$.

\subsection{Length functions}
We want to understand the asymptotic behavior of formula (\ref{convolution}). The leading order
in $N$ of the right hand side is given by
$$-n+2\#\cV-\#\pi+2\#\cW-\#\sigma=n+(\vert\pi\vert-2\vert\cV\vert)
+(\vert\sigma\vert-2\vert\cW\vert),$$ whereas the leading order of
the left hand side is given by
$$2\#\cU-\#\gamma=2\#(\cV\vee\cW)-\#(\sigma\pi)=n+(\vert\pi\sigma
\vert-2\vert \cV\vee\cW\vert).$$
This suggests the introducing of the
following ``length functions" for permutations, partitions, and
partitioned permutations.
\begin{notation} \
\begin{enumerate}
\item For $\cV\in\cP(n)$ and $\pi\in S_n$ we put
\begin{align*}
\vert\cV\vert&:=n-\#\cV\\
\vert\pi\vert&:=n-\#\pi.
\end{align*}
\item For any $(\cV,\pi)\in\cPP(n)$ we put
$$\vert(\cV,\pi)\vert:=2\vert\cV\vert-\vert\pi\vert=n-(2\#\cV-\#\pi).$$
\end{enumerate}
\end{notation}

Let us first observe that these quantities behave
actually like a length. It is clear from the definition that they are always non-negative;
that they also obey a
triangle inequality is the content of the next lemma.

\begin{lemma}\label{lemma:triangle} \
\begin{enumerate} \item For all $\pi,\sigma\in S_n$ we have
$$\vert\pi\sigma\vert\leq \vert\pi\vert+\vert\sigma\vert.$$
\item For all $\cV,\cW\in\cP(n)$ we have
$$\vert \cV\vee\cW\vert\leq \vert\cV\vert+\vert\cW\vert .$$
\item For all partitioned permutations
$(\cV,\pi),(\cW,\sigma)\in\cPP(n)$ we have
$$\vert(\cV\vee\cW,\pi\sigma)\vert\leq \vert(\cV,\pi)\vert
+\vert(\cW,\sigma)\vert.$$
\end{enumerate}
\end{lemma}

\begin{proof}
(1) This is well-known, since $\vert\pi\vert$ is the minimal number
of factors needed to write $\pi$ as a product of transpositions.

(2) Each block $B$ of $\cW$ can glue at most $\#B-1$ many blocks of
$\cV$ together, i.e., $\cW$ can glue at most $n-\#\cW$ many blocks
of $\cV$ together, thus the difference between $\vert\cV\vert$ and
$\vert\cV\vee\cW\vert$ cannot exceed $n-\#\cW$ and hence
$$\#\cV-\#(\cV\vee\cW)\leq n-\#\cW.$$
This is equivalent to our assertion.

(3) We prove this, for fixed $\pi$ and $\sigma$ by induction over
$\vert\cV\vert+\vert\cW\vert$. The smallest possible value of the
latter appears for $\vert\cV\vert=\vert\pi\vert$ and
$\vert\cW\vert=\vert\sigma\vert$ (i.e., $\cV=0_\pi$ and
$\cW=0_\sigma$). But then we have (since $\cV\vee\cW\geq \pi\sigma$)
$$2\vert\cV\vee\cW\vert-\vert\pi\sigma\vert
\leq \vert\cV\vee\cW\vert\leq \vert\cV\vert+\vert\cW\vert,$$ which
is exactly our assertion for this case. For the induction step, on
the other side, one only has to observe that if one increases
$\vert\cV\vert$ (or $\vert\cW\vert$) by one then
$\vert\cV\vee\cW\vert$ can also increase by at most 1.
\end{proof}

\begin{remark}
1) Note that the triangle inequality for partitioned
permutations together with \eqref{convolution} implies the
following. Given random matrices  $A=(A_N)_{N\in\NN}$ and
$B=(B_N)_{N\in\NN}$ which have limit distributions of all orders.
If $A$ and $B$  are independent and at least one of them is
unitarily invariant, then $C=(C_N)_{N\in\NN}$ with $C_N := A_N
B_N$ also has limit distributions of all orders. 

2) Since we know that Gaussian and  Wishart random matrices have
limit distributions of all orders (see e.g. \cite[Thm. 3.1 and
Thm. 3.5]{HigherOrderFreeness1}), and since they are unitarily
invariant, it follows by induction from the previous part that
any polynomial in independent Gaussian and Wishart matrices has
limit distributions of all orders.
\end{remark}

\subsection{Multiplication of partitioned permutations}

Suppose $\{ B^{(N)}_1,\ab \dots , B^{(N)}_n \}$ has limit
distributions of all orders. Then the left hand side of equation 
\eqref{ffkk} has order $N^{2\#(\cU) - \#(\gamma)}$ and the right
hand side of equation \eqref{ffkk} has order $N^{-n + 2 \#(\cV) -
\#(\pi) + |\gamma\pi^{-1}|}$. Thus the only terms of the right
hand side that have order $N^{2 \#(\cU) - \#(\gamma)}$ are those
for which 
$$
2 \#(\cU) - \#(\gamma) =
-n + 2 \#(\cV) - \#(\pi) + |\gamma\pi^{-1}|
$$
i.e. for which $|(\cU, \gamma)| = |(\cV, \pi)| + |\gamma
\pi^{-1}|$. Hence
\begin{align*}
\ff^{(N)}(\cU,\gamma) [B^{(N)}_1, & \dots, B^{(N)}_n] \\ \notag
= \kern-2em &
\sum_{(\cV,\pi)\in\cPP(n)\atop
{\cV \vee \gamma\pi^{-1} = \cU 
\atop
\vert(\cU, \gamma)\vert = \vert(\cV, \pi)\vert
+\vert\gamma\pi^{-1}\vert}} \kern-2em
\KK^{(N)}(\cV, \pi)[B^{(N)}_1,  \dots,B^{(N)}_n] \cdot
N^{|\gamma\pi^{-1}|} \\ 
& \qquad\qquad\qquad\qquad\qquad\mbox{} + O(N^{ 2\#(\cU) -
\#(\gamma) -2}) \notag
\end{align*}
Thus after taking limits we have
\begin{equation}\label{ffkk-leading}
\ff(\cU, \gamma)[b_1, \dots , b_n] 
= \kern-0.5em
\sum_{(\cV,\pi)\in\cPP(n)} \kern-0.5em
\KK(\cV, \pi)[b_1,  \dots,b_n]
\end{equation}
where the sum is over all $(\cV, \pi)$ in $\cPP(n)$ such that 
$\cV \vee \gamma\pi^{-1} = \cU$ and $\vert(\cU, \gamma)\vert =
\vert(\cV, \pi)\vert +\vert\gamma\pi^{-1}\vert$.

A similar analysis of equation \eqref{convolution} gives that for
independent $\{A^{(N)}_1,\ab \dots , A^{(N)}_n\}$ and
$\{B^{(N)}_1,  \dots , B^{(N)}_n\}$ with the $A^{(N)}_i$'s
unitarily invariant and both having limit distributions of all
orders we have

\begin{multline*}
\ff^{(N)}(\cU,\gamma)[A^{(N)}_1 B^{(N)}_1, 
                      \dots, A^{(N)}_n B^{(N)}_n]   \\ 
= \kern-2.5em
\sum_{(\cV,\pi),(\cW,\sigma) \in \cPP(n)
\atop 
{\cV \vee \cW = \cU,\, \pi\sigma = \gamma
\atop 
\vert(\cV, \pi)\vert + \vert(\cW, \sigma)\vert
= \vert( \cV \vee \cW, \pi \sigma )\vert}} \kern-3em
\KK^{(N)}(\cV,\pi)[A^{(N)}_1,\dots,A^{(N)}_n]\cdot
\ff^{(N)}(\cW,\sigma)[B^{(N)}_1, \dots, B^{(N)}_n] \\
\quad \mbox{} + O(N^{2 \#(\cU) - \#(\gamma) -2})
\end{multline*}
and again after taking limits
\begin{multline}\label{convolution-leading}
\ff(\cU,\gamma)[a_1 b_1, 
                      \dots, a_n b_n]   \\ 
= \kern-2.5em
\sum_{(\cV,\pi),(\cW,\sigma) \in \cPP(n)}
\KK(\cV,\pi)[a_1,\dots, a_n]\cdot
\ff(\cW,\sigma)[b_1, \dots, b_n] 
\end{multline}
where the sum is over all $(\cV,\pi),(\cW,\sigma) \in \cPP(n)$
such that 
\begin{itemize}
\item[$\diamond$] $\cV \vee \cW = \cU$
\item[$\diamond$] $\pi \sigma = \gamma$
\item[$\diamond$] $|(\cV, \pi)| + |(\cW, \sigma)| = |(\cU,
\gamma)|$
\end{itemize}

In order to write this in a more compact form it is convenient to
define a multiplication for partitioned permutations (in $\CC\cPP(n)$) as follows.

\begin{definition}\label{def:mult}
For $(\cV,\pi),(\cW,\sigma)\in\cPP(n)$ we define their product as
follows.
\begin{multline}
(\cV,\pi)\cdot(\cW,\sigma):=\\
=\begin{cases} (\cV\vee\cW,\pi\sigma) &\text{if
$\vert(\cV,\pi)\vert+ \vert(\cW,\sigma)\vert=
\vert (\cV\vee\cW,\pi\sigma)\vert$,}\\
0 &\text{otherwise.}
\end{cases}
\end{multline}
\end{definition}

\begin{proposition}
The multiplication defined in Definition \ref{def:mult} is
associative.
\end{proposition}

\begin{proof}
We have to check that
\begin{equation}\label{eq:assoc}
\bigl((\cV,\pi)\cdot (\cW,\sigma)\bigr)\cdot (\cU,\tau)=
(\cV,\pi)\cdot\bigl( (\cW,\sigma)\cdot (\cU,\tau)\bigr).
\end{equation}
Since both sides are equal to $(\cV\vee\cW\vee\cU,\pi\sigma\tau)$ in
case they do not vanish, we have to see that the conditions for
non-vanishing are for both sides the same.

The conditions for the left hand side are
$$\vert(\cV,\pi)\vert+\vert(\cW,\sigma)\vert=\vert(\cV\vee\cW,\pi\sigma)
\vert$$ and
$$\vert(\cV\vee\cW,\pi\sigma)\vert+\vert(\cU,\tau)\vert=
\vert(\cV\vee\cW\vee\cU,\pi\sigma\tau)\vert.$$ These imply
\begin{align*}
\vert(\cV,\pi)\vert+\vert(\cW,\sigma)\vert+\vert(\cU,\tau)\vert&=
\vert(\cU\vee\cW\vee\cU,\pi\sigma\tau)\vert\\
&\leq \vert(\cV,\pi)\vert+\vert(\cW\vee\cU,\sigma\tau)\vert,
\end{align*}
However, the triangle inequality
$$\vert(\cW\vee\cU,\sigma\tau)\vert\leq \vert(\cW,\sigma)\vert+
\vert(\cU,\tau)\vert$$ yields that we have actually equality in the
above inequality, thus leading to
$$\vert(\cW,\sigma)\vert+\vert(\cU,\tau)\vert=\vert(\cW\vee\cU,\sigma\tau)
\vert$$ and
$$\vert(\cV,\pi)\vert +
\vert(\cW\vee\cU,\sigma\tau)\vert=
\vert(\cV\vee\cW\vee\cU,\pi\sigma\tau)\vert.$$ These are exactly the
two conditions for the vanishing of the right hand side of
(\ref{eq:assoc}). The other direction goes analogously
\end{proof}

Now we can write formulas (\ref{ffkk-leading}) and
(\ref{convolution-leading}) in convolution form
\begin{equation}\label{ffkk-leading-compact}
\ff(\cU,\gamma)[b_1,\dots,b_n] 
= \kern-2em \sum_{(\cV,\pi)\in\cPP(n)\atop
(\cV,\pi)\cdot(0,\gamma\pi^{-1})=(\cU,\gamma)}
\kern-2em
\KK(\cV,\pi)[b_1,\dots,b_n] 
\end{equation}
and
\begin{multline}\label{convolution-leading-compact}
\ff(\cU,\gamma) [a_1 b_1, \dots,
                       a_n b_n]   \\ 
= \kern-1em
\sum_{(\cU,\pi),(\cW,\sigma)\in\cPP(n)\atop
(\cV,\pi)\cdot(\cW,\sigma)=(\cU,\gamma)}
\kern-2em
\KK(\cV,\pi)[a_1,\dots, a_n] \cdot
\ff(\cW,\sigma)[b_1, \dots, b_n] 
\end{multline}
Note that both $\ff(\cV,\pi)$ and $\KK(\cV,\pi)$ are multiplicative
in the sense that they factor according to the decomposition of
$\cV$ into blocks.

The philosophy for our definition of higher order freeness
will be that equation \eqref{ffkk-leading-compact} is the
analogue of the moment-cumulant formula and shall be used to
define the quantities
$\KK$, which will thus take on the role of cumulants in our
theory --  whereas the $\ff$ are the moments (see Definition
\ref{def:freecumulants}). We shall define higher order freeness
by requiring the vanishing of mixed cumulants, see Definition
\ref{def:higherfreeness}. On the other hand, equation
\eqref{convolution-leading-compact} would be another
way of expressing the fact that the $a$'s are free from
the $b$'s. Of course, we will have to prove that those two
possibilities are actually equivalent (see Theorem
\ref{thm:charact-freeness}).

\section{Multiplicative functions on partitioned permutations
and their convolution}
\label{sec:multiplicative}
\subsection{Convolution of multiplicative functions}
Formulas \eqref{ffkk-leading-compact} and
\eqref{convolution-leading-compact} above are a
generalization of the formulas describing first order
freeness in terms of cumulants and convolution of
multiplicative functions on non-crossing partitions. Since
the dependence on the random matrices is irrelevant for
this structure we will free ourselves in this section from
the random matrices and look on the combinatorial heart of
the observed formulas. In Section
\ref{sec:higher}, we will return to the more general
situation involving multiplicative functions which depend
also on random matrices or more generally elements from
an algebra.

\begin{definition}\  \nopagebreak
\begin{enumerate} \item  We denote by $\cPP$ the set of partitioned permutations on an
arbitrary number of elements, i.e.,
$$\cPP=\bigcup_{n\in\NN}\cPP(n).$$
\item For two functions
$$f,g:\cPP\to\CC$$
we define their convolution
$$f*g:\cPP\to\CC$$
by $$(f*g)(\cU,\gamma):=\sum_{(\cV,\pi),(\cW,\sigma)\in\cPP(n)\atop
(\cV,\pi)\cdot (\cW,\sigma)=(\cU,\gamma)} f(\cV,\pi)\ g(\cW,\sigma)
$$
for any $(\cU,\gamma)\in\cPP(n)$.
\end{enumerate}
\end{definition}

\begin{definition} A function $f:\cPP\to\CC$ is called
\emph{multiplicative} if $f(1_n, \pi)$ depends
only on the conjugacy class of $\pi$ and we have
$$f(\cV,\pi)=\prod_{V\in\cV}f(1_V,\pi\vert_V)$$
\end{definition}

Our main interest will be in multiplicative functions. It is easy to
see that the convolution of two multiplicative functions is again
multiplicative. It is clear that a multiplicative function is
determined by the values of $f(1_n,\pi)$ for all $n\in\NN$ and all
$\pi\in S_n$.

An important example of a multiplicative function is the
$\delta$-function presented below.

\begin{notation}
The $\delta$-function on $\cPP$ is the multiplicative function
determined by
$$\delta(1_n,\pi)=\begin{cases}
1,&\text{if $n=1$}\\
0,&\text{otherwise}.
\end{cases}$$
\end{notation}
Thus for $(\cU, \pi) \in \cPP(n)$
$$\delta(\cU,\pi)=\begin{cases}
1,&\text{if $(\cU,\pi)=\bigl(0_n,(1)(2)\dots (n)\bigr)$ for some $n$}\\
0,&\text{otherwise},
\end{cases}$$

\begin{proposition}
The convolution of multiplicative functions on $\cPP$ is
commutative and $\delta$ is the unit element.
\end{proposition}

\begin{proof}
It is clear that $\delta$ is the unit element. For commutativity, we
note that for multiplicative functions we have
$$f(\cV,\pi)=f(\cV,\pi^{-1}),$$
and thus
$$(g*f)(\cU,\gamma)=(g*f)(\cU,\gamma^{-1})=
\sum_{(\cV,\pi),(\cW,\sigma)\in\cPP(n)\atop (\cV,\pi)\cdot
(\cW,\sigma)=(\cU,\gamma^{-1})} g(\cV,\pi)f(\cW,\sigma).
$$
Since the condition $(\cV,\pi)\cdot (\cW,\sigma)=(\cU,\gamma^{-1})$
is equivalent to the condition $(\cW,\sigma^{-1})\cdot
(\cV,\pi^{-1})=(\cU,\gamma)$ we can continue with
$$(g*f)(\cU,\gamma)=
\sum_{(\cV,\pi),(\cW,\sigma)\in\cPP(n)\atop (\cW,\sigma^{-1})\cdot
(\cV,\pi^{-1})=(\cU,\gamma)} f(\cW,\sigma^{-1}) g(\cV,\pi^{-1})
=(f*g)(\cU,\gamma).$$

\end{proof}

\subsection{Factorizations}
Let us now try to characterize the non-trivial factorizations
$(\cU,\gamma)=(\cV,\pi)\cdot (\cW,\sigma)$ appearing in the
definition of our convolution. Let us first observe some simple
general inequalities.

\begin{lemma}\label{lemma:comp} \
\begin{enumerate} \item  For permutations $\pi,\sigma\in S(n)$ we have
$$
\vert\pi\vert+\vert\sigma\vert+\vert\pi\sigma\vert \geq
2\vert\pi\vee\sigma\vert.
$$
\item For partitions $\cV_2\leq \cV_1$ and $\cW_2\leq \cW_1$ we have
$$
\vert\cW_1\vert+\vert\cV_1\vert+  \vert\cV_2\vee\cW_2\vert\geq \vert
\cV_1\vee \cW_1\vert+\vert\cW_2\vert +\vert\cV_2\vert
$$
and
$$
\vert\cV_1\vee\cW_2\vert+\vert\cV_2\vee\cW_1\vert \geq \vert
\cV_1\vee\cW_1\vert+\vert\cV_2\vee\cW_2\vert.
$$
\end{enumerate}
\end{lemma}

\begin{proof}
(1) By the triangle inequality for partitioned permutations we have
$$\vert(0_\pi\vee 0_{\sigma},\pi\sigma)\vert\leq
\vert(0_\pi,\pi)\vert+\vert(0_{\sigma},\sigma)\vert,$$ i.e.,
\begin{equation}\label{eq:hier}
2\vert\pi\vee\sigma\vert-\vert\pi\sigma\vert\leq
\vert\pi\vert+\vert\sigma\vert.
\end{equation}

(2) Consider first the special case $\cW_1=\cW_2=\cW$. Then we
clearly have
$$\#(\cV_2\vee\cW)-\#(\cV_1\vee\cW)\leq \#\cV_2-\#\cV_1,$$
which leads to
$$\vert\cV_1\vee\cW\vert-\vert\cV_2\vee\cW\vert\leq \vert\cV_1\vert-
\vert\cV_2\vert.$$ {}From this the general case follows by
\begin{align*}
\vert\cV_1\vee\cW_1\vert-\vert\cV_2\vee\cW_2\vert&=
\vert\cV_1\vee\cW_1\vert-\vert\cV_1\vee\cW_2\vert+\vert\cV_1\vee\cW_2\vert
-\vert\cV_2\vee\cW_2\vert\\
&\leq
\vert\cW_1\vert-\vert\cW_2\vert+\vert\cV_1\vert-\vert\cV_2\vert.
\end{align*}
The second inequality follows from this as follows:
\begin{align*}
\vert\cV_1\vee\cW_1\vert-\vert \cV_1\vee\cW_2\vert&=
\vert\cV_1\vee(\cV_2\vee\cW_1)
\vert-\vert \cV_1\vee (\cV_2\vee\cW_2)\vert\\
&\leq \vert\cV_2\vee\cW_1\vert-\vert\cV_2\vee\cW_2\vert.
\end{align*}
\end{proof}

\begin{theorem}\label{four-conditions}
\label{the:NAME} For $(\cV,\pi),(\cW,\sigma)\in\cPP(n)$ the equation
$$(\cV,\pi)\cdot (\cW,\sigma)=(\cV\vee\cW,\pi\sigma)$$
is equivalent to the conjunction of the following four conditions:
\begin{align*}
\vert\pi\vert+\vert\sigma\vert+\vert\pi\sigma\vert&= 2\vert\pi\vee\sigma\vert,\\
\vert\cV\vert+\vert\pi\vee\sigma\vert& =\vert\pi\vert+\vert\cV\vee\sigma\vert,\\
\vert\cW\vert+\vert\pi\vee\sigma\vert&= \vert\sigma\vert+\vert\pi\vee\cW\vert,\\
\vert\cV\vee\sigma\vert+\vert\pi\vee\cW\vert&=\vert\cV\vee\cW\vert+\vert\pi
\vee\sigma\vert.
\end{align*}
\end{theorem}

\begin{proof}
Adding the four inequalities given by Lemma \ref{lemma:comp}
\begin{align*}
\vert\pi\vert+\vert\sigma\vert+\vert\pi\sigma\vert&\geq 2\vert\pi\vee\sigma\vert,\\
2\vert\cV\vert+2\vert\pi\vee\sigma\vert&\geq 2\vert\pi\vert+2\vert\cV\vee\sigma\vert,\\
2\vert\cW\vert+2\vert\pi\vee\sigma\vert&\geq2\vert\sigma\vert+2\vert\pi\vee\cW\vert,\\
2\vert\cV\vee\sigma\vert+2\vert\pi\vee\cW\vert&\geq2\vert\cV\vee\cW\vert+
2\vert\pi\vee\sigma\vert
\end{align*}
gives
$$2\vert\cV\vert-\vert\pi\vert+2\vert\cW\vert-\vert\sigma\vert\geq
2\vert\cV\vee\cW\vert-\vert\pi\sigma\vert,
$$
i.e.,
$$\vert(\cV,\pi)\vert+\vert(\cW,\sigma)\vert\geq\vert(\cV\vee\cW,\pi\sigma).$$
Since $(\cV,\pi)\cdot (\cW,\sigma)=(\cV\vee\cW,\pi\sigma)$ means
that we require equality in the last inequality, this is equivalent
to having equality in all the four inequalities.
\end{proof}

The conditions describing our factorizations have a quite
geometrical meaning. Let us elaborate on this in the
following.

\begin{definition}\label{minimal-planar-definition}
Let $\gamma\in S(n)$ be a fixed permutation.
\begin{enumerate} \item
A permutation $\pi\in S(n)$ is called \emph{$\gamma$-planar} if
$$\vert\pi\vert+\vert\pi^{-1}\gamma\vert+\vert\gamma\vert=2\vert\pi\vee
\gamma\vert.$$ \item A partitioned permutation $(\cV,\pi)\in\cPP(n)$
is called \emph{$\gamma$-minimal} if
$$\vert\cV\vee\gamma\vert-\vert\pi\vee\gamma\vert=\vert\cV\vert-
\vert\pi\vert.$$ \end{enumerate}
\end{definition}

\begin{remark} \label{minimal-planar-remark}
{\em i}\,)
It is easy to check (for example, by
calculating the Euler characteristic) that
$\gamma$-planarity of $\pi$ corresponds indeed to a planar
diagram, i.e.\ one can draw a planar graph representing
permutations $\gamma$ and
$\pi$ without any crossings. The most important cases are
when $\gamma$ consists of a single cycle
\cite{Biane1997crossings} and when $\gamma$ consists of
two cycles \cite{MingoNica2004annular}. \\
{\em ii}\,) 
The notion of $\gamma$-minimality of $(\cV,\pi)$
means that $\cV$ connects only blocks of $\pi$ which are
not already connected by $\gamma$. \\
{\em iii}\,)
If $(\cV, \pi)$ satisfies both (1) and (2) of
Definition \ref{minimal-planar-definition} then $(\cV, \pi) (0,
\pi^{-1} \gamma) = (1, \gamma)$, by Theorem
\ref{four-conditions}.
\end{remark}

\begin{corollary}
Assume that we have the equation
$$(\cU,\gamma)=(\cV,\pi)\cdot (\cW,\sigma).$$
Then $\pi$ and $\sigma$ must be $\gamma$-planar and $(\cV,\pi)$ and
$(\cW,\sigma)$ must be $\gamma$-minimal.
\end{corollary}

\subsection{Factorizations of disc and tunnel permutations}

\begin{notation}
{\em i }\,)
We call $(\cV,\pi)\in \cPP_n$ a \emph{disc permutation} if
$\cV=0_\pi$; the latter is equivalent to the condition
$\vert\cV\vert=\vert\sigma\vert$. For
$\pi\in S_n$, by $(0,\pi)$ we will always mean the disc
permutation
$$(0,\pi):=(0_\pi,\pi)\in\cPP(n).$$
{\em ii}\,)
We call $(\cV, \pi) \in \cPP_n$ a
\emph{tunnel permutation} if $\vert \cV \vert = \vert \pi
\vert + 1$. This means that $\cV$ is obtained from $\pi$
by joining a pair of cycles; i.e. one block of $\cV$
contains exactly two cycles of $\pi$ and all other blocks
contain only one cycles of $\pi$. 
\end{notation}

A motivation for those names comes from the identification between
partitioned permutations and so-called surfaced permutations; see
the Appendix for more information on this.

Our goal is now to understand more explicitly the factorizations
of disc and tunnel permutations. (It will turn out that those are
the relevant ones for first and second order freeness). For this, note that
we can rewrite the crucial condition for our product of partitioned
permutations,
$$2\vert\cV\vert-\vert\pi\vert+
2\vert\cW\vert-\vert\sigma\vert= 2\vert
\cV\vee\cW\vert-\vert\pi\sigma\vert,$$ in the form
$$\bigl(\vert\cV\vert-\vert\pi\vert\bigr)+
\bigl(\vert\cW\vert-\vert\sigma\vert\bigr)+
\bigl(\vert\cV\vert+\vert\cW\vert-\vert\cV\vee\cW\vert\bigr)= \bigl(\vert
\cV\vee\cW\vert-\vert\pi\sigma\vert\bigr).$$ Since all terms in brackets are non-negative
integers this formula can be used to obtain explicit solutions to our factorization problem for
small values of the right hand side. Essentially, this tells us that factorizations of a disc
permutation can only be of the form $\text{disc}\times\text{disc}$; and factorizations of a
tunnel permutation can only be of the form $\text{disc}\times\text{disc}$,
$\text{disc}\times\text{tunnel}$, and
$\text{tunnel}\times\text{disc}$. Of course, one can generalize
the following arguments to higher order type permutations,
however, the number of possibilities grows quite quickly.

\begin{proposition}
\label{pro:types-of-factorizations} \
\begin{enumerate}
\item The solutions to the equation
$$(1_n,\gamma_n)=(0,\gamma_n)=(\cV,\pi)\cdot(\cW,\sigma)$$
are exactly of the form
$$(1_n,\gamma_n)=(0,\pi)\cdot (0,
\pi^{-1}\gamma_n),$$ for some $\pi\in NC(n)$. \item The
solutions to the equation
$$(1_{m+n},\gamma_{m,n})=(\cV,\pi)\cdot (\cW,\sigma)$$
are exactly of the following three forms:
\begin{enumerate}
\item
$$(1_{m+n},\gamma_{m,n})=(0,\pi)\cdot (0,
\pi^{-1}\gamma_{m,n}),$$ where $\pi\in S_{NC}(m,n)$;
\item
$$(1_{m+n},\gamma_{m,n})=(0,\pi)\cdot (\cW,
\pi^{-1}\gamma_{m,n}),$$ where $\pi\in NC(m)\times NC(n)$ and
$\vert\cW\vert=\vert \pi^{-1}\gamma_{m,n}\vert+1$;
\item
$$(1_{m+n},\gamma_{m,n})=(\cV,\pi)\cdot (0,
\pi^{-1}\gamma_{m,n}),$$ where $\pi\in NC(m)\times NC(n)$ and
$\vert\cV\vert=\vert \pi\vert+1$.
\end{enumerate}
\end{enumerate}
\end{proposition}

\begin{proof}
(1) The correspondence between non-crossing partitions and
permutations was studied in detail by Biane
\cite{Biane1997crossings}. In this case we have
$$\bigl(\vert\cV\vert-\vert\pi\vert\bigr)+
\bigl(\vert\cW\vert-\vert\sigma\vert\bigr)+
\bigl(\vert\cV\vert+\vert\cW\vert-\vert\cV\vee\cW\vert\bigr)= \vert
1_n\vert-\vert\gamma_n\vert=0.$$ Since all three terms in brackets
are greater or equal to zero, all of them must vanish, i.e.,
$$
\vert\cV\vert=\vert\pi\vert,\qquad\text{thus}\qquad \cV=0_\pi$$
$$\vert\cW\vert=\vert\sigma\vert,\qquad\text{thus}\qquad
\cW=0_\sigma$$ and
$$\vert\pi\vert+\vert\sigma\vert=
\vert\cV\vert+\vert\cW\vert=\vert\cV\vee\cW\vert=
\vert\gamma\vert=n-1.
$$

(2) Now we have
$$\bigl(\vert\cV\vert-\vert\pi\vert\bigr)+
\bigl(\vert\cW\vert-\vert\sigma\vert\bigr)+
\bigl(\vert\cV\vert+\vert\cW\vert-\vert\cV\vee\cW\vert\bigr)=
\bigl(\vert \cV\vee\cW\vert-\vert\pi\sigma\vert\bigr)=1,$$ which
means that two of the terms on the left-hand side must be equal to
0, and the other term must be equal to 1. Thus we have the following
three possibilities.
\def\theenumi{\alph{enumi}}
\begin{enumerate}
\item
\begin{align*} \vert\cV\vert=&\vert\pi\vert,\qquad\text{thus}
&
\cV&=0_\pi, \\
\vert\cW\vert=&\vert\sigma\vert,\qquad\text{thus} & \cW&=0_\sigma
\end{align*}
and
$$\vert\pi\vert+\vert\sigma\vert=
\vert\cV\vert+\vert\cW\vert=\vert\cV\vee\cW\vert+1=m+n.$$ Note that
$$\pi\vee\sigma=\cV\vee\cW=1_{m+n},$$
and thus $\pi$ connects the two cycles of $\gamma_{m,n}$. This means
that $\pi$ is a non-crossing $(m,n)$-permutation.
\item
$$\vert\cV\vert=\vert\pi\vert,\qquad\text{thus}\qquad
\cV=0_\pi,$$
$$\vert\cW\vert=\vert\sigma\vert+1,$$
and
$$\vert\cV\vert+\vert\cW\vert=\vert\cV\vee\cW\vert=m+n-1.$$
This implies
$$\vert\pi\vert+\vert \gamma_{m,n}\pi^{-1}\vert=m+n-2,$$
which means that $\pi$ must be a disconnected non-crossing
$(m,n)$-annular permutation, i.e.,
$$\pi=\pi_1\times\pi_2\qquad\text{with}\qquad
\pi_1\in NC(m), \pi_2\in NC(n).$$
\item
$$\vert\cV\vert=\vert\pi\vert+1,$$
$$\vert\cW\vert=\vert\sigma\vert+1,,\qquad\text{thus}\qquad
\cW=0_\sigma$$ and
$$\vert\cV\vert+\vert\cW\vert=\vert\cV\vee\cW\vert=m+n-1.$$
This implies
$$\vert\pi\vert+\vert \gamma_{m,n}\pi^{-1}\vert=m+n-2,$$
which means that $\pi$ must be a disconnected non-crossing
$(m,n)$-annular permutation, i.e.,
$$\pi=\pi_1\times\pi_2\qquad\text{with}\qquad
\pi_1\in NC(m), \pi_2\in NC(n).$$
\end{enumerate}
\end{proof}
\def\theenumi{\arabic{enumi}}

\begin{example}
We can now use the previous description of factorizations of disc
and tunnel permutations to write down explicit first and second
order formulas for our convolution of multiplicative functions.

1) In the first order case we have
\begin{equation}
(f*g)(1_n,\gamma_n)=(f*g)(0,\gamma_n)= \sum_{\pi\in
NC(n)}f(0,\pi)g(0, \pi^{-1}\gamma_n).
\end{equation}
This equation is exactly the formula for the convolution of
multiplicative functions on non-crossing partitions, which is the
cornerstone of the combinatorial description of first
order freeness \cite{NicaSpeicher1997Fourier}.
(Note that $\pi^{-1}\gamma_n$ is in this case
the Kreweras complement of $\pi$.)

2) In the second order case we have
\begin{align*}
&(f*g)(1_{m+n},\gamma_{m,n})=\sum_{\pi\in\SNC(m,n)}f(0,\pi)g(0,
\pi^{-1}\gamma_{m,n})\\&\quad+\sum_{\pi\in NC(m)\times NC(n)\atop
\vert\cV\vert=\vert\pi\vert+1} \bigl(f(0,\gamma_{m,n}\pi^{-1})
g(\cV,\pi)+ f(\cV,\pi)g(0,\pi^{-1}\gamma_{m,n})\bigr).
\end{align*}
We should expect that this formula is the combinatorial key for the
understanding of second order freeness. However, in this form it
does not match exactly the formulas appearing in
\cite{HigherOrderFreeness2}. Let us, however, for a multiplicative
function $f$ put, for $\pi\in NC(n)$,
\begin{equation}
\tilde f_1(\pi):=f(1_n,\pi)\qquad (\pi\in NC(n))
\end{equation}
and, for $\pi_1\in NC(m)$ and $\pi_2\in NC(n)$,
\begin{equation}
\tilde f_2(\pi_1,\pi_2)=\!\!\!
\sum_{\substack{\cV\geq\pi_1\times\pi_2, \\
\vert\cV\vert=\vert\pi\vert+1,\\ \cV\vee(\pi_1\times\pi_2)=1_{m+n}}}
\!\!\! f(\cV,\pi_1\times \pi_2).
\end{equation}
Note that in the definition of $\tilde f_2$ the sum is running over all
$\cV$ which connect exactly one cycle of $\pi_1$ with one cycle of
$\pi_2$.

Then, with $h=f*g$, we have
\begin{align*}
&\tilde h_2(1_m,1_n)=\sum_{\pi\in S_{NC}(m,n)}\tilde f_1(\pi)
\tilde g_1(\pi^{-1}\gamma_{m,n})\\
&\qquad+\sum_{\pi_1,\pi_2\in NC(m)\times NC(n)}\bigl( \tilde
f_2(\pi_1,\pi_2)\tilde g_1(\pi_1^{-1}\times\pi_2^{-1}\gamma_{m.n})
\\&\quad\qquad\qquad\qquad\qquad\qquad\qquad
+\tilde f_1(\pi_1\times \pi_2)\tilde
g_2(\pi_1^{-1}\gamma_m,\pi_2^{-1}\gamma_n) \bigr).
\end{align*}
In this form we recover exactly the structure of the formula (10)
from \cite{HigherOrderFreeness2},
which describes second order freeness. The
descriptions in terms of $f$ and in terms of $\tilde f_2$ are
equivalent. Whereas $f$ is multiplicative, $\tilde f_2$ satisfies a
kind of cocycle property. From our present perspective the
description of second (and higher) order freeness in terms of
multiplicative functions seems more natural. In any case, we see
that our convolution of multiplicative functions on partitioned
permutations is a generalization of the structure underlying first
and second order freeness.
\end{example}

\subsection{Zeta and M\"obius function}
In the definition of our convolution we are running over
factorizations of $(\cU,\gamma)$ into products
$(\cV,\pi)\cdot(\cW,\sigma)$. In the first order case the second
factor is determined if the first factor is given. In the general
case, however, we do not have such a uniqueness of the
decomposition; if we fix $(\cV,\pi)$ there might be different
choices for $(\cW,\sigma)$. For example, this situation was
considered in Proposition \ref{pro:types-of-factorizations} in
the case (2b). However, in the case when $(\cW,\sigma)$ is a disc
permutation, it must be of the form
$(0_{\pi^{-1}\gamma},\pi^{-1}\gamma)$ and is thus uniquely
determined. Note that factorizations of such a special form appear
in our formula (\ref{ffkk-leading-compact}) and thus deserve special
attention.

\begin{notation}
Let $(\cU,\gamma)\in\cPP$ be a fixed partitioned permutation. We say
that $(\cV,\pi)\in\cPP$ is \emph{$(\cU,\gamma)$--non-crossing} if
$$ (\cV,\pi)\cdot (0,\pi^{-1}\gamma)= (\cU,\gamma).$$
The set of $(\cU,\gamma)$--non-crossing partitioned permutations will be
denoted by $\cPP_{NC}(\cU,\gamma)$, see Remark
\ref{minimal-planar-remark}.
\end{notation}

To justify this notation we point out that
$(1_n,\gamma_n)$--non-crossing partitioned permutations can be
identified with non-crossing permutations; to be precise
$$\cPP_{NC}(1_n,\gamma_n)=\{(0_\pi,\pi)\mid\pi\in NC(n)\}.$$
Furthermore,
\begin{multline*}
\cPP_{NC}(1_{m+n},\gamma_{m,n})=
\{(0_\pi,\pi)\mid\pi\in\SNC(m,n)\}\cup\\
\cup\{(\cV,\pi_1\times\pi_2)\mid \pi_1\in \NC(m),\pi_2\in NC(n),
\cV\geq\pi, \vert\cV\vert=\vert\pi\vert+1\}.
\end{multline*}
We can now also use a special multiplicative function, which we will
call Zeta-function $\zeta$, to single out such factorizations. It
will be useful to be able to invert formula
(\ref{ffkk-leading-compact}), which means we need also the inverse
of $\zeta$ under our convolution. This inverse, called the
M\"obius-function $\mu$, is a key object in the theory and contains
a lot of important information.

\begin{notation}\
\begin{enumerate} \item The \emph{Zeta-function} $\zeta$ is the multiplicative function on
$\cPP$ which is determined by
$$\zeta(1_n,\pi)=\begin{cases}
1&\text{if $(1_n,\pi)$ is a disc permutation, i.e., if $1_n=0_\pi$,}\\
0&\text{otherwise.}
\end{cases}$$
\item The \emph{M\"obius function} $\mu$ is the inverse of $\zeta$ under
convolution, i.e., it is determined by
$$\zeta*\mu=\delta=\mu*\zeta.$$
\end{enumerate}
\end{notation}

Note that in general
$$\zeta(\cV,\pi)=\begin{cases}
1,& \text{if $\cV=0_\pi$}\\
0,&\text{if $\cV>0_\pi$}.
\end{cases}
$$
It is also quite easy to see that the M\"obius function exists and is uniquely determined
as the inverse of the Zeta-function ---  the determining
equations can be solved recursively. Indeed letting $\mu_n =
\mu(1_n, \gamma_n)$ and $\mu_{m,n} = \mu(1_{m+n},
\gamma_{m,n})$ we have
$$
0 = \mu_{1,1} + \mu_2
$$
$$
0 = \mu_{1,2} + 2 \mu_1 \mu_{1,1} + 
2 \mu_3 + 2 \mu_1 \mu_2 
$$
$$
0 = \mu_{2,2} + 4 \mu_1 \mu_{2,1} 
+ 4 \mu_1^2 \mu_{1,1} + 4 \mu_4 + 8 \mu_1
\mu_3  + 2 \mu_2^2 + 4 \mu_1^2 \mu_2
$$
$$
0 = \mu_{1,3} + 3 \mu_1 \mu_{2,1} + 3
\mu_2 \mu_{1,1} + 3 \mu_4 + 6 \mu_1 \mu_3 +
3 \mu_2^2 + 3 \mu_1^2 \mu_2 
$$
$$
0 = \mu_{2,3} +  2 \mu_1 \mu_{1,3} 
+ 3 \mu_1 \mu_{2,2} + 3 \mu_2 \mu_{1,2} 
+ 9 \mu_1^2 \mu_{1,2} + 6 \mu_1 \mu_2 \mu_{1,1}
+ 6 \mu_1^3 \mu_{1,1}
$$
$$\mbox{}
+ 6 \mu_5 + 18 \mu_1 \mu_4 
+ 12 \mu_2 \mu_3 + 18 \mu_1^2 \mu_3
+12 \mu_1 \mu_2^2  + 6 \mu_1^3 \mu_2
$$
$$
0 =
\mu_{3,3} 
+ 6 \mu_1 \mu_{2,3} 
+ 6 \mu _2 \mu _{1,3}
+ 6 \mu_1^2 \mu_{1,3} 
+ 9 \mu _1^2 \mu_{2,2} 
+ 18 \mu _1 \mu _2 \mu_{1,2}
+ 18  \mu_1^3 \mu_{1,2} 
$$
$$\mbox{}
+ 9 \mu _2^2 \mu_{1,1}
+ 18 \mu_1^2 \mu_2 \mu_{1,1}
+ 9 \mu_1^4 \mu_{1,1}
+ 9 \mu_6 
+ 36 \mu_1 \mu_5
+ 27 \mu_2 \mu_4
+ 54 \mu_1^2 \mu_4
$$
$$\mbox{}
+ 9 \mu_3^2
+ 72 \mu_1 \mu _2 \mu_3
+ 36 \mu_1^3 \mu_3
+ 12 \mu_2^3
+ 36\mu_1^2 \mu_2^2 
+ 9 \mu_1^4 \mu_2    
$$
This shows how, knowing the first order M\"obius function $\mu_n$, the
second order M\"obius function $\mu_{m,n}$ can be calculated recursively. 

One should observe that with these notations we have
$$(f*\zeta)(\cU,\gamma)=\sum_{(\cV,\pi)\in\cPP_{NC}(\cU,\gamma)}f(\cV,\pi).$$

In the following we will use the notation
$$\zeta^{*p}=\underbrace{\zeta*\dots *\zeta}_{\text{$p$-times}}.$$
It is clear, by definition, that $\zeta^{*p}$ counts factorizations into the product of $p$
disc permutations, thus we have the following result.
\begin{proposition}
For $(\cU,\gamma)\in\cPP$ and $p\geq 1$ we have
$$\zeta^{*p}(\cU,\gamma):=\#\{(\pi_1,\dots,\pi_p)\mid
(\cU,\gamma)=(0,\pi_1)\cdots (0,\pi_p)\}.$$
\end{proposition}

Of special interest for us is the case $p=2$.

\begin{proposition}\label{zeta-numbering}
We have for all $r\geq 1$ and $n(1),\dots,n(r)\in\NN$, $n:=n(1)+\cdots+n(r)$ that
$$(\zeta*\zeta)(1_n,\gamma_{n(1),\dots,n(r)})=\#\SNC(n(1),\dots,n(r)).$$
\end{proposition}
\begin{proof}
As noted above, $(\zeta*\zeta)(1_n,\gamma_{n(1),\cdots,n(r)})$
counts the number of factorizations of
$(1_n,\gamma_{n(1),\dots,n(r)})$ into a product of two disc
permutations, i.e., the number of factorizations of the form
$$(1_n,\gamma_{n(1),\dots,n(r)})=(0,\pi)\cdot (0,\pi^{-1}\gamma_{n(1),\dots,
n(r)}),$$ with
$$\vert\pi\vert+\vert \pi^{-1}\gamma\vert=\vert\gamma\vert=n-r$$
and $\pi\vee\gamma=1_n$. But this describes exactly connected
$(n(1),\dots,n(r))$-annular permutations
$\pi\in\SNC(n(1),\dots,n(r))$.
\end{proof}

\begin{notation}
We put
$$c_{n(1),\dots,n(r)}:=\#\SNC(n(1),\dots,n(r)).$$
\end{notation}

Note in particular that $c_n$ counts the number of non-crossing
partitions of $n$ elements and thus is the Catalan number
$$c_n=\frac{1}{n+1} \binom{2n}{n},$$
and that $c_{m,n}$ counts the number of non-crossing $(m,n)$-annular
permutations, and thus \cite{MingoNica2004annular}
$$c_{m,n}=\frac{2mn}{m+n} \binom{2m-1}{m} \binom{2n-1}{n}.$$

More generally, an explicit formula for the number of factorizations into $p$ factors was
derived by Bousquet-M\'elou and Schaeffer \cite{MR1761777}, namely one has (with
$n:=n(1)+\cdots+n(r)$)
$$\zeta^{*p}(1_n,\gamma_{n(1),\cdots,n(r)})=
p\frac{[(p-1)n-1]!}{[(p-1)n-r+2]!}\prod_{i=1}^{r} \left[ n(i){pn(i)-1 \choose n(i)}\right],$$

and thus in particular
$$c_{n(1),\dots,n(r)}=
2\frac{(n-1)!}{(n-r+2)!}\prod_{i=1}^{r}  \left[ n(i){2n(i)-1
\choose n(i)}\right].
$$

For our purposes, however, the following recursive formula for the
number of factorizations is more interesting.

In the next theorem we will show how to reduce the
problem of counting the number of disc
factorizations on $[n]$ to counting the
factorizations on $[n-1]$. This will enable of to
obtain a recursive formula for $c_{n_1, \dots,
n_r}$. 

\begin{notation}\label{factorization-notation}
Let $(\cU, \gamma)$ be a
partitioned permutation of $[n]$ with $\gamma(1)
\not = 1$. Let $\hat \gamma_k$ be the
restriction of $(1,k) \gamma (1, \gamma^{-1}(k))$ to
the invariant subset $[2,n] := \{2,3,4, \dots, n\}$.
Then 
$$
\vert\hat\gamma_k \vert  =
\begin{cases}  
\vert\gamma\vert&\text{if $1$ and $k$ are in
different cycles of $\gamma$,}\\
\vert\gamma\vert-1 & \text{if $k=1$ or $\gamma(1)$ }
\\
\vert\gamma\vert-2 &
\text{if $1$ and $k$ are in the same cycle of $\gamma$}, \\
&\text{but $k \not =1$ and $k \not = \gamma(1)$,}
\\
\end{cases}
$$
Let $\ol \cU = \cU|_{[2,n]}$ be the
restriction of $\cU$ to $[2,n]$, \ie if the blocks
of $\cU$ are $U_1, \dots , U_r$ and $1 \in U_1$, then
the blocks of $\ol \cU$ are $\ol U_1, U_2, \dots,
U_r$ where $\ol U_1 = U_1 \cap [2, n]$. In the
theorem below we sum over a set of partitions
$\cP_k$ of $[2, n]$ described as follows.

For $k=1$,
$\gamma(1)$ or $k$ not in the $\gamma$-orbit of 1,
$\cP_k = \{\, \ol \cU \,\}$ \ie\/ $\cP_k$ consists of
the single partition $\ol \cU$. 

For $k$ in the $\gamma$-orbit of 1 but $k \not =1$, $\gamma(1)$,
$\cP_k = \{ \wh \cU \mid \hat \gamma_k \leq \wh \cU, |\wh \cU| =
|\cU | -2$, and $\ol \cU = \wh \cU \vee (k, \gamma^{-1}(k)) \}$.
In words this means $\ol U_1$ is split into two blocks: 

\begin{list}{$\circ$}{\itemsep=0em\leftmargin=1em}
\item the first containing
the cycle of $\hat \gamma_k$ containing $\gamma^{-1}(k)$ and
some (possibly none) of the other cycles of $\gamma$ contained
in $U_1$
\item the second containing the cycle of $\hat \gamma_k$
containing $k$ and the remaining (possibly none) cycles of
$\gamma$ contained in $U_1$.
\end{list}

More explicitly, in the case $k$ is in the
$\gamma$-orbit of 1 but $k \not =1$, $\gamma(1)$, let
us write $\gamma$ as as a product of cycles $d_1
\cdots d_s$ where $d_1 = (1, \gamma(1), \dots ,
\gamma^t(1))$ is the cycle that contains 1. Let
$d'_1 = (\gamma(1), \gamma^2(1),\ab \dots,
\gamma^{-1}(k))$ and $d''_1 = ( k, \dots ,
\gamma^t(1) )$. Then $\hat
\gamma_k = d'_1 d''_1 d_2 \cdots d_s$. $\cP_k$
consists of all partitions $\wh \cU$ of $[2,n]$ such
that $\wh \cU = \{U'_1, U''_1, U_2,\ab \dots ,
U_r \}$ where $U'_1 \cup U''_1 = \ol U_1$, $U'_1
\cap U''_1 = \emptyset$, $U'_1$ contains $d'_1$,
$U''_1$ contains $d''_1$, and each cycle of
$\gamma$ that was in $U_1$ is now in either $U'_1$
or $U''_1$, \ie\ $\hat\gamma_k \leq \wh \cU$ and $|
\wh \cU| = |\cU| - 2$. 
\end{notation}

\begin{theorem}\label{thm:rec-fact}
\begin{equation}\label{eq:zeta-recursion}
\zeta^{\ast^2}(U, \gamma) =
\sum_{k=1}^n \sum_{\wh \cU \in \cP_k} 
\zeta^{\ast^2}(\wh \cU, \wh{\gamma_k})
\end{equation}
\end{theorem}

\begin{proof}
We must show that for each factorization $(0, \pi)
\cdot (0, \sigma)$ of $(\cU, \gamma)$ there are $k
:= \pi(1)$, $\wh \cU \in \cP_k$, and permutations of
$[2, n]$, $\hat\pi$ and $\hat\sigma$ such that $(0,
\hat\pi) \cdot (0, \hat\sigma) = (\wh \cU ,
\hat\gamma_k)$.
Conversely we must show that given
$k$, $\wh \cU \in \cP_k$ and a factorization $(0,
\hat\pi) \cdot (0,\hat\sigma)$ of $(\wh \cU, \hat
\gamma_k)$  there are $\pi$ and $\sigma$ such that
$(0, \pi) \cdot (0, \sigma) = (\cU, \gamma)$ and
$\pi(1) = k$. Moreover we must show that these two
maps are inverses of each other. The relation between
$\pi$,
$\sigma$ and $\hat \pi$, $\hat \sigma$ is given by
$\hat \pi = (1,k) \pi|_{[2,n]}$, $\hat \sigma =
\sigma (1, \gamma^{-1}(k))|_{[2,n]}$. So on the
level of permutations we have a bijection. The main
work of the proof is to show that starting with
$\pi$ and $\sigma$ we have $\wh \cU := \hat \pi
\vee \hat \sigma \in \cP_k$ and $2 | \wh \cU| -
|\hat \gamma_k| = |\hat \pi| + |\hat \sigma|$; and
then conversely starting with $\wh \cU \in \cP_k$
and a factorization $(0, \hat \pi) \cdot (0, \hat
\sigma)$ of $(\wh \cU, \hat \gamma_k)$ then $2|\cU| -
|\gamma| = |\pi| + |\sigma|$ and $\pi \vee \sigma =
\cU$. 

Note that we have for all $k$
$$
\vert\hat\pi\vert =
\begin{cases}
\vert\pi\vert-1 & k \not =1 \\
|\pi| & k = 1 \\
\end{cases} 
$$ $$
\vert\hat\sigma\vert  =
\begin{cases}
\vert\sigma\vert-1 & k \not =\gamma(1) \\
|\sigma| & k = \gamma(1) \\
\end{cases} 
$$

It is necessary to break the proof into four cases:
$k$ is not in the $\gamma$-orbit of 1; $k$ is in the
$\gamma$-orbit of 1 but $k \not = 1, \gamma(1)$;
$k=1$; and $k = \gamma(1)$.

Suppose we have a factorization
$$(\cU,\gamma)=(0,\pi)\cdot (0,\sigma),$$ i.e.,
$\gamma=\pi\sigma$, $\cU=\pi\vee\sigma$, and
$$2 \vert \cU \vert - \vert \gamma \vert = \vert \pi
\vert + \vert \sigma \vert$$ with $k := \pi(1)$ not
in the $\gamma$-orbit of 1. Then$|\hat\gamma_k| =
|\gamma|$ and $\cP_k$ contains only the partition of
$[2,n]$ which results from
$\cU$ by removing
$1$, \ie\/ $\wh \cU = \ol \cU$. Then we have
$\vert\hat\cU\vert=\vert\cU\vert-1$. Hence $|\hat\pi| + |\hat
\sigma| = |\pi| + |\sigma| - 2 = 2|\cU| - |\gamma| - 2 = |\wh
\cU| - |\gamma| = |\wh \cU| - |\hat \gamma_k |$.

Also $0_\pi|_{[2,n]} = 0_{\hat \pi}$ and
$0_\gamma|_{[2,n]} \leq 0_{\hat \gamma_k}$. Thus
$\wh\cU = (\pi \vee \gamma)|_{[2,n]}
\leq \hat\pi \vee \hat\gamma_k$. On the other hand the
difference between $0_\gamma|_{[2,n]}$ and $0_{\hat\gamma_k}$ is
that the blocks containing 1 and $k$ have been joined. However
these points were already connected by $\pi$. Thus $\hat\pi \vee
\hat\gamma_k \leq \wh \cU$, and so
$\hat\cU=\hat\pi\vee\hat\sigma$, and thus
$$(\hat\cU,\hat\gamma)=(0,\hat\pi)\cdot (0,\hat\sigma).$$

Conversely, given a factorization $(0, \hat\pi) \cdot (0,
\hat\sigma)$ of $(\wh \cU, \hat\gamma_k)$, let $\pi = (1,k)
\hat\pi$ and $\sigma = \hat\sigma (1, \gamma^{-1}(k))$. Then $\pi
\vee \sigma = \cU$ because 1 has been connected to the block of
$\wh\cU$ containing $k$. Also $\#(\pi) = \#(\hat\pi)$ and
$\#(\sigma) = \#(\hat\sigma)$; thus $|\pi| = |\hat\pi|-1$ and
$|\sigma| = |\hat\sigma| - 1$, and so $|\pi| + |\sigma| = 2
|\cU| - |\gamma|$. This establishes the bijection when $k$ is not in
the $\gamma$-orbit of 1.

Let us now consider the case that $1$ and $k$ are in
the same cycle of $\gamma$, but $k \not= 1,
\gamma(1)$. Again suppose that $(0, \pi) \cdot (0,
\sigma)$ is a factorization of $(\cU, \gamma)$ with
$\pi(1) = k$. In this case we have that
$|\hat \gamma_k| = |\gamma|-2$ and so by the triangle
inequality, Lemma 4.7

\begin{align*}
2\vert\hat\pi\vee\hat\sigma\vert-\vert\gamma\vert+2&=
2\vert\hat\pi\vee\hat\sigma\vert-\vert\hat\gamma\vert\\&=
\vert(\hat\pi\vee\hat\sigma,\hat\pi\hat\sigma)\vert\\
&\leq
\vert(0,\hat\pi)\vert+\vert(0,\hat\sigma)\vert\\
&=\vert\hat\pi\vert+\vert\hat\sigma\vert\\
&=\vert\pi\vert+\vert\sigma\vert-2\\
&=2\vert\cU\vert-\vert\gamma\vert-2,
\end{align*}
and thus
$$\vert\hat\pi\vee\hat\sigma\vert\leq\vert\cU\vert-2.$$
On the other hand, let us compare
$$\hat\pi\vee\hat\sigma=\hat\pi\vee
\hat\gamma\qquad\text{with}\qquad \cU=\pi\vee\gamma.$$ Note that all
our changes of the permutations affected only what
happens on the first cycle of $\gamma$. Since the
transition from $\gamma$ to
$\hat\gamma$ consists in removing the point $1$ and splitting the
first cycle of $\gamma$ into two cycles, we can lose at most one
block by going over from $\hat\pi\vee\hat\gamma$ to $\pi\vee\gamma$.
Thus
$$\vert\hat\pi\vee\hat\sigma\vert=(n-1)-\# (\hat\pi\vee\hat\sigma)\geq
(n-1)-(\#\cU+1)=\vert\cU\vert-2,$$ so that we necessarily have the
equality
$$\vert\hat\pi\vee\hat\sigma\vert=\vert\cU\vert-2.$$
Thus $\wh \cU := \hat \pi \vee \hat \sigma \in \cP_k$
and 
$2 |
\hat \pi \vee \hat \sigma| - |\hat \gamma_k | = |\hat
\pi| + |\hat \sigma |$. Hence $(0, \hat \pi) \cdot
(0, \hat \sigma)$ is a factorization of $(\wh \cU,
\hat \gamma_k)$. 

Conversely let us suppose that $k$ is in the
$\gamma$-orbit of 1 but $k \not =1$ or $\gamma(1)$
and $\wh \cU \in \cP_k$ and $(0, \hat \pi) \cdot
(0, \hat \sigma)$ is a factorization of $(\wh \cU,
\hat \gamma_k)$. We must  show that $\pi \vee \sigma
= \cU$ and that $|\pi| + |\sigma| = 2|\cU| -
|\gamma|$. 1 and $k$ are in the same orbit of $\pi$
and 1 and $\gamma^{-1}(k)$ are in the same orbit of
$\sigma$. So the blocks of $\wh \cU$ containing
$d'_1$ and $d''_1$ are joined in $\pi \vee \sigma$.
Thus $\pi \vee \sigma = \cU$. Also $|\wh \cU| =
|\cU| - 2$, so $|\pi| + |\sigma| = |\hat \pi| +
|\hat \sigma| + 2 = 2 |\wh \cU| - |\hat \gamma_k| +2
= |\cU| - |\hat \gamma_k| -2 = 2|\cU| - |\gamma|$.
Thus $(0, \pi) \cdot (0, \sigma)$ is a factorization
of $(\cU, \gamma)$. This establishes the bijection
in the case $k$ is in the $\gamma$-orbit of 1 but $k
\not = 1$ or $\gamma(1)$. 

Next suppose that $k =1$ and $(0, \pi) \cdot (0,
\sigma)$ is a factorization of $(\cU, \gamma)$ with
$\pi(1) = 1$. Then $|\hat \pi| + |\hat \sigma| =
|\pi| + |\sigma| - 1 = 2 |\cU| - |\gamma| -1 = 2
|\wh \cU| - |\gamma| + 1 = 2 |\wh \cU| -
|\hat \gamma_k|$.  Let $U_1$ be the block of $\cU$
containing 1 and $\ol U_1 = U_1 \cap [2,n]$. We must
show that $\ol U_1$ is a block of $\hat \pi \vee
\hat \gamma_k$. Since $\pi \vee \gamma = \cU$ we
know that if $d_i$ and $d_j$ are cycles of $\gamma$
contained in $U_1$ then $\pi$ must connect them.
Since $\pi |_{\ol U_1} = \hat \pi|_{\ol U_1}$ we see
that $\hat \pi$ connects the corresponding cycles of
$\hat \gamma_k$ (which are unchanged except for the
cycle containing 1). Similarly if $f_1$ and $f_2$
are cycles of $\pi$ contained in $U_1$ and neither
is a singleton then they are connected by $\gamma$
and thus by $\hat \gamma_k$. Thus $(0, \hat\pi)
\cdot (0, \hat \sigma)$ is a factorization of $(\wh
\cU, \hat\gamma_k)$. 

Conversely suppose that $k=1$, $\wh \cU \in \cP_k$,
and $(0, \hat \pi) \cdot (0, \hat \sigma)$ is a
factorization of $(\wh \cU, \hat\gamma_k)$. We must
show that $\pi(1) =1$ and $(0, \pi) \cdot (0,
\sigma)$ is a factorization of $(\cU, \gamma)$.
Since $\hat \pi \vee \hat \gamma_k = \wh \cU$ and
$\gamma$ connects 1 to $\gamma(1) \in \ol U_1$, we
have that $\pi \vee \gamma = \cU$. Also $|\pi| +
|\sigma| = |\hat \pi| + |\hat \sigma| + 1= 2 |\wh
\cU| - |\hat \gamma_k| + 1 = 2 |\cU| - |\hat
\gamma_1| -1 = 2 |\cU| - |\gamma|$. Thus $(0, \pi)
\cdot (0, \sigma)$ is a factorization of $(\cU,
\gamma)$. This completes the case when $k =1$. The
proof in the case $k = \gamma(1)$ is exactly the
same except that the roles of $\pi$ and $\sigma$ are
reversed.  \end{proof}

Let us take a closer look at the meaning of Theorem
\ref{thm:rec-fact} for the case $(\cU, \gamma) =
(1_n, \gamma_{n(1), \dots, n(r)})$. To reduce the
depth of subscripts we shall write $c(n_1, \dots ,
n_r)$ for $c_{n(1), \dots, n(r)}$.

\begin{proposition} \label{prop:faktoryzacjena2}
We have for all $r, n_1, \dots, n_r \in \NN$ the
recursion
\begin{multline}\label{eq:rekurencja1}
c(n_1, \dots, n_r) = \sum_{l=2}^r
n_l \cdot c(n_1 + n_l - 1, n_2,
\dots, n_{l-1}, n_{l+1}, \dots, n_r)
\\ + 
\sum_{k=1}^{n_1}
\mathop{\sum_{A = \{i_1, \dots, i_s \}}}_{B = \{
j_1, \dots, j_t \}} 
c(k-1,n_{i_1}, \dots, n_{i_s})
c(n_1 - k, n_{j_1}, \dots , n_{j_t}) \\
\end{multline}
where the sum is over all pairs of subsets $A, B
\subset [2,r]$ such that $A \cap B = \emptyset$ and
$A \cup B = [2,r]$ including the possibility that
either $A$ or $B$ could be empty. We
have for all $m,n\geq 1$
$$ c_n= \sum_{1\leq k\leq n} c_{k-1} c_{n-k},
$$
and
\begin{equation}
\label{eq:rekurencja2} c_{m,n}=  \sum_{1\leq k\leq n} \big( c_{k-1}
c_{m,n-k}+c_{m,k-1} c_{n-k} \big) + mc_{m+n-1}, 
\end{equation} where
we use the convention that $c_0=1$ but $c(n_1,
\dots, n_r) =0$ if $r > 1$ and for some $i$, $n_i
=0$.
\end{proposition}

\begin{proof}
Let $n = n_1 + \cdots + n_r$. By
Proposition \ref{zeta-numbering}
$c(n_1,\dots , n_r) = \zeta^{\ast^2}(1_n,\ab
\gamma_{n_1, \dots , n_r})$. So we must give the
correspondence between the terms on the right hand
side of (\ref{eq:zeta-recursion}) and the right hand
side of (\ref{eq:rekurencja1}). In this case $\cU =
1_n$ and $\ol \cU = 1_{n-1}$ (in the notation of
\ref{factorization-notation}). Thus $\cP_k = \{
1_{n-1}\}$. Also for $n_1 + \cdots + n_{l-1} < k
\leq n_1 + \cdots + n_l$, $\zeta^{\ast^2}(1_{n-1},
\hat \gamma_k) = c(n_1 + n_l -1, n_2, \dots,
n_{l-1}, n_{l+1}, \dots , n_r)$. Thus
\begin{eqnarray}\label{second-part}\lefteqn{
\sum_{k = n_1 +1}^n
\sum_{\wh \cU \in \cP_k} \zeta^{\ast^2}(\wh \cU,
\hat \gamma_k) = \sum_{k = n_1 +1}^n
\zeta^{\ast^2}(1_{n-1},\hat \gamma_k) }\\ \notag
&=&
\sum_{l=2}^r 
\sum_{k = n_1 + \cdots + n_{l-1} + 1}^{n_1 + \cdots +
n_l} c(n_1 + n_l -1, n_2, \dots, n_{l-1}, n_{l+1},
\dots , n_r) \\ \notag
&=&
\sum_{l=2}^r  n_l \cdot
c(n_1 + n_l -1, n_2, \dots, n_{l-1}, n_{l+1},
\dots , n_r) \notag
\end{eqnarray}
For $k \leq n_1$, $\hat \gamma_k = d'_1 d''_1 d_2
\cdots d_r$, with $d'_1$ a cycle of length $k-1$ and
$d''_1$ a cycle of length $n_1 - k$. $\cP_k$ is the
set of all partitions of the cycles of $\hat
\gamma_k$ into two blocks such that $d'_1$ and
$d''_1$ are in different blocks. Hence
$$
\sum_{\wh \cU \in \cP_k}
\zeta^{\ast^2}(\wh \cU, \hat \gamma_k) =
\mathop{\sum_{A = \{ i_1, \dots , i_s \} }}%
_{B = \{j_1, \dots , j_t \}}
c(n_1 - k, n_{i_1}, \dots , n_{i_s})
c(k - 1, n_{j_1}, \dots , n_{j_t})
$$
where the sum is over all pairs of subsets $A, B
\subset [2,r]$ such that $A \cap B = \emptyset$ and
$A \cup B = [2,r]$ including the possibility that
either $A$ or $B$ could be empty. Thus
\begin{multline}\label{first-part}
\sum_{k=1}^{n_1} \sum_{\wh \cU \in \cP_k}
\zeta^{\ast^2}(\wh \cU, \hat \gamma_k) \\
= \sum_{k=1}^{n_1}
\mathop{\sum_{A = \{ i_1, \dots , i_s \} }}%
_{B = \{j_1, \dots , j_t \}}
c(n_1 - k, n_{i_1}, \dots , n_{i_s})
c(k - 1, n_{j_1}, \dots , n_{j_t})
\end{multline}
Assembling equations (\ref{second-part}) and
(\ref{first-part}) gives the result. 
\end{proof}

In
\cite{OBrienZuber84}, O'Brien and Zuber used a similar
formula of this kind in order
to compute the asymptotics of, so called, external field matrix
integral. See also \cite{MR1761777} and Theorem \ref{theo:rekurencja}.

Clearly, our notions around the convolution of functions on $\cPP$
are analogous to (and motivated by) the convolution of functions on
posets. Even though we are not able to put the above theory into the
framework of posets, it seems that this analogy goes quite far.
The following description of the M\"obius functions is an
instance of this---its poset analogue is due to Hall (see
\cite{rota}). It is essentially the
simple observation that one can expand the M\"obius function in
terms of a geometric series as
$$\mu=\zeta^{*-1}=\bigl(\delta+(\zeta-\delta)\bigr)^{*-1}=
\sum_{k=0}^\infty (-1)^k (\zeta-\delta)^{*k}.$$

\begin{proposition}\label{geometric-series}
We have for any $(\cU,\gamma)\in\cPP$ that
$$\mu(\cU,\gamma)=\delta(\cU,\gamma)+
\sum_{k=1}^\infty\sum_{(\cU,\gamma)=(0,\pi_1)\cdots (0,\pi_k)\atop
\pi_i\not=e\;\forall i} (-1)^k.$$
\end{proposition}

\begin{proof}
As noted above this is just the geometric series for
$$(\delta+(\zeta-\delta))^{*-1}.$$ (Note that we are working for
this in the algebra of functions on $\cPP$ with the pointwise sum
and the convolution as sum and product---we are not bothering about
multiplicativity.) The only thing to check is that the sum
is finite, and this is the case because the number of factors
$k$ is bounded by $\vert(\cU,\gamma)\vert$, since
$\vert(0,\pi)\vert\geq 1$ for any $\pi\not= e$.
\end{proof}

This description of the M\"obius function allows us now to derive a
recursive formula for $\mu$.

\begin{theorem}
\label{theo:rekurencja} Consider $(\cU,\gamma)\in\cPP$ such that
$\gamma(1)\not=1$. Then we have
\begin{equation} \label{eq:moebous}
\mu(\cU,\gamma)= (-1) \sum_{(0,(1,k))\cdot(\cV,\pi)=(\cU,\gamma)\atop
k\not=1} \mu(\cV,\pi),
\end{equation}
where the sum runs over all decompositions of $(\cU,\gamma)$ into a product of a disc
transposition $(0,(1,k))$ (with $k\geq 2$) and a $(\cV,\pi)\in \cPP$.
\end{theorem}

The proof of this theorem will rely on the following lemma.

\begin{lemma}
\label{lem:faktoryzacje} Let $(\cU,\gamma)\in\cPP$ such that
$\gamma(1)\not=1$. For $p\in\NN$, we denote by $\cS_p$ the set
consisting of all tuples $(\pi_1,\dots,\pi_p)$ of permutations such
that $\pi_i\not=e$ for all $i=1,\dots,p$ and
$$(0,\pi_1) \cdots (0,\pi_p)=(\cU,\gamma).$$

We consider now the two sums
\begin{equation}
\label{eq:Asuma1} S_1:=\sum_{p=1}^\infty
\sum_{(\pi_1,\dots,\pi_p)\in\cS_p} (-1)^p
\end{equation}
and
\begin{equation}
\label{eq:Asuma2} S_2:=\sum_{p=1}^\infty
\sum_{(\pi_1,\dots,\pi_p)\in\cS_p\atop \text{$\pi_1=(1,k)$ for
$k\not=1$}} (-1)^p
\end{equation}
where the second sum $S_2$ is over all tuples $(\pi_1,\dots,\pi_p)$
as for the first sum $S_1$, but now with the additional property
that $\pi_1$ is a transposition interchanging the element $1$ with
some other element.

Then the two sums \eqref{eq:Asuma1} and \eqref{eq:Asuma2} are equal,
$$S_1=S_2.$$
\end{lemma}

\begin{proof}
Let $\pi=(\pi_1,\dots,\pi_p) \in \cS_p$. Let $1\leq q\leq p$
denote the smallest index for which $1$ is not a fixed point of 
$\pi_q$; note that such a $q$ necessarily exists since
$\gamma(1) \not = 1$. We shall group all factorizations into
three classes: 1a), 1b) and 2). Class 1) consists of
factorizations for which
$\pi_{q}$ is a transposition interchanging $1$ with some other
element. The subclass 1a) consists of factorizations for which
$q=1$ and subclass 1b) of those for which $q\geq 2$. Class 2)
consists of all other factorizations.

Let $\Pi=(\pi_1,\dots,\pi_p)$ be  a factorization from the class
1b). We define
$$\Pi'=(\pi'_1,\dots,\pi'_{p-1})=(\pi_1,\dots,\pi_{q-2},\pi_{q-1}
\pi_{q},\pi_{q+1},\dots,\pi_p).$$ In the following we shall prove
that $f:\Pi\mapsto \Pi'$ is a bijection between factorizations of
class $1b)$ and factorizations of class $2)$.

Firstly, we prove that $\Pi' \in \cS_p$ and is of class $2)$.
Clearly,
$\pi'_{q-1}=\pi_{q-1} \pi_q$ is a permutation which does not fix
$1$, it is not a transposition interchanging $1$ with some other
element, and we have
$$(0,\pi_{q-1})\cdot (0,\pi_{q})=(0,\pi'_{q-1}).$$

In order to show that $f$ is a bijection we shall describe its
inverse. If $\Pi'=(\pi'_1,\dots,\pi_{p-1}') \in \cS_p$ and is of
class 2), we define $1\leq q\leq p-1$ to be the smallest number
for which
$\pi'_{q-1}$ does not fix $1$. There is a unique decomposition
$\pi'_{q-1}=\pi_{q-1}
\pi_q$ such that $1$ is a fixed point of $\pi_{q-1}$ and $\pi_q$ is
a transposition interchanging $1$ with some other element. 
Thus $|\pi_{q-1}| + |\pi_q| = |\pi'_{q-1}|$. The
assumption that the factorization $\Pi'$ is of class 2) implies that
$\pi_{q-1}\neq e$. For $1\leq i\leq q-2$ we set $\pi_i=\pi'_i$ and
for $q+1\leq i \leq p$ we set $\pi_i=\pi'_{i-1}$. In this way we
defined $\Pi=(\pi_1,\dots,\pi_p)$. Now it is easy to check that
$g:\Pi'\mapsto \Pi$ is a left and right inverse of $f$.

Since the factorization $\Pi$ and the corresponding $\Pi'$
contribute to \eqref{eq:Asuma1} with the opposite signs, the
contribution of all factorizations of class 1b) cancels with the
contribution of factorizations of class 2).
\end{proof}

\begin{proof}[Proof of \ref{theo:rekurencja}] In the proof we
will consider all factorizations $(0, \pi_1) \cdot (0, \pi_2)
\cdots (0, \pi_p) = (\cU, \gamma)$ with the requirement that
$\pi_i \not = e$ for all $i$, \ie\/ $(\pi_1, \dots, \pi_p) \in
\cS_p$, as in the proof of Lemma \ref{lem:faktoryzacje}.
Sometimes we will require in addition that
$\pi_1 = (1, k)$ with $k \not = 1$. To simplify the notation we
will not explicitly state every time that $\pi_i \not = e$. 
Since $\gamma(1) \not = 1$ we have $\delta(\cU, \gamma) = 0$.
When $\gamma$ is a transposition the right hand side of equation
(\ref{eq:moebous}) is $-1$; so we can assume that $\gamma$ is not
a transposition. So by Proposition
\ref{geometric-series} we have
\begin{eqnarray*}
\mu(\cU, \gamma) &=&
\sum_{p =1}^\infty 
\mathop{\sum_{(0, \pi_1) \cdots (0, \pi_p)}}%
_{= (\cU, \gamma)} (-1)^p 
\stackrel{(\ref{lem:faktoryzacje})}{=}
\sum_{p =1}^\infty 
\mathop{\sum_{(0, (1,k)) \cdots (0, \pi_p)}}%
_{\mbox{} = (\cU, \gamma)} (-1)^p \\
&=&
\sum_{p=2}^\infty
\mathop{\sum_{(0, (1,k)), (\cV, \pi)}}%
_{(0, (1,k))\cdot (\cV, \pi) = (\cU, \gamma)}
\mathop{\sum_{(0, \pi_2) \cdots (0, \pi_p)}}%
_{\mbox{} = (\cV, \pi)} (-1)^p \\
&=&
- \kern-2em\mathop{\sum_{(0, (1,k)), (\cV, \pi)}}%
_{(0, (1,k))\cdot (\cV, \pi) = (\cU, \gamma)}
\sum_{p=2}^\infty
\mathop{\sum_{(0, \pi_2) \cdots (0, \pi_p)}}%
_{= (\cV, \pi)} (-1)^{p-1} \\
&=&
- \kern-2em\mathop{\sum_{(0, (1,k)), (\cV, \pi)}}%
_{(0, (1,k))\cdot (\cV, \pi) = (\cU, \gamma)}
\sum_{p=2}^\infty
\mathop{\sum_{(0, \pi_2) \cdots (0, \pi_p)}}%
_{= (\cV, \pi)} (-1)^{p-1} \\
&=&
- \kern-2em\mathop{\sum_{(0, (1,k)), (\cV, \pi)}}%
_{(0, (1,k))\cdot (\cV, \pi) = (\cU, \gamma)}
\mu( \cV, \pi)
\end{eqnarray*}
\end{proof}

One observes that the recursion formulas for the M\"obius function
and for $\zeta^{*2}$ look very similar. However, there are some
significant differences. The recursion for $\zeta^{*2}$
effectively expresses $\zeta^{*2}$ for $n$ points in terms of
$\zeta^{*2}$ for
$n-1$ points. The recursion for the M\"obius function does not
reduce the number of points. Nevertheless, at least for first and
second order one can match the two recursions and connect the values
of the M\"obius function with the values of the
function $\zeta^{*2}$ (i.e., with the number of non-crossing
partitions and non-crossing annular permutations). In order to see
this let us first specify the meaning of Theorem
\ref{theo:rekurencja} for first and second order. In first order we
get
$$
\mu(1_n,\gamma_n) = - \sum_{1\leq k\leq n-1} \mu(1_{k},\gamma_k)
\mu(1_{n-k},\gamma_{n-k}),
$$
which shows that $(-1)^{n} \mu(1_{n+1},\gamma_{n+1})$ and
$\zeta^{*2}(1_n,\gamma_n)$ satisfy the same recursion (namely the
one for the Catalan numbers). This is, of course, just the
well-known fact
\cite{Kreweras,Speicher1994}
that the M\"obius function on
non-crossing partitions is given by the signed and shifted Catalan
numbers. In second order our recursion reads
\begin{align*}
&(-1)\mu(1_{m+n},\gamma_{m,n})= m \cdot\mu(1_{m+n},\gamma_{m+n})\\
&\quad+\sum_{1\leq k\leq n-1} \big( \mu(1_{m+k},\gamma_{m,k})
\mu(1_{n-k},\gamma_{n-k}) + \mu(1_{m+n-k},\gamma_{m,n-k}) \mu(1_k,
\gamma_{k}) \big),
\end{align*}
which we recognize --- by taking into account the shifted
relation between $\mu$ and $\zeta^{*2}$ on the first level --- as
the recursion for $(-1)^{m+n} \zeta^{*2}(1_{m+n}, \ab
\gamma_{m,n})$. Let us collect these explicit results about the
M\"obius function in the following theorem.

\begin{theorem}
We have for $m,n\in\NN$ that
$$ \mu(1_n,\gamma_n)= (-1)^{n-1}\cdot \# NC(n-1)=(-1)^{n-1}\cdot c_{n-1}$$
and
$$\mu(1_{m+n},\gamma_{m,n})=(-1)^{m+n}\cdot \#\SNC(m,n)=(-1)^{m+n}\cdot c_{m,n}.$$
\end{theorem}

For higher orders we were not able to match the values of $\mu$ with
those of $\zeta^{*2}$.

\section{R-transform formulas}\label{sec:R-transform}
Let us consider the situation that two multiplicative functions $f$ and $h$ on $\cPP$ are
related by $h=f*\zeta$. We want to understand what this means for the relations between the
numbers $\kk_n:=f(1_n,\gamma_n)$ and $\kk_{m,n}:=f(1_{m+n},\gamma_{m,n})$ on one side and the
numbers $\alpha_n:=h(1_n,\gamma_n)$ and $\alpha_{m,n}:=h(1_{m+n},\gamma_{m,n})$ on the other
side. In particular, we want to express this in terms of the generating power series of these
numbers,
$$C(x):=1+\sum_{n\geq 1} \kk_nx^n,\qquad
C(x,y):=\sum_{m,n\geq 1}\kk_{m,n}x^my^n$$ and
$$M(x):=1+\sum_{n\geq 1} \alpha_nx^n,\qquad
M(x,y):=\sum_{m,n\geq 1}\alpha_{m,n}x^my^n.$$ (Note that the above summation corresponds to
putting formally
$$f(1_0,\gamma_0):=1\qquad\text{and}\qquad
f(1_0,\gamma_{0,0}):=0$$ for a  multiplicative $f$. Our notation is motivated by the fact that
the most important realization of the relation $h=f*\zeta$ will be the situation where the
$\alpha$'s are the correlation moments and the $\kk$'s the corresponding cumulants, thus $M$ is
a moment series and $C$ is a cumulant series.) On the first order level we have
$$\alpha_n=\sum_{\pi\in NC(n)}f(0_\pi,\pi),$$
which is the usual moment-cumulant formula of free probability theory, and it is well-known
\cite{Speicher1994} that this is equivalent to
$$C\bigl(xM(x)\bigr)=M(x).$$

Our main goal now is to derive the analogue of this for the second order level. There we have
$$\alpha_{m,n}=\sum_{\pi\in\SNC(m,n)}f(0_\pi,\pi)+
\sum_{\pi_1\times\pi_2\in NC(m)\times NC(n) \atop \vert\cV\vert=\vert\pi_1\times\pi_2\vert+1}
f(\cV,\pi_1\times\pi_2).$$ 

It turns out that the second term, the sum over disconnected
partitions, is quite easy to deal with. The first term, the sum
over connected annular permutations, looks much more involved,
however, one can handle this also if one realizes that one can
reduce this first term to the second one. Namely, one can sum
over all connected annular permutations by first bundling all
through-cycles into one through-cycle and secondly decomposing
this through-cycle into sub-cycles all of which are
through-cycles. In this way one can reduce the problem of
dealing with all annular non-crossing permutations to the
problem of considering permutations with exactly one
through-cycle and the problem of considering permutations where
all cycles are through-cycles. The first problem corresponds
exactly to the above sum over disconnected partitions. So we can
write

$$\sum_{\pi\in\SNC(m,n)}f(0_\pi,\pi)=\sum_{\pi_1\times\pi_2\in
NC(m)\times NC(n)
\atop \vert\cV\vert=\vert\pi_1\times\pi_2\vert+1} \tilde f(\cV,\pi_1\times\pi_2),$$ where
$\tilde f$ is now the multiplicative function corresponding to
$$\tilde f(1_n,\gamma_n)=\tilde
\kk_n,\qquad \tilde f(1_{m+n},\gamma_{m,n})=\tilde \kk_{m,n}$$ with
$$\tilde\kk_n:=\kk_n$$
and
$$\tilde \kk_{m,n}:=\sum_{\pi\in\SNC^{all}(m,n)} f(0_\pi,\pi).$$
Thus we can combine this to get finally
\begin{align*}
\alpha_{m,n}&=\sum_{\pi_1\times\pi_2\in NC(m)\times NC(n) \atop
\vert\cV\vert=\vert\pi_1\times\pi_2\vert+1} \bigl(f(\cV,\pi_1\times\pi_2)+ \tilde
f(\cV,\pi_1\times\pi_2)\bigr)\\&= \sum_{\pi_1\times\pi_2\in NC(m)\times NC(n) \atop
\vert\cV\vert=\vert\pi_1\times\pi_2\vert+1} g(\cV,\pi_1\times\pi_2),
\end{align*}
where $g$ is the multiplicative function corresponding to
$$g(1_n,\gamma_n)=\tilde
\alpha_n,\qquad g(1_{m+n},\gamma_{m,n})=\tilde \alpha_{m,n}$$ with
$$\tilde\alpha_n=\tilde \kk_n=\kk_n$$
and
$$\tilde \alpha_{m,n}=\kk_{m,n}+\tilde\kk_{m,n}.$$
So we have to translate the relation between $\tilde \kk_{m,n}$ and $f$ and the relation
between $\alpha_{m,n}$ and $g$ into relations between the corresponding formal power series.

\begin{proposition}\label{prop:first-step}
Let $f$ be a multiplicative function on $\cPP$ with
$$f(1_n,\gamma_n)=:\kk_n\qquad
\text{and}\qquad C(x):=1+\sum_{n\geq 1}\kk_nx^n.$$ Put
$$\tilde \kk_{m,n}:=\sum_{\pi\in\SNC^{all}(m,n)} f(0_\pi,\pi),$$
where $\SNC^{all}(m,n)$ denotes the permutations in $\SNC(m,n)$ for which all cycles are
through-cycles. Consider the corresponding generating power series
$$\tilde C(x,y):=\sum_{m,n\geq 1} \tilde\kk_{m,n}x^my^n.$$
Then we have
$$\tilde C(x,y)=
-xy\frac{\partial^2}{\partial x\partial y} \log\bigl( \frac{xC(y)-yC(x)}{x-y}),$$ or
equivalently
$$\tilde C(x,y)=
-xy\Bigl(\frac{\bigl(C(x)-xC'(x)\bigr)\bigl(C(y)-yC'(y)\bigr)}
{\bigl(xC(y)-yC(x)\bigr)^2}-\frac 1{(x-y)^2}\Bigr).$$
\end{proposition}

\begin{proof}
Note that we can parametrize an element $\pi\in \SNC^{all}(m,n)$ in a bijective way by
specifying the number of cycles, the number of elements on each circle for all cycles, the
position of a fixed element (let's say 1) in its cycle and the first element on the other
circle of this cycle. Let us denote the number of cycles by $r$, the number of elements of the
cycles on the first circle by $i_1,\dots,i_r$ and the number of elements of those cycles on the
other circle by $j_1,\dots,j_r$. Thus the $l$-th cycle contains $i_l+j_l$ elements and makes
the contribution $\kk_{i_l+j_l}$ in the calculation of $\tilde\kk_{m,n}$. We normalize things
so that the first cycle contains the element 1. Fixing $i_1,\dots,i_r$ and $j_1,\dots,j_r$ we
thus have $i_1$ possibilities for where 1 sits in the first cycle and $n$ possibilities for the
first element of this cycle on the other circle. This means we have
$$\tilde \kk_{m,n}=\sum_{r\geq 1}\sum_{i_1,\dots,i_r\geq 1\atop
i_1+\cdots+ i_r=m}\sum_{j_1,\dots,j_r\geq 1\atop j_1+\cdots+j_r=n} i_1 n \kk_{i_1+j_1}\cdots
\kk_{i_r+j_r}$$ and thus
\begin{align*}
&\tilde C(x,y)=\sum_{r\geq 1}\sum_{i_1,\dots,i_r\geq 1} \sum_{j_1,\dots,j_r\geq 1} i_1
(j_1+\cdots j_r) \kk_{i_1+j_1}\cdots
\kk_{i_r+j_r}x^{i_1}\cdots x^{i_r}y^{j_1}\cdots y^{j_r}\\
&=\sum_{r\geq 1}\sum_{i_1,\dots,i_r\geq 1}\sum_{j_1,\dots,j_r\geq 1} i_1 y\frac
\partial{\partial y}\bigl( \kk_{i_1gb+j_1}\cdots
\kk_{i_r+j_r}x^{i_1}\cdots x^{i_r}y^{j_1}\cdots y^{j_r}\bigr)\\
&=\sum_{r\geq 1}y\frac \partial{\partial y}\Bigl(\bigl( \sum_{i_1,j_1\geq 1}
i_1\kk_{i_1+j_1}x^{i_1}y^{j_1}\big)\cdot \bigl(\sum_{i_2,j_2\geq
1}\kk_{i_2+j_2}x^{i_2}y^{j_2}\bigr)\cdots \bigl(\sum_{i_r,j_r\geq
1}\kk_{i_r+j_r}x^{i_r}y^{j_r}\bigr)\Bigr)
\end{align*}
Let us now use the notation
$$
\hat C(x,y):=\sum_{i,j\geq 1}\kk_{i+j}x^{i}y^{j}.
$$
Then we can continue with
\begin{align*}
\tilde C(x,y) &=\sum_{r\geq 1}y\frac \partial{\partial y}\Bigl(
\bigl(x\frac\partial{\partial x} \hat C(x,y)\bigr)\cdot \hat C(x,y)^{r-1}\Bigr)\\
&=\sum_{r\geq 1}xy\frac \partial{\partial y}\Bigl(
\frac 1r\frac\partial{\partial x} \bigl(\hat C(x,y)^r\bigr)\Bigr)\\
&=xy\frac \partial{\partial y}\frac\partial{\partial x} \Bigl(\sum_{r\geq 1}\frac 1r
\hat C(x,y)^r\Bigr)\\
&=-xy\frac \partial{\partial y}\frac\partial{\partial x} \log\bigl(1-\hat C(x,y)\bigr)
\end{align*}
The assertions follow now by noting that
$$\hat C(x,y)=1-\frac{xC(y)-yC(x)}{x-y}$$
and by working out the partial derivatives.
\end{proof}

\begin{proposition}\label{prop:second-step}
Let $g$ be a multiplicative function on $\cPP$. Put
$$\tilde \alpha_{m,n}:=g(1_{m+n},\gamma_{m,n})$$
and denote its generating power series of second order by
$$H(x,y):=\sum_{m,n\geq 1}\tilde \alpha_{m,n} x^my^n.$$
Put
$$\alpha_n:=(g*\zeta)(1_n,\gamma_n)$$
and
$$\alpha_{m,n}:=\sum_{(\cV,\pi_1\times\pi_2)\atop
\vert \cV\vert=\vert \pi_1\times\pi_2\vert+1} g(\cV,\pi)$$ and denote the corresponding
generating functions by
$$M(x):=1+\sum_{n\geq 1} \alpha_n x^n\qquad\text{and}\qquad
M(x,y):=\sum_{m,n\geq 1} \alpha_{m,n}x^my^n.$$ Then we have the relation
$$M(x,y)=H(xM(x),yM(y))\cdot
\Bigl(1+ x\frac {M'(x)}{M(x)}\Bigr)\cdot \Bigl(1+ y\frac {M'(y)}{M(y)}\Bigr).$$
\end{proposition}

\begin{proof}
Let us do the summation in the definition of $\alpha_{m,n}$ in the way that we first fix the
two cycles $V_1\in\pi_1$ and $V_2\in\pi_2$ which are connected by $\cV$ and sum over all
possibilities for fixed $V_1,V_2$. If $V_1$ has $k$ elements and $V_2$ has $l$ elements then
this contributes the factor $\tilde \alpha_{k,l}$. Furthermore, $\pi_1\backslash V_1$
decomposes into $k$ independent non-crossing partitions and the summations over them (for fixed
$V_1$) gives the $\alpha_i$ for the intervals between consecutive elements from $V_1$. (Of
course, we are counting here modulo $m$.) For the final summation over $V_1$ we have to notice
that there are two different possibilities: either a fixed number (let's say 1) is an element
of $V_1$ - in which case we can specify the situation by prescribing the number $k$ of elements
of $V_1$ and the differences $i_1,\dots,i_k$ between consecutive elements in $V_1$ - or 1 is
not an element of $V_1$, --- in which case we need an extra
factor $i_1$, because we have now
$i_1$ different possibilities how 1 can lie between two consecutive elements of $V_1$. Since we
have the same situation for $V_2$ we can thus write $\alpha_{m,n}$ in the
form
\begin{align*}
\alpha_{m,n}&=\sum_{k,l\geq 1}\sum_{i_1,\dots,i_k\geq 0\atop k+i_1+\cdots+
i_k=m}\sum_{j_1,\dots,j_l\geq 0\atop l+j_1+\cdots+j_l=n}
\tilde\alpha_{k,l}\alpha_{i_1}\cdots\alpha_{i_k}\alpha_{j_1}\cdots\alpha_{j_l}
\Bigl(1+i_1+j_1+i_1 j_1\Bigr).
\end{align*}
Translating this into generating power series gives the assertion.
\end{proof}

The combination of the previous two propositions, with
$$H(x,y)=C(x,y)+\tilde C(x,y),$$
gives now our main result.

\begin{theorem}\label{second-R}
Let $f$ and $h$ be multiplicative functions on $\cPP$ which are related by
$$h=f*\zeta.$$
Denote
$$\kk_n:=f(1_n,\gamma_n),\qquad \kk_{m,n}:=f(1_{m+n},\gamma_{m,n})$$
and
$$\alpha_n:=h(1_n,\gamma_n),\qquad \alpha_{m,n}:=h(1_{m+n},\gamma_{m,n})$$
and define the corresponding generating power series
$$C(x):=1+\sum_{n\geq 1} \kk_nx^n,\qquad
C(x,y):=\sum_{m,n\geq 1}\kk_{m,n}x^my^n$$ and
$$M(x):=1+\sum_{n\geq 1} \alpha_nx^n,\qquad
M(x,y):=\sum_{m,n\geq 1}\alpha_{m,n}x^my^n.$$ Then we have as formal power series the first
order relation
\begin{equation}\label{eq:M-C-first-order}
C(xM(x))=M(x)
\end{equation}
and for the second order
\begin{equation}\label{eq:M-C-eins}
M(x,y)=H\bigl(xM(x),yM(y)\bigr)\cdot \frac{\frac d{dx}(xM(x))}{M(x)}\cdot \frac{ \frac
d{dy}(yM(y))}{M(y)},
\end{equation}
where
\begin{equation}\label{eq:M-C-einsa}
H(x,y):=C(x,y)-xy\frac{\partial^2}{\partial x\partial
y}\log\Bigl(\frac{xC(y)-yC(x)}{x-y}\Bigr),
\end{equation}
or equivalently,
\begin{multline}\label{eq:M-C}
M(x,y)=C\bigl(xM(x),yM(y)\bigr)\cdot \frac{\frac d{dx}(xM(x))}{M(x)}\cdot \frac{
\frac d{dy}(yM(y))}{M(y)}\\
+xy\Bigl( \frac{\frac d{dx}(xM(x))\cdot \frac d{dy}(yM(y))}{(xM(x)-yM(y))^2}-\frac 1{(x-y)^2}
\Bigr).\end{multline}
\end{theorem}

\begin{proof}
The formulation (\ref{eq:M-C-eins}) and (\ref{eq:M-C-einsa}) follows directly from a
combination of Propositions  \ref{prop:first-step} and \ref{prop:second-step}. In order to
reformulate this to (\ref{eq:M-C}) one uses the equivalence of the two formulas in Proposition
\ref{prop:first-step} and the fact that $C(xM(x))=M(x)$ yields
$$1-x C'(xM(x))=\frac{M(x)}{\frac d{dx}(xM(x))}.$$
\end{proof}
If we go over from the moment generating series $M$ to a kind of Cauchy transform like
quantity $G$, then these formulas take on a particularly nice
form.

\begin{corollary}
Consider the same situation and notations as in Theorem 
\ref{second-R}. In terms of 
$$G(x):=\frac 1x M(1/x),\quad G(x,y):=\frac 1{xy} M(1/x,1/y),
\quad \cR(x,y):=\frac 1{xy}C(x,y)$$ the Equation (\ref{eq:M-C})
can be written as
\begin{equation}
G(x,y)=G'(x)G'(y)\Bigl\{\cR(G(x),G(y))+\frac
1{(G(x)-G(y))^2}\Bigl\}-\frac 1{(x-y)^2}.
\end{equation}
\end{corollary}

$\cR(x,y)$ is the second order $R$-transform. Note that
Voiculescu's first order
$R$-transform $\cR$ is defined by the relation $C(x)=1+z\cR(x)$, and equation
(\ref{eq:M-C-first-order}) says for this
$$\frac 1{G(x)}+\cR(G(x))=x,$$
i.e., that $G(x)$ and $K(x):=\frac 1x+\cR(x)$ are inverses of each other under composition.

\begin{example}
Let us apply our formulas to some examples.

1) If we put $f$ to be the multiplicative function with $\kk_2=1$ and all other $\kk_n$ and all
$\kk_{m,n}$ vanishing, then $h=f*\zeta$ counts the non-crossing pairings, i.e., in this case
$M(x)$ is the generating function of the number of non-crossing pairings (on one circle) and
$M(x,y)$ is the generating function of the number of non-crossing annular pairings (on two
circles). Let us calculate it by using the above theorem.
\\
We have
$$C(x)=1+x^2, \qquad C(x,y)=0$$
and we know that $M$ is the generating function of number of non-crossing pairings on a circle.
In this case
$$\hat C(x,y)=xy,$$
and thus
$$
H(x,y)=-xy \frac{\partial^2}{\partial x\partial y}\log(1-xy) = \frac{xy}{(1-xy)^2},
$$
which yields the result
$$M(x,y)=xy\cdot\frac  {\frac d{dx}(xM(x))\cdot \frac d{dy}(yM(y))}
{\bigl(1-xy M(x)M(y)\bigr)^2}.$$
Related formulas are known in the physical literature,
see, e.g, \cite{fmp}, \cite{bz}, \cite{kkp}.

2) If we put $f=\zeta$ then $h=\zeta*\zeta$ counts the non-crossing permutations, i.e., in this
case $M$ is the generating function of the number of non-crossing permutations (which is the
same as non-crossing partition) on one circle and $M_2$ is the generating function of the
number of annular non-crossing permutations (on two circles).
\\
We have
$$C(x)=\frac 1{1-x},\qquad C(x,y)=0.$$
In this case
$$\hat C(x,y)=\frac{1-x-y}{(1-x)(1-y)},$$
and thus
$$H(x,y)=-xy \frac{\partial^2}{\partial x\partial y}\log(1-xy)
=\frac{xy}{(1-x-y)^2},$$ which yields
$$M(x,y)=xy\cdot\frac  {\frac d{dx}(xM(x))\cdot \frac d{dy}(yM(y))}
{\bigl(1- x M(x)-y M(y)\bigr)^2}.$$

3) Let us finally see whether we can extract the value of the M\"obius function from our
formula. Since we have $\delta=\mu*\zeta$, our formula with
$$M(x)=1+x,\qquad M(x,y)=0$$
should allow to solve for $C(x,y)$ which is then the generating function for the annular
M\"obius function. Note that we already know $M(x)$ in this case to be the generating function
of the disc M\"obius function.
\\
If $M(x,y)$ vanishes identically this implies that $H(x,y)$ vanishes identically, leading to
the identity
\begin{align*}
C(x,y)&= xy\frac{\partial^2}{\partial x\partial y}\log(\frac{xC(y)-yC(x)}{x-y})\\
&=xy\Bigl(\frac{(C(x)-xC'(x))\cdot (C(y)-yC'(y))}{(xC(y)-yC(x))^2}-\frac 1{(x-y)^2}\Bigr)
\end{align*}
\end{example}

\begin{remark}
Equation \ref{eq:M-C} gives the second order version of {\em
moment-cumulant} relations. 
$$
\alpha_{1,1} = \kappa_{1,1} + \kappa_2
$$
$$
\alpha_{2,1} = \kappa_{1,2} + 2 \kappa_1 \kappa_{1,1} + 
2 \kappa_3 + 2 \kappa_1 \kappa_2 
$$
$$
\alpha_{2,2} = \kappa_{2,2} + 4 \kappa_1 \kappa_{2,1} 
+ 4 \kappa_1^2 \kappa_{1,1} + 4 \kappa_4 + 8 \kappa_1
\kappa_3  + 2 \kappa_2^2 + 4 \kappa_1^2 \kappa_2
$$
$$
\alpha_{1,3} = \kappa_{1,3} + 3 \kappa_1 \kappa_{2,1} + 3
\kappa_2 \kappa_{1,1} + 3 \kappa_4 + 6 \kappa_1 \kappa_3 +
3 \kappa_2^2 + 3 \kappa_1^2 \kappa_2 
$$
$$
\alpha_{2,3} = \kappa_{2,3} +  2 \kappa_1 \kappa_{1,3} 
+ 3 \kappa_1 \kappa_{2,2} + 3 \kappa_2 \kappa_{1,2} 
+ 9 \kappa_1^2 \kappa_{1,2} + 6 \kappa_1 \kappa_2 \kappa_{1,1}
+ 6 \kappa_1^3 \kappa_{1,1}
$$
$$\mbox{}
+ 6 \kappa_5 + 18 \kappa_1 \kappa_4 
+ 12 \kappa_2 \kappa_3 + 18 \kappa_1^2 \kappa_3
+12 \kappa_1 \kappa_2^2  + 6 \kappa_1^3 \kappa_2
$$
$$
\alpha_{3,3} =
\kappa_{3,3} 
+ 6 \kappa_1 \kappa_{2,3} 
+ 6 \kappa _2 \kappa _{1,3}
+ 6 \kappa_1^2 \kappa_{1,3} 
+ 9 \kappa _1^2 \kappa_{2,2} 
+ 18 \kappa _1 \kappa _2 \kappa_{1,2}
+ 18  \kappa_1^3 \kappa_{1,2} 
$$
$$\mbox{}
+ 9 \kappa _2^2 \kappa_{1,1}
+ 18 \kappa_1^2 \kappa_2 \kappa_{1,1}
+ 9 \kappa_1^4 \kappa_{1,1}
+ 9 \kappa_6 
+ 36 \kappa_1 \kappa_5
+ 27 \kappa_2 \kappa_4
+ 54 \kappa_1^2 \kappa_4
$$
$$\mbox{}
+ 9 \kappa_3^2
+ 72 \kappa_1 \kappa _2 \kappa_3
+ 36 \kappa_1^3 \kappa_3
+ 12 \kappa_2^3
+ 36\kappa_1^2 \kappa_2^2 
+ 9 \kappa_1^4 \kappa_2    
$$

\medskip\hrule\medskip

$$
\kappa_{1, 1} = \alpha_{1}^2 - \alpha_{2} + \alpha_{1, 1}
$$
$$
\kappa_{1, 2} = -4 \alpha_{1}^3 + 6 \alpha_{1} \alpha_{2} 
- 2 \alpha_{3} - 2 \alpha_{1} \alpha_{1, 1} + \alpha_{1, 2}
$$
$$
\kappa_{2, 2} = 
18 \alpha_{1}^4 - 36 \alpha_{1}^2 \alpha_{2} 
+ 6 \alpha_{2}^2 + 16 \alpha_{1} \alpha_{3} - 4 \alpha_{4} 
+ 4 \alpha_{1}^2 \alpha_{1, 1} 
- 4 \alpha_{1} \alpha_{1, 2} + \alpha_{2, 2}
$$
$$
\kappa_{1, 3} = 
15 \alpha_{1}^4 - 30 \alpha_{1}^2 \alpha_{2} 
+ 6 \alpha_{2}^2 + 12 \alpha_{1} \alpha_{3} - 3 \alpha_{4} 
+ 6 \alpha_{1}^2 \alpha_{1, 1} - 3 \alpha_{2} \alpha_{1, 1} 
- 3 \alpha_{1} \alpha_{1, 2} + \alpha_{1, 3}
$$
$$
\kappa_{2, 3} = -72 \alpha_{1}^5 
+ 180\ \alpha_{1}^3\ \alpha_{2} - 72 \alpha_{1} \alpha_{2}^2 
- 84\ \alpha_{1}^2 \alpha_{3} + 24 \alpha_{2} \alpha_{3} 
+ 30 \alpha_{1} \alpha_{4} - 6 \alpha_{5} 
$$
$$\mbox{}
- 12 \alpha_{1}^3 \alpha_{1, 1} 
+ 6 \alpha_{1} \alpha_{2} \alpha_{1, 1} 
+ 12 \alpha_{1}^2 \alpha_{1, 2} - 3 \alpha_{2} \alpha_{1, 2} 
- 2 \alpha_{1} \alpha_{1, 3} - 3 \alpha_{1} \alpha_{2, 2} + \alpha_{2, 3}
$$
$$
\kappa_{3, 3} = 
300 \alpha_{1}^6 - 900 \alpha_{1}^4 \alpha_{2} 
+ 576 \alpha_{1}^2 \alpha_{2}^2 - 48 \alpha_{2}^3 
+ 432 \alpha_{1}^3 \alpha_{3} 
- 288 \alpha_{1} \alpha_{2} \alpha_{3} + 18 \alpha_{3}^2 
$$
$$\mbox{}
- 180 \alpha_{1}^2 \alpha_{4} + 45 \alpha_{2} \alpha_{4} 
+ 54 \alpha_{1} \alpha_{5} - 9 \alpha_{6} 
+ 36 \alpha_{1}^4 \alpha_{1, 1} 
- 36 \alpha_{1}^2 \alpha_{2} \alpha_{1, 1} 
+ 9 \alpha_{2}^2 \alpha_{1, 1} 
$$
$$\mbox{}
- 36 \alpha_{1}^3 \alpha_{1, 2} 
+ 18 \alpha_{1} \alpha_{2} \alpha_{1, 2} 
+ 12 \alpha_{1}^2 \alpha_{1, 3} 
- 6 \alpha_{2} \alpha_{1, 3} 
+ 9 \alpha_{1}^2 \alpha_{2, 2} 
- 6 \alpha_{1} \alpha_{2, 3} + \alpha_{3, 3}
$$
\end{remark}

\section{Higher order freeness and corresponding cumulants}
\label{sec:higher}
\subsection{Abstract framework}
\begin{definition}
A \emph{higher-order (non-commutative) probability space}, or
briefly HOPS, $(\cA,\ff)$ consists of a unital algebra
$\cA$ and a collection $\ff=(\ff_n)_{n\in\NN}$ of maps
($n\in\NN$)
$$\ff_n:\underbrace{\cA\times\dots\times\cA}_{\text{$n$ times}}
\to\CC,$$ which are linear and tracial in each of its $n$ arguments
and which are symmetric under exchange of its $n$ arguments and
which satisfy
$$\ff_1(1)=1$$
and
$$\ff_n(1,a_2,\dots,a_n)=0$$
for all $n\geq 2$ and all $a_2,\dots,a_n\in\cA$.
\end{definition}

Of course, we can include the usual (first order) non-commutative
probability space $(\cA,\ff_1)$ into this framework by putting
all higher $\ff_n$ equal to zero. In the same way we recover a
second order non-commutative probability space
$(\cA,\ff_1,\ff_2)$ by putting $\ff_n=0$ for all $n\geq 3$.

\begin{definition}
1) We denote by $\cPP(\cA)$ the set of partitioned permutations
decorated with elements from $\cA$, i.e.,
$$\cPP=\bigcup_{n\in\NN}\bigl(\cPP(n)\times\cA^n\bigr).$$
2) For a function
\begin{align*}
f:\cPP(\cA)&\to\CC\\
(\cV,\pi)\times (a_1,\dots,a_n)&\mapsto f(\cV,\pi)[a_1,\dots,a_n]
\end{align*}
and a function
$$g:\cPP\to\CC$$
we define their convolution
$$f*g:\cPP(\cA)\to\CC$$
by $$(f*g)(\cU,\gamma)[a_1,\dots,a_n]:=
\sum_{(\cV,\pi),(\cW,\sigma)\in\cPP(n)\atop (\cV,\pi)\cdot
(\cW,\sigma)=(\cU,\gamma)} f(\cV,\pi)[a_1,\dots,a_n]\cdot
g(\cW,\sigma)$$ for all $(\cU,\gamma)\in\cPP(n)$ and all
$a_1,\dots,a_n\in\cA$.
\end{definition}

\begin{definition}
A function $f:\cPP(\cA)\to\CC$ is called \emph{multiplicative} if we
have
$$f(\cV,\pi)[a_1,\dots,a_n]=\prod_{B\in\cV} f(1_B,\pi\vert_B)
[(a_1,\dots,a_n)_B]$$ and
$$f(1_n,\sigma^{-1}\pi\sigma)[a_{\sigma(1)},\dots,a_{\sigma(n)}]=
f(1_n,\pi)[a_1,\dots,a_n]$$ for all $a_1,\dots,a_n\in\cA$ and all
$\pi,\sigma\in S(n)$.
\end{definition}

Note that this extension of our formalism on multiplicative
functions on $\cPP$ and their convolution from the last section is
not changing the results from the last section. The structure of all
formulas remains the same; one just has to insert the
$a_1,\dots,a_n$ as dummy variables at the right positions. Thus, in
particular, $\delta$ is still the unit for this extended convolution
and $f=g*\zeta$ is equivalent to $g=f*\mu$ for multiplicative $f,g$
on $\cPP(\cA)$. And again, the convolution of a multiplicative
function on $\cPP(\cA)$ with a multiplicative function on $\cPP$
gives a multiplicative function on $\cPP(\cA)$.

It is clear that a multiplicative function $f$ on $\cPP(\cA)$ is
uniquely determined by the values of
$f(1_{n},\gamma_{n(1),\dots,n(r)})[a_1,\dots,a_{n}]$ (where we put
$n:=n(1)+\dots+n(r)$) for all $r\in\NN$, all $n(1),\dots,n(r)\in\NN$
and all $a_1,\dots,a_n\in \cA$.

\subsection{Moment and cumulant functions}
Let us now apply this formalism to get moment and cumulant functions for higher order
probability spaces. So let a HOPS $(\cA,\ff)$ be given. We will use the $\ff_n$ to produce a
multiplicative ``moment" function on $\cPP(\cA)$, which we will also denote by $\ff$. Namely,
we put
\begin{multline*}
\ff(1_{n},\gamma_{n(1),\dots,n(r)})[a_1,\dots,a_{n}]\\:=
\ff_r(a_1\cdots a_{n(1)};\dots;a_{n(1)+\cdots n(r-1)+1}\cdots a_{n})
\end{multline*}
and extend this by multiplicativity. (Note that we need the $\ff_n$
to be tracial in their arguments for this extension.)

Here is an example for our function $\ff$.
\begin{align*}
\ff\bigl(\{1,3,4\}\{2\},(1,3)(2)(4)\bigr)[a_1,a_2,a_3,a_4]=
\ff_2(a_1a_3,a_4)\cdot \ff_1(a_2)
\end{align*}

\begin{definition}\label{def:freecumulants}
For a given HOPS $(\cA,\ff)$ we define the corresponding
\emph{(higher order) free cumulants} as a function on $\cPP(\cA)$ by
$$\kk=\ff*\mu,$$
or more explicitly
$$\kk(\cU,\gamma)[a_1,\dots,a_n]:=\sum_{(\cV,\pi),(\cW,\sigma)\in
\cPP(n)\atop (\cV,\pi)\cdot(\cW,\sigma)=(\cU,\gamma)}
\ff(\cV,\pi)[a_1,\dots,a_n]\cdot\mu(\cW,\sigma),$$ for all $n\in\NN$, $(\cU,\gamma)\in\cPP(n)$,
$a_1,\dots,a_n\in\cA$.
\end{definition}
As we noted before the definition above is equivalent to
the statement $\ff=\kk*\zeta$, i.e.,
$$\ff(\cU,\gamma)[a_1,\dots,a_n]=\sum_{(\cV,\pi)\in \cPP_{NC}(\cU,\gamma)}
\kk(\cV,\pi)[a_1,\dots,a_n]$$ for all
$(\cU,\gamma)[a_1,\dots,a_n]\in\cPP(\cA)$.

Furthermore, as with $\ff$, $\kk$ is also a multiplicative
function on
$\cPP(\cA)$. Thus in the same way as all $\ff(\cU,\gamma)$ are
determined by the knowledge of all
\begin{multline*}
\ff(1_{n},\gamma_{n(1),\dots,n(r)})[a_1,\dots,a_n]
\\=\ff_r(a_1 \cdots a_{n(1)};\dots;
a_{n(1)+\cdots + n(r-1)+1}\cdots a_{n(1) + \cdots + n(r)})
\end{multline*}
the free cumulants $\kk(\cU,\gamma)$ are determined by the values of
\begin{multline*}
\kk(1_{n},\gamma_{n(1),\dots,n(r)})[a_1,\dots,a_n]
\\=:
\kk_{n(1),\dots,n(r)}(a_1, \dots ,a_{n(1)};\dots; 
a_{n(1) + \cdots + n(r-1)+1},\dots ,a_{n(1) + \cdots + n(r)}).
\end{multline*}

\begin{remark}
Note that whereas on the level of $\ff$ we also know (by definition)
that we can multiply elements along the cycles of $\pi$ (and thus we
do not need a comma as separator for those elements along a cycle),
this is not true for $\kk$. Thus we have, e.g.,
$$\ff(1_3,(1,2)(3))[a_1,a_2,a_3]=\ff_2(a_1a_2;a_3)=
\ff(1_2,(1),(2))[a_1a_2;a_3],$$ but no clear relation exists
among
$$\kk(1_3,(1,2)(3))[a_1,a_2,a_3]=\kk_{2,1}(a_1,a_2;a_3)$$
and
$$\kk(1_2,(1),(2))[a_1a_2;a_3]=\kk_{1,1}[a_1a_2;a_3].$$
\end{remark}

Note also that since our convolution on $\cPP$ coincides on the
first level with the usual convolution of multiplicative functions
on non-crossing partitions, the above definition of cumulants
reduces on the first level to the usual free cumulants.

\subsection{Higher order freeness}
Equipped with the notion of cumulants we can now define ``freeness"
by the requirement of vanishing of mixed cumulants.
\begin{definition}\label{def:higherfreeness}
We say that a family $(\cX_i)_{i\in I}$ of subsets of $\cA$ is \emph{free (of all orders)} if
we have the following vanishing of mixed cumulants: For all $n\geq 2$ and all $a_k\in
\cX_{i(k)}$ ($1\leq k\leq n$) such that $i(p)\not= i(q)$ for some $1\leq p,q\leq n$ we have
$$\kk(1_n,\pi)[a_1,\dots,a_n]=0$$
for all $\pi\in S(n)$.
\end{definition}

\begin{example}
Let us see that this definition includes the definition of Voiculescu
\cite{VoiculescuDykemaNica} for (first order) freeness and the definition of Mingo and Speicher
\cite{HigherOrderFreeness1} for second order freeness.

1) On the first level this follows from the fact that our cumulants
reduce then to the usual free cumulants and it is well-known that
freeness is equivalent to the vanishing of mixed cumulants. One can
see it directly as follows: Let us consider $a_k\in\cX_{i(k)}$ with
$i(k)\not=i(k+1)$ and $\ff_1(a_k)=0$ for all $k=1,\dots,n$. Then we
have
$$\ff_1(a_1\cdots a_n)=\ff(1_n,\gamma_n)[a_1,\dots,a_n]
=\sum_{\pi\in NC(n)}\kk(0_\pi,\pi)[a_1,\dots,a_n].$$ However the vanishing of mixed moments
means now that the only $\pi$ which contribute are those which do not connect elements from
different sets. Furthermore, the fact that all our variables are centered excludes singletons.
But then it is easy to see that there are no such $\pi$ at all, so the sum is zero.

2) Now we have to consider two cyclically alternating and centered
tuples $a_1,\dots,a_m$ and $b_1,\dots,b_n$. Then we have
\begin{align*}
\ff_2(a_1\cdots a_m;b_1\cdots b_n)&=
\ff(1_{m+n},\gamma_{m,n})[a_1,\dots,a_m,b_1,\dots,b_n]\\
&=\sum_{(\cV,\pi)\in\cPP_{NC}(m,n)}\kk(\cV,\pi)[a_1,\dots,a_m,
b_1,\dots,b_n].
\end{align*}
Again, the vanishing of mixed moments requires that $(\cV,\pi)$ connects only elements from the
same set and the centredness of the elements excludes singletons. It is then easy to see that,
for $n\geq 2$, the only possibilities for such $(\cV,\pi)$ arise for $m=n$ and they have to be
disc permutations $(0_\pi,\pi)$ which are pairings $(a_1,b_{1+s})(a_2,b_{2+s})\cdots
(a_n,b_{n+s})$ for some $s$. The factors $k(1_2,(...))[a_k,b_{k+s}]$ are just
$\ff_1(a_2b_{2+s})$, so that one finally gets, for $n\geq 2$, the formula
\begin{align*}
\ff_2(a_1\cdots a_m;b_1\cdots b_n)&= \delta_{mn}\sum_{k=1}^n
\ff_1(a_1b_{1+s})\cdots \ff_1(a_nb_{n+s}).
\end{align*}
For $n=m=1$ one gets with
$$\ff_2(a_1;b_1)=k_2(a_1,a_2)+k_{1,1}(a_1;b_1)$$
the conclusion that $\ff_2(a_1;b_1)$ has to vanish if $a_1$ and
$b_1$ are from different sets. Nothing is required if both are from
the same set. We see that we get exactly the defining properties for
second order freeness from \cite{HigherOrderFreeness1}.

3) It would be nice to be able to reformulate in a similar way the
definition of higher order freeness in terms of the $\ff$ instead of
the cumulants. However, the situation with more than two circles is
getting much more involved and we are not aware of such a
reformulation for third and higher order freeness.
\end{example}

As in the case of the first order freeness, one sees immediately
that constants are free from everything.

\begin{proposition}\label{prop:one-free-all}
Let $(\cA,\ff)$ be a HOPS. Then $\{1\}$ is free of all orders from every subset
$\cX\subset\cA$.
\end{proposition}

\begin{proof}
We have to prove that
$$\kk(1_n,\gamma_{n(1),\dots,n(r)})[1,a_2,\dots,a_n]=0,$$
unless $n=1$. We will do this by induction on $n$. The case $n=2$ is
clear because
$$\kk(1_2,(12))[1,a_2]=\ff_1(1\cdot a_2)-\ff_1(1)\cdot \ff_1(a_2)=0$$
and
$$\kk(1_2,(1)(2))[1,a_2]= \ff_2(1;a_2)=0.$$
In general, one has
\begin{align*}
\ff(1_n,\gamma_{n(1),\dots,n(r)})&[1,a_2,\dots,a_n]\\&=
\sum_{(\cV,\pi)\in\cPP_{NC}(n(1),\dots,n(r))}\kk(\cV,\pi)[1,a_2,\dots,a_n]\\
&=\kk(1_n,\gamma_{n(1),\dots,n(r)})[1,a_2,\dots,a_n]\\
&\qquad\qquad+ \sum_{(\cV,\pi)\in\cPP_{NC}(n(1),\dots,n(r))\atop
\vert(\cV,\pi)\vert<\vert(1_n,\gamma_{n(1),\dots,n(r)})
\vert}\kk(\cV,\pi)[1,a_2,\dots,a_n]
\end{align*}
By induction hypothesis, in the later sum exactly terms of the form
$(\{1\}\cup\tilde\cV,(1)\cup\tilde \pi)$ with
$$(\tilde\cV,\tilde\pi)\in\cPP_{NC}(n(1)-1,n(2),\dots,n(r))$$
contribute. In the case $n(1)>1$ the sum over those yields
$$\ff(1_{n-1},\gamma_{n(1)-1,n(2),\dots,n(r)}[a_2,\dots,a_n].$$
In this case, also
$$\ff(1_n,\gamma_{n(1),\dots,n(r)})[1,a_2,\dots,a_n]=
\ff(1_{n-1},\gamma_{n(1)-1,n(2),\dots,n(r)}[a_2,\dots,a_n],$$ and
thus $\kk(1_n,\gamma_{n(1),\dots,n(r)})[1,a_2,\dots,a_n]=0$. If, on
the other side, $n(1)=1$ (i.e., $1$ is the only element on its
circle), then we have to set
$$\cPP_{NC}(0,n(2),\dots,n(r))=\emptyset,$$
because then the first circle cannot be connected to the others if
we ask $1$ to be a cycle of its own. But this means that in this
case
$$\kk(1_n,\gamma_{n(1),\dots,n(r)})[1,a_2,\dots,a_n]=
\ff(1_n,\gamma_{n(1),\dots,n(r)})[1,a_2,\dots,a_n]$$ However, for
$n(1)=1$ and $n>1$ we have
$$\ff(1_n,\gamma_{1,\dots,n(r)})[1,a_2,\dots,a_n]=0.$$
\end{proof}

Note that our definition of freeness behaves clearly very nicely
with respect to decompositions of our sets. For example, we have
that $\cX_1,\cX_2,\cX_3$ are free if and only if $\cX_1$ and
$\cX_2\cup\cX_3$ are free and $\cX_2$ and $\cX_3$ are free. Thus we
can reduce the investigation of freeness to the understanding of
freeness for the case of two sets. A characterization for this is
given in the next theorem.

\begin{theorem}\label{thm:charact-freeness}
Let $(\cA,\ff)$ be a higher order probability space and consider two
subsets of $\cX_1,\cX_2\subset\cA$. Then the following are
equivalent.
\begin{enumerate}
\item
The sets $\cX_1,\cX_2$ are free of all orders.
\item
The sets $\cX_1\cup\{1\}$, $\cX_2\cup\{1\}$ are free of all orders.
\item
We have
\begin{multline*}
\qquad\qquad\ff(\cU,\gamma)[a_1b_1,\dots,a_nb_n]\\=
\sum_{(\cV,\pi)\cdot (\cW,\sigma)=(\cU,\gamma)}
\kk(\cV,\pi)[a_1,\dots,a_n]\cdot \ff(\cW,\sigma)[b_1,\dots,b_n]
\end{multline*}
for all $n\in\NN$, all $(\cU,\gamma)\in\cPP(n)$ and all
$a_1,\dots,a_n\in\cX_1\cup\{1\}$, $b_1,\dots,b_n\in\cX_2\cup\{1\}$.
\item
We have
\begin{multline*}
\qquad\qquad\ff(\cU,\gamma)[a_1b_1,\dots,a_nb_n]\\=
\sum_{(\cV,\pi)\cdot (\cW,\sigma)=(\cU,\gamma)}
\ff(\cV,\pi)[a_1,\dots,a_n]\cdot \kk(\cW,\sigma)[b_1,\dots,b_n]
\end{multline*}
for all $n\in\NN$, all $(\cU,\gamma)\in\cPP(n)$ and all
$a_1,\dots,a_n\in\cX_1\cup\{1\}$, $b_1,\dots,b_n\in\cX_2\cup\{1\}$.
\item
We have
\begin{multline*}
\qquad\qquad\kk(\cU,\gamma)[a_1b_1,\dots,a_nb_n]\\=
\sum_{(\cV,\pi)\cdot (\cW,\sigma)=(\cU,\gamma)}
\kk(\cV,\pi)[a_1,\dots,a_n]\cdot \kk(\cW,\sigma)[b_1,\dots,b_n]
\end{multline*}
for all $n\in\NN$, all $(\cU,\gamma)\in\cPP(n)$ and all
$a_1,\dots,a_n\in\cX_1\cup\{1\}$, $b_1,\dots,b_n\in\cX_2\cup\{1\}$.

\end{enumerate}
\end{theorem}

In order to prove this we would like to write
$\ff(\cU,\gamma)[a_1b_1,\dots,a_nb_n]$ in the form
$\ff(\hat\cU,\hat\gamma) [a_1,b_1,\dots,a_n,b_n]$. Let us
introduce the following formalism for this. Let
$(\cU,\gamma)\in\cPP(n)$ be a partitioned permutation of the
numbers $1,2,3,\dots,n$. Double now this set of numbers by
introducing a copy $\bar 1,\bar 2,\bar 3,\dots,\bar n$ and
interleave the new and old numbers as follows:
$$1, \bar 1, 2,\bar 2, 3, \bar 3, \dots ,n,\bar n.$$
If we induce now $(\cU,\gamma)$ on $1,2,\dots,n$ to
$(\hat\cU,\hat\gamma)$ on $1,\bar1, \dots, n,\bar n$ by putting
$$\hat\gamma(k)=\bar k\qquad\text{and}\qquad
\hat\gamma(\bar k)=\gamma(k),$$ then this has exactly the wanted property. The vanishing of
mixed cumulants means that in the factorizations of $(\hat\cU,\hat\gamma)$ in $(\cV,\pi)$ times
a disc permutation we are only interested in $(\cV,\pi)$ which have the property that each
block of $\cV$ contains either only unbarred numbers or only bared numbers, i.e., $(\cV,\pi)$
must be of the form $(\cV_a\cup\cV_b,\pi_a\cup\pi_b)$ with
$$(\cV_a,\pi_a)\in\cPP(1,\dots,n)\qquad\text{and}\qquad
(\cV_b,\pi_b)\in\cPP(\bar 1,\dots,\bar n).$$

Let us first observe some simple relations between the quantities on
$1,\dots,n$ and their relatives on $1,\bar 1,\dots,n,\bar n$.

\begin{lemma}
\label{lem-corresponding-length}
1) We have
$$\vert\hat\gamma\vert=n+\vert\gamma\vert,\qquad
\vert\hat\cU\vert=n+\vert\cU\vert,$$ and thus
$$\vert(\hat\cU,\hat\gamma)\vert=n+\vert(\cU,\gamma)\vert.$$
2) We have
$$\vert\pi_a\cup\pi_b\vert=\vert\pi_a\vert+\vert\pi_b\vert,\qquad
\vert\cV_a\cup\cV_b\vert=\vert\cV_a\vert+\vert\cV_b\vert,$$ and thus
$$\vert(\cV_a\cup\cV_b,\pi_a\cup\pi_b)\vert=\vert(\cV_a,\pi_a)\vert
+\vert(\cV_b,\pi_b)\vert.$$ 3) We have that
$(\pi_a\cup\pi_b)\hat\gamma$ maps unbarred to bared and bared to
unbarred elements and, for all $k=1,\dots,n$,
$$[(\pi_a\cup\pi_b)\hat\gamma]^2(\bar k)=\pi_b\pi_a\gamma (k),$$
thus
$$\vert(\pi_a\cup\pi_b)\hat\gamma\vert=n+
\vert\pi_b\pi_a\gamma\vert$$
\end{lemma}

\begin{proof}
Only the third part is non-trivial. To see this observe
$$(\pi_a\cup\pi_b)\hat\gamma(\bar k)=\pi_a(\gamma(k))$$
and thus
$$[(\pi_a\cup\pi_b)\hat\gamma]^2(\bar k)=\pi_b \overline{
\pi_a(\gamma(k))},$$ which is our first equation, with the
identification of $\pi_b\in S(\bar1,\dots,\bar 1)$ with the
corresponding permutation in $S(1,\dots,n)$. Since the mapping
between bared and unbarred elements is clear, this yields that
$(\pi_a\cup\pi_b)\hat\gamma$ and $\pi_b\pi_a \gamma$ have the same
number of orbits which gives the last equation.
\end{proof}

This lemma allows us to characterize the contributing factorizations
in $(\hat\cU,\hat\gamma)$ in terms of special factorizations of
$(\cU,\gamma)$.

\begin{proposition}
\label{prop:special-factorizations}
The statement
$$(\cV_a\cup\cV_b,\pi_a\cup\pi_b)\in\cPP_{NC}(\hat\cU,\hat\gamma)$$
is equivalent to the statement
$$(\cV_a,\pi_a)\cdot (\cV_b,\pi_b)\in\cPP_{NC}(\cU,\gamma),$$
where in the last product we identify $(\cV_b,\pi)\in\cPP(\bar 1,
\dots,\bar n)$ with the corresponding element in $\cPP(1,\dots,n)$.
\end{proposition}

\begin{proof}
Note that
$(\cV_a\cup\cV_b,\pi_a\cup\pi_b)\in\cPP_{NC}(\hat\cU,\hat\gamma)$ is
equivalent to
\begin{equation}\label{eq:equiv1}
\vert(\cV_a\cup\cV_b,\pi_a\cup\pi_b)\vert+\vert(\pi_a\cup\pi_b)^{-1}
\hat\gamma\vert=\vert(\hat\cU,\hat\gamma)\vert
\end{equation}
and
\begin{equation}\label{eq:equiv3}
\hat\cU=(\cV_a\cup\cV_b)\vee\hat\gamma.
\end{equation}
On the other hand, $(\cV_a,\pi_a)\cdot
(\cV_b,\pi_b)\in\cPP_{NC}(\cU,\gamma)$, means
$$(\cV_a,\pi_a)\cdot (\cV_b,\pi_b)\cdot
(0_{\pi_b^{-1} \pi_a^{-1}\gamma},\pi_b^{-1}\pi_a^{-1}\gamma)=
(\cU,\gamma),$$ which is equivalent to
\begin{equation}\label{eq:equiv2}
\vert(\cV_a,\pi_a)\vert+\vert(\cV_b,\pi_b)\vert+ \vert
\pi_b^{-1}\pi_a^{-1}\gamma\vert=\vert(\cU,\gamma)\vert
\end{equation}
and
\begin{equation}\label{eq:equiv4}
\cU=\cV_a\vee\cV_b\vee\gamma.
\end{equation}
Equations (\ref{eq:equiv1}) and (\ref{eq:equiv2}) are, by  Lemma
\ref{lem-corresponding-length}, equivalent.

The equivalence between (\ref{eq:equiv3}) and (\ref{eq:equiv4}) is also easily checked.
\end{proof}

Equipped with these tools we can now prove our main Theorem
\ref{thm:charact-freeness}.

\begin{proof}
The equivalences between (3), (4), and (5) follow by convolving with
the $\zeta$ or the $\mu$ function. That (2) is actually the same as
(1) follows from Prop. \ref{prop:one-free-all}.

$(1)\Longrightarrow (3)$: We have
\begin{align*}
\ff(\cU,\gamma)&[a_1b_1,\dots,a_nb_n]=
\ff(\hat\cU,\hat\gamma)[a_1,b_1,\dots,a_n,b_n]\\
&= \sum_{(\cV,\pi)\cdot(\cW,\sigma)=(\hat\cU,\hat\gamma)}
\kk(\cV,\pi)[a_1,b_1,\dots,a_n,b_n]\cdot
\zeta(\cW,\sigma)\\
&= \sum_{(\cV,\pi)\in\cPP_{NC}(\hat\cU,\hat\gamma)}
\kk(\cV,\pi)[a_1,b_1,\dots,a_n,b_n]
\end{align*}
By our assumption on the vanishing of mixed cumulants, only
$(\cV,\pi)$ of the form $(\cV_1\cup\cV_2,\pi_a\cup\pi_b)$ with
$$(\cV_a,\pi_a)\in\cPP(1,\dots,n)\qquad\text{and}\qquad
(\cV_b,\pi_b)\in\cPP(\bar 1,\dots,\bar n)$$ contribute and, by the
above Proposition
\ref{prop:special-factorizations},
$$(\cV_a\cup\cV_b,\pi_a\cup\pi_b)\in\cPP_{NC}(\hat\cU,\hat\gamma)$$
is equivalent to
$$(\cV_a,\pi_a)\cdot (\cV_b,\pi_b)\in\cPP_{NC}(\cU,\gamma).$$
Thus we can continue with
\begin{align*}
&\ff(\cU,\gamma)[a_1b_1,\dots,a_nb_n]\\
&=\sum_{(\cV_a\cup\cV_b,\pi_a\cup\pi_b)\in\cPP_{NC}(\hat\cU,\hat\gamma)}
\kk(\cV_a,\pi_a)[a_1,a_2,\dots,a_n]\cdot
\kk(\cV_b,\pi_b)[b_1,b_2,\dots,b_n]\\
&=\sum_{(\cV_a,\pi_a)\cdot (\cV_b,\pi_b)\in\cPP_{NC}(\cU,\gamma)}
\kk(\cV_a,\pi_a)[a_1,a_2,\dots,a_n]\cdot
\kk(\cV_b,\pi_b)[b_1,b_2,\dots,b_n]\\
&=\sum_{(\cV_a,\pi_a)\cdot
(\cV_b,\pi_b)\cdot(\cW,\sigma)=(\cU,\gamma)}
\kk(\cV_a,\pi_a)[a_1,a_2,\dots,a_n]\cdot\\
&\qquad\qquad\qquad\qquad\qquad\qquad\qquad\cdot
\kk(\cV_b,\pi_b)[b_1,b_2,\dots,b_n]\cdot \zeta(\cW,\sigma)\\
&=\sum_{(\cV_a,\pi_a)\cdot (\cV,\pi)=(\cU,\gamma)}
\kk(\cV_a,\pi_a)[a_1,\dots,a_n]\cdot \ff(\cV,\pi)[b_1,\dots,b_n].
\end{align*}
$(3)\Longrightarrow (1)$: Note that (3) allows us to calculate all
moments of elements from $\cX_1\cup\cX_2$ out of the moments of
elements from $\cX_1$ and the moments of elements from $\cX_2$. (In
order to do so, we also have to allow some of the $a$'s or $b$'s to
be equal to the unit 1.) Since this calculation rule is the same as
for free sets, this shows that the sets $\cX_1$ and $\cX_2$ must be
free.
\end{proof}

This theorem is now the key ingredient to transfer freeness from
sets to their generated algebras.

\begin{theorem}
Let $(\cA,\ff)$ be a HOPS and consider subsets $(\cX_i)_{i\in I}$.
For each $i\in I$, let $\cA_i$ be the unital algebra generated by
elements from $\cX_i$. Then the following are equivalent.
\begin{enumerate}
\item
The subsets $(\cX_i)_{i\in I}$ are free of all orders.
\item
The subalgebras $(\cA_i)_{i\in I}$ are free of all orders.
\end{enumerate}
\end{theorem}

\begin{proof}
Since the cumulant $\kk(\cV,\pi)[a_1,\dots,a_n]$ is a multi-linear
functional in the $n$ variables $a_1,\dots,a_n$, it is clear that
taking sums of elements within the sets $\cX_i$ preserves freeness.
What we have to see is that also taking products preserves freeness.
Since we can iterate our arguments, it suffices to see the
following: if $\cX_1$ and $\cX_2$ are free, then also
$\cX_1\cup\{a_0a_1\mid a_0,a_1\in\cX_1\}$ and $\cX_2$ are free.
Adding one product after the other to $\cX_1$ and by Theorem
\ref{thm:charact-freeness} it is enough to show that
\begin{multline*}
\ff(\cU,\gamma)[a_0a_1b_1,a_2b_2\dots,a_nb_n]\\=
\sum_{(\cV,\pi)\cdot (\cW,\sigma)=(\cU,\gamma)}
\ff(\cV,\pi)[a_0a_1,a_2\dots,a_n]\cdot
\kk(\cW,\sigma)[b_1,\dots,b_n]
\end{multline*}
for all $n\in\NN$, all $(\cU,\gamma)\in\cPP(n)$ and all
$a_0,a_1,\dots,a_n\in\cX_1\cup\{1\}$,
$b_1,\dots,b_n\in\cX_2\cup\{1\}$. Let us induce
$(\cU,\pi)\in\cPP(1,\dots,n)$ to
$(\hat\cU,\hat\pi)\in\cPP(0,1,\dots,n)$ by requiring that $\hat\cW$
and $\hat \pi$ restricted to $1,\dots,n$ agree with $\cW$ and $\pi$,
respectively, and that $0$ and $1$ are in the same block of
$\hat\cW$ and $\hat\pi(0)=1$. Then we can calculate
\begin{align*}
&\ff(\cU,\gamma)[a_0a_1b_1,a_2b_2\dots,a_nb_n]=
\ff(\hat\cU,\hat\pi)[a_0 1,a_1b_1,a_2b_2,\dots,a_nb_n]\\
&= \sum_{(\cV,\pi)\cdot (\cW,\sigma)=(\hat\cU,\hat\gamma)}
\ff(\cV,\pi)[a_0,a_1,a_2\dots,a_n]\cdot
\kk(\cW,\sigma)[1,b_1,\dots,b_n].
\end{align*}
By Proposition \ref{prop:one-free-all} we know that
$\kk(\cW,\sigma)[1,b_1,\dots,b_n]$ is only different from zero if
$\cW$ has $0$ as a singleton, i.e., $(\cW,\sigma)$ has to be of the
form
$$\cW=\{0\}\cup \tilde\cW,\qquad \sigma=(0)\tilde\sigma,$$
with
$$(\tilde\cW,\tilde\sigma)\in\cPP(1,\dots,n).$$
But then we must have that $\pi(0)=1$ and $0$ and $1$ must be in the
same block of $\cV$. Thus there is a unique $(\cV',\pi')$ so that
$(\cV,\pi)=(\hat\cV',\hat\pi')$ and
$$(\cV,\pi)\cdot (\cW,\sigma)=(\hat\cU,\hat\gamma)$$ is equivalent to
$$(\cV',\pi')\cdot (\tilde\cW,\tilde\sigma)=(\cU,\gamma).$$
Note also that in this situation
$$\kk(\cW,\sigma)[a_0,a_1,,\dots,a_n]=
\kk(\tilde\cW,\tilde\sigma)[b_1,b_2,\dots,b_n].$$ and
$$\ff(\hat\cV',\hat\pi')[a_0,a_1,,\dots,a_n]=
\ff(\cV',\pi')[a_0a_1,a_2,\dots,a_n].$$

So we can continue the above calculation as follows
\begin{align*}
&\ff(\cU,\gamma)[a_0a_1b_1,a_2b_2\dots,a_nb_n]\\
&= \sum_{(\cV',\pi')\cdot (\tilde\cW,\tilde\sigma)=(\cU,\gamma)}
\ff(\cV',\pi')[a_0a_1,a_2\dots,a_n]\cdot
\kk(\tilde\cW,\tilde\sigma)[b_1,\dots,b_n],
\end{align*}
which is exactly what we had to show.
\end{proof}

\subsection{Distribution of one random variable}
For the case where we restrict our attention to just one random variable $a\in\cA$ we introduce
the following notation.

\begin{notation}
Let $(\cA,\ff)$ be a HOPS and let $a\in\cA$.

1) For, $(\cV,\pi)\in\cPP(n)$, we will write
$$\ff^a(\cV,\pi):=\ff(\cV,\pi)\underbrace{[a,\dots,a]}_{\text{$n$-times}}$$
and
$$\KK^a(\cV,\pi):=\KK(\cV,\pi)\underbrace{[a,\dots,a]}_{\text{$n$-times}}.$$

2) A \emph{Young diagram} is a
$\lambda=(\lambda_{1}, \ldots,\lambda_{l} )$
for some $l\in\NN$ and $\lambda_1,\dots,\lambda_l\in\NN$ with $\lambda_1\geq\lambda_2
\geq\cdots\geq \lambda_l$. We put
$\vert\lambda\vert:=
\lambda_1+\cdots+\lambda_l$ (the total number of boxes of the Young diagram
$\lambda$). The set of all Young diagrams will be denoted by $\Young$.

3) The information about the higher order moments of $a$ can also be parametrized by Young
diagrams as follows: for $\lambda=(\lambda_{1},\ldots,\lambda_{l})$
we put
$$\ff^a(\lambda):=
\ff(1_{\vert\lambda\vert},\pi)\underbrace{[a,\dots,a]}_{\text{$n$-times}}=\ff_l(a^{\lambda_1},
\dots,a^{\lambda_l})$$ where $\pi$ is any permutation whose conjugacy class
corresponds to $\lambda$ (i.e., $\pi\in S_{\vert\lambda\vert}$ has cycles of length
$\lambda_1,\dots,\lambda_l$.
The collection of all higher order moments $\left(\ff^a(\lambda)\right)_{\lambda\in\Young}$
is called the \emph{(higher order) distribution of $a$}.

4) Similarly as for moments, we put
$$\KK^a(\lambda):=
\KK(1_{\vert\lambda\vert},\pi)\underbrace{[a,\dots,a]}_{\text{$n$-times}},$$
where $\pi$ is any permutation whose conjugacy class
corresponds to $\lambda$.

\end{notation}

\begin{remark}
For first and second order moments and cumulants,
we used in Section \ref{section2} also the following notations:
$$\alpha_n:=\ff^a(1_n,\gamma_n) \qquad \alpha_{m,n}^a:=\ff^a(1_{m+n},\gamma_{m,n}),$$
and
$$\KK^a_n:=\KK^a(1_n,\gamma_n) \qquad \KK_{m,n}:=\KK^a(1_{m+n},\gamma_{m,n}),$$
where $\gamma_n$ and $\gamma_{m,n}$ are permutations with one cycle and two cycles,
respectively.
\end{remark}

The vanishing of mixed cumulants translates in this framework
into the additivity of the cumulants for sums of free variables.

\begin{theorem}\label{thm:additivity-of-cumulants}
Let $(\cA,\ff)$ be a HOPS and $a,b\in\cA$ free of all orders. Then we have
$$\KK^{a+b}(\lambda)=\KK^a(\lambda)+\KK^b(\lambda)$$
for all $\lambda\in\Young$.
\end{theorem}

\begin{proof}
By the multilinearity of the cumulants and the vanishing of mixed cumulants for
free variable,
we have for any $n\in\NN$ and $\pi\in S_n$:
\begin{align*}
\KK^{a+b}(1_n,\pi)&=\KK(1_n,\pi)[a+b,\dots,a+b]\\
&=\KK(1_n,\pi)[a,\dots,a]+\KK(1_n,\pi)[b,\dots,b]\\
&=\KK^a(1_n,\pi)+\KK^b(1_n,\pi).
\end{align*}
\end{proof}


\section{Random matrices, Itzykson-Zuber integrals and higher order freeness}
\label{sec:RMT-IZ}

\subsection{Asymptotic higher order freeness of random matrices}

Let us now come back to our original motivation for our theory -- the asymptotic behavior of
random matrices. In order to reformulate our calculations from Section \ref{sec:correlation} in
our language of higher order freeness we still need to define the notion of ``asymptotic
freeness".

\begin{definition}\label{def:limithops}
1) Let $(\cA,\ff)$ and, for each $N\in\NN$,
$(\cA_{N},\ff^{(N)})$ be HOPSs. Let $I$ be an index set and for
each $i\in I$, $a_i\in\cA$ and $a_i^{(N)}\in\cA_N$ ($N\in\NN$).
We say that the family $(a_i^{(N)}\mid i\in I)$ converges, for
$N\to\infty$, to $(a_i\mid i\in I)$, denoted by
$$(a_i^{(N)})_{ i\in I} \tovert (a_i)_{i\in I},$$
if we have for all $n\in\NN$ and all polynomials $p_1,\dots,p_n$
in $\vert I\vert$-many non-commuting indeterminates that
\begin{multline}
$$\lim_{N\to\infty} \ff^{(N)}_n\left(p_1\bigl((a_i^{(N)})_{i\in I}\bigr),\dots,
p_n\bigl((a_i^{(N)})_{i\in I}\bigr)\right)\\= \ff_n\left(p_1\bigl((a_i)_{i\in I}\bigr),\dots,
p_n\bigl((a_i)_{i\in I}\bigr)\right).
\end{multline}

2) Let, for each $N\in\NN$, $(\cA_{N},\ff^{(N)})$ be
HOPSs. Let $I$ be an index set and, for each $i\in I$ and
$N\in\NN$, $a_i^{(N)}\in\cA_N$. We say that the sequence
of families $(a_i^{(N)})_{i\in I}$ has a \emph{limit
distribution of all orders} if there exists a HOPS
$(\cA,\ff)$ such that
$$(a_i^{(N)})_{ i\in I} \tovert (a_i)_{i\in I},$$
for some $a_i\in\cA$ ($i\in I$)

3) Let, for each $N\in\NN$, $(\cA_{N},\ff^{(N)})$ be HOPSs. Let $I$ be an index set and, for
each $i\in I$ and $N\in\NN$, $a_i^{(N)}\in\cA_N$. Let $I=I_1\cup\cdots \cup I_k$ be a
decomposition of $I$ into $k$ disjoint subsets. We say that the sets $\{a_i^{(N)}\mid i\in
I_1\},\dots,\{a_i^{(N)}\mid i\in I_k\}$ are \emph{asymptotically free of all orders} if there
exists a HOPS $(\cA,\ff)$ such that
$$(a_i^{(N)})_{ i\in I} \tovert (a_i)_{i\in I},$$
for some $a_i\in\cA$ ($i\in I$) and such that the sets $\{a_i\mid i\in I_1\},\dots, \{a_i\mid
i\in I_k\}$ are free of all orders in $(\cA,\ff)$.
\end{definition}

With this notation and by invoking Theorem \ref{thm:charact-freeness}
we can reformulate our main result on random matrices, Theorem \ref{distrib-product},
in the following form.

\begin{theorem}
Let $\cM_N:=M_N\otimes L^{\infty-}(\Omega)$ be an ensemble of $N\times
N$-random matrices. Define rescaled correlation functions $\tilde\ff^{(N)}=
(\tilde \ff_n^{(N)})_{n\in\NN}$ on
$\cM_N$ by ($n\in\NN$, $D_1,\dots,D_n\in\cM_N$)
\begin{equation}
\tilde\ff^{(N)}_n(D_1,\dots,D_n):= 
\cc_n(\Tr(D_1),\dots,\Tr(D_n))\cdot N^{2-n}.
\end{equation}
Assume that we have, for each $N\in\NN$, subalgebras $\cA_N,\cB_N\in \cM_N$
such that
\begin{enumerate}
\item
$\cA_N$ is a unitarily invariant ensemble,
\item
$\cA_N$ and $\cB_N$ are independent.
\end{enumerate}
Let $(A_i^{(N)})_{i\in I}$ be a family of elements in $(\cA_N,\tilde\ff^{(N)})$ which has a
higher order limit distribution and let $(B_j^{(N)})_{j\in J}$ ($N\in\NN$) be a family of
elements in $(\cB_N,\tilde\ff^{(N)})$ which has a higher order limit distribution. Then the
families $\{A_i^{(N)}\mid i\in I\}$ and $\{B_j^{(N)}\mid j\in J\}$ are asymptotically free of
all orders.
\end{theorem}

\subsection{Itzykson-Zuber integrals}

\begin{definition}
For $N\times N$ matrices
$A_N,B_N$ their \emph{Itzykson-Zuber integral} is defined as
the following function in $z\in\CC$:
$$\IZ(z,A_{N},B_{N}):=N^{-2}\log E(e^{zN\Tr(A_{N}UB_{N}U^*)}),$$
where $U$ denotes a Haar unitary $N\times N$-random matrix.
\end{definition}

Consider now a sequence of such matrices $A_N$ and $B_N$. Note that $A_N$ and $B_N$ are
non-random, thus all distributions of order higher than 1 vanish identically. If we assume that
$A_N$ and $B_N$ have a first order (eigenvalue) limit distribution for $N\to\infty$, then it is
known (see \cite{Collins2002}) that each Taylor coefficient about zero of $z\to
IZ(z,A_{N},B_{N})$ admits a limit as $N\to\infty$. Note that the effect of the Haar unitary
random matrix in the above Itzykson-Zuber integral was to make $A_N$ and $UB_NU^*$
asymptotically free of all orders. We show now that this kind of result extends also to the
case of random matrices $A_N$ and $B_N$, and that our theory allows to identify the limit of
the Taylor coefficients very precisely.

\begin{theorem}
Let $A=(A_N)_{N\in\NN}$ and $B=(B_N)_{N\in\NN}$ be two ensembles of $N\times N$-random matrices
which are asymptotically free of all orders with respect to the rescaled correlation functions
$\tilde \ff^{(N)}$. Denote the corresponding limiting distribution of $(A_N)_{N\in\NN}$ by
$\ff^a$ and the corresponding limit distribution of $(B_N)_{N\in\NN}$ by $\ff^b$. Then, as
formal power series in $z$, we have
\begin{equation}\label{izlim}
\lim_{N\to\infty} N^{-2}\log\EE[e^{zN\Tr(A_NB_N)}] = \sum_{n= 1}^\infty
\frac{z^{n}}{n!}\sum_{(\cV,\pi),(\cW,\sigma)\in \cPP(n)\atop
 (\cV,\pi)\cdot(\cW,\sigma)=(1_{n},e)}
 \KK^a(\cV,\pi)\cdot \ff^b(\cW,\sigma).
\end{equation}
\end{theorem}

\begin{proof}
Recall that the logarithm of the exponential generating series of the moments
of a random variable is the exponential generating series of the classical cumulants
of that variable. Thus we have
\begin{align*}
N^{-2}\cdot\log& \EE[e^{zNA_NB_N}]\\
&=N^{-2}\sum_{n=1}^\infty
\cc_n(N\Tr(A_NB_N),\dots,N\Tr(A_NB_N))\cdot\frac {z^n}{n!}\\
&= \sum_{n=1}^\infty N^{n-2}\cdot \ff^{(N)}
(1_n,e)[A_NB_N,\dots,A_NB_N]\cdot\frac {z^n}{n!}.
\end{align*}
By our assumption that $A_N$ and $B_N$ are asymptotically free with respect to
$\tilde \ff_n^{(N)}=N^{n-2}\ff_n^{N}$, this converges to
$$\sum_{n=1}^\infty \ff(1_n,e)[ab,\dots,ab]\cdot\frac {z^n}{n!},$$
where $a$ and $b$ are free of all orders with respect to $\ff$. Theorem
\ref{thm:charact-freeness} yields then the assertion.
\end{proof}

The main result of this part is the following theorem, which shows that the higher order
analogue of the
Itzykson-Zuber integral behaves like a kind of $R$-transform in one matrix argument if
the other is restricted to non-random matrices.
This result can be seen as a higher order version of a result of Zinn-Justin \cite{MR2000h:82043}.

\begin{theorem}
Let $C=(C_N)_{N\in\NN}$ be a sequence of non-random $N\times N$ matrices which has a first
order limit distribution. Then, for any sequence of $N\times N$-random matrices
$A=(A_N)_{N\in\NN}$ for which a higher order limit distribution exists, we define as a formal
power series in the moments of $C$
$$R^A(C):=\lim_{N\to\infty}N^{-2}\log E[e^{N\Tr(A_N U_NC_NU_N^*)}],$$
where $U_N$ are Haar unitary $N\times N$ random matrices which are independent
from $A_N$.
Then we have the following:

1) If $A=(A_N)_{N\in\NN}$ and $B=(B_N)_{N\in\NN}$ are random matrix ensembles which
are asymptotically free, then we
have as formal power series
\begin{equation}\label{eq:add-of-IZ-R}
R^{A+B}(C)=R^A(C)+R^B(C).
\end{equation}

2) More precisely, if we denote the limit moments of $C$ by
$x_k:=\lim_{N\to\infty}
\tr(C_N^k)$ ($k\in\NN$), then one has as a formal power series
\begin{equation}\label{eq:IZ-R}
R^A\left((x_k)_{k\in\NN}\right)=\sum_{\lambda\in\Young} \frac{x^{\lambda}c_{\lambda}}{|\lambda| !}
\kk^A(1_{|\lambda|},\lambda),
\end{equation}
where
$c_{\lambda}$ is the number of permutations in the conjugation class in
${S}_{|\lambda |}$ corresponding to $\lambda$, and where
$$x^\lambda:=x_1^{\lambda_1}x_2^{\lambda_2}\cdots,\qquad
\text{for $\lambda=(\lambda_1\geq \lambda_2\geq\dots)$.}$$

\end{theorem}

\begin{proof}
In order to get \eqref{eq:IZ-R} (which implies \eqref{eq:add-of-IZ-R})
we have to specialize \eqref{izlim} to
the situation that $B=C$ has all moments of order higher than 1 equal to zero.
This means that $\ff^b(\cW,\sigma)=0$ unless $\cW=0_\sigma$, in which case we
have
$$\ff^b(0_\sigma,\sigma)=x^\lambda,$$
where $\lambda$ is the Young diagram corresponding to $\sigma$.
\eqref{eq:IZ-R} follows then from the simple
observation that, for fixed $\sigma$,
the only solution of $(\cV,\pi)\cdot (0_\sigma,\sigma)=(1_n,e)$ is
given by $(\cV,\pi)=(1_n,\sigma^{-1})$.
\end{proof}

This theorem tells us that Itzykson-Zuber type of integrals contain
the whole data to linearize higher order freeness.
Also, the technology introduced in this paper, especially
Theorem \ref{second-R}
gives methods to refine asymptotics of spherical integrals \`a la Guionnet and
Ma\"\i{}da (see \cite{GuionnetMaidaRtransform}).

As a foreshadowing of such applications we close with the
following proposition. More details will be provided in a
forthcoming paper.

Let $C_N$ be a non-random matrix of rank $2$ with eigenvalues $x$ and $y$ (which will
be considered as indeterminates in the following).
Then we have that $\ff^{(N)}(\cW,\sigma)[C,\dots,C]$ is only different from zero
for $\cW=0_\sigma$, in which case it is
$$\ff^{(N)}(0_\sigma,\sigma)[C,\dots,C]=(x^{\lambda_1}+y^{\lambda_1})\cdots
(x^{\lambda_l}+y^{\lambda_l}),$$
where $\lambda(\sigma):=
\lambda=(\lambda_1,\dots,\lambda_l)$ is the Young diagram encoding the conjugacy
class of $\sigma$.
Thus Theorem \ref{distrib-product} yields in this case
\begin{align*}
&N^{-1}\log E[e^{N\Tr(A_NUB_NU^*}]=\sum_{n=1}^\infty
N^{n-1}\ff^{(N)}(1_n,e)[A_NC_N,\dots,A_NC_N]\cdot
\frac{z^n}{n!}\\
&=\sum_{n=1}^\infty\frac{z^n}{n!} \sum_{\sigma\in
S_n}\KK^{(N)}(1_n,\sigma^{-1})[A_N,\dots,A_N]\cdot
(x^{\lambda_1(\sigma)}+y^{\lambda_1(\sigma)})\cdots
(x^{\lambda_l(\sigma)}+y^{\lambda_l(\sigma)})
\end{align*}
Now we invoke the assumption on the existence of a limit distribution for $A_N$ and
the fact (which we have never used up to now) that the orders in the Weingarten function
and thus also in $\KK^{(N)}$ decreases in steps of 2. This
allows us to recognize the two leading orders of the above
quantity and we get the following result.

\begin{proposition}
Let $C_{N}$ be a matrix of rank $2$ and eigenvalues $x,y$. Consider a sequence of $N\times
N$-random matrices $A=(A_N)_{N\in\NN}$ for which a higher order limit distribution $\ff^a$
exists. Denote the limiting first and second order cumulants by
$$\KK_n:=\KK^a(1_n,\gamma_n) \qquad \KK_{m,n}:=\KK
^a(1_{m+n},\gamma_{m,n}),$$
where $\gamma_n$ and $\gamma_{m,n}$ are permutations with one cycle and two cycles,
respectively.
Then, as a formal power series in $x,y$ we have
\begin{multline*}
N^{-1}\lim_{N\to\infty}\log \EE[e^{N\Tr A_{N}UC_{N}U^{*}}] \\=\sum_{n\geq 1} \frac
{\kk_n}{n}\cdot (x^n+y^n) +N^{-1}\sum_{m,n\geq 1}\frac {\kk_{m,n}}{mn}\cdot
(x^{m}+y^{m})(x^{n}+y^{n}) +O(N^{-2}).
\end{multline*}
\end{proposition}

This expansion extends results of \cite{Collins2002}. It is more general than that obtained in
\cite{GuionnetMaidaRtransform}, because it also handles the case when $A_{N}$
has asymptotic fluctuations. However, unlike in \cite{GuionnetMaidaRtransform},
the convergence that we obtain is formal and it would be very interesting to
check if it still holds at an analytic level.

\section{Appendix: Surfaced permutations}
\label{sec:surfaced}

\newcommand{\support}[1]{|#1|}

In this appendix we will present a more geometrical view on partitioned
permutations.
As we shall see in the following,
partitioned permutations are just special cases of ``surfaced
permutations"; in particular the results of this article can be
equivalently formulated in the language of surfaced permutations. On
the other hand, for the purpose of this article we do not need
anything more than just partitioned permutations and the Reader not
interested in surfaced permutations may skip this Section without
much harm.

\subsection{Motivations}

Our goal is to study factorizations of permutations, i.e.\ solutions
$(\pi_1,\dots,\pi_k)$ of the equation
$$ \gamma=\pi_1 \cdots \pi_k,$$
where $\gamma\in S_n$ is some fixed permutation and
$\pi_1,\dots,\pi_k\in S_n$ are subject to some additional
constraints, depending on a particular context. Typically, one of
these constrains concerns $\vert\pi_1\vert+\cdots+\vert\pi_k\vert$,
the other one concerns the orbits of the action of
$\pi_1,\dots,\pi_k$.

It would be very useful to equip permutations $\pi_1,\dots,\pi_k$
with some additional structure in such a way that the product
$\tilde{\pi}_1\cdots\tilde{\pi}_k$ of the resulting enriched
permutations $\tilde{\pi}_1,\dots,\tilde{\pi}_k$ would carry both
the information about the product $\pi_1 \cdots \pi_k $ of
permutations and the information about
$\vert\pi_1\vert+\cdots+\vert\pi_k\vert$. As we shall see in the
following, surfaced permutations provide an appropriate tool.

\subsection{Definition}
\begin{figure}[tb]
\includegraphics{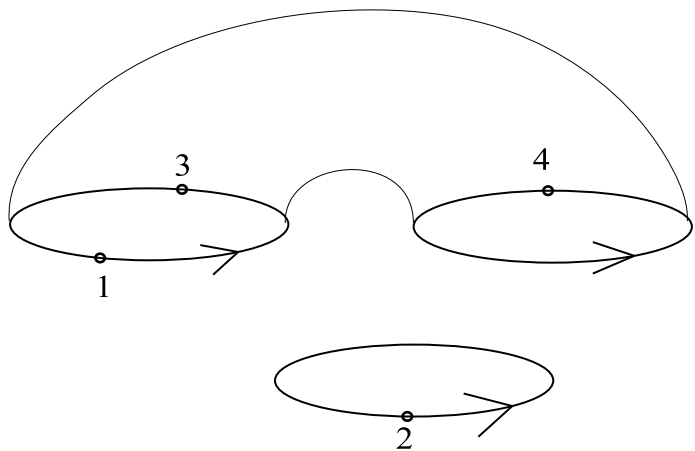}
\caption{Example of a surfaced permutation. Its support is equal to
$(1,3)(2)(4)\in S_4$. This surfaced permutation corresponds to a
partitioned permutation $\big(\{1,3,4\}\{2\},(1,3)(2)(4) \big)$}.
\label{fig:young1}
\end{figure}

We say that $\sigma=(S,j)$ is a \emph{surfaced permutation} of some
finite set $A$ if $S$ is a two--dimensional surface with a fixed
orientation and with a boundary $\partial S$ and if $j:A\rightarrow
\partial S$ is a injection. We can think about the information
carried by $j$ as follows: some of the points on the boundary
$\partial S$ are distinguished and carry different labels from the
set $A$. We also require that every connected component of $\partial
S$ carries at least one distinguished point. An example of a
surfaced permutation is presented on Figure \ref{fig:young1}.

We identify surfaced permutations $(S_1,j_1)$, $(S_2,j_2)$ of the
same set $A$ if there exists a orientation preserving
homeomorphism $f:S_1\rightarrow 
S_2$ such that $f\circ j_1=j_2$. The set of surfaced
permutations of set
$\{1,\dots,n\}$ will be denoted by $\Sur_n$.


\subsection{Surfaced permutations and the usual permutations}
%

Let $(S,j)\in\Sur_n$; the boundary $\partial S$ with the inherited
orientation from $S$ is just a collection of oriented circles with
some distinguished points labeled $1,\dots,n$ marked on them. In
this way we can define a permutation $\sigma\in S_n$, called 
the support of $(S,j)$, the cycles of which correspond to
connected components of $\partial S$, as it can be seen on Figure
\ref{fig:young1}. It is therefore a good idea to think that a
surfaced permutation is just a (usual) permutation $\sigma\in
S_n$ equipped with some additional information carried by the
surface $S$.

A surfaced permutation $(S,j)\in\Sur_n$ can be uniquely specified
(up to the equivalence relation) by its support $\sigma\in S_n$ and
by specifying the shape of the connected components of $S$. The
latter information is given by an equivalence relation on cycles of
$\sigma$ (each class corresponds to a connected component of $S$)
and furthermore for each class of this relation we should specify
the genus of the corresponding connected component of $S$. Above it
should be understood that the genus of a surface $S$ with a boundary
is by definition equal to the genus of a surface $S'$ without
boundary obtained from $S$ by gluing a disc to every connected
component of $\partial S$; for example both a disc and the lateral
surface of a cylinder have genus zero.

\subsection{Surfaced permutations and partitioned
permutations} \qquad

\noindent
Among surfaced permutations a special class will
be very important for our purposes, namely surfaced permutations
$(S,j)$ such that each connected component of $S$ has genus
zero. It is easy to see that there is a bijection between such
surfaced permutations $(S,j)$ and partitioned permutations
$(\cV,\sigma)$ given as follows:
$\sigma$ is the support of $(S,j)$ and $\cV$ is the partition given
by connected components of $S$.

\subsection{Products of surfaced permutations}

\begin{figure}[tb]
  \includegraphics{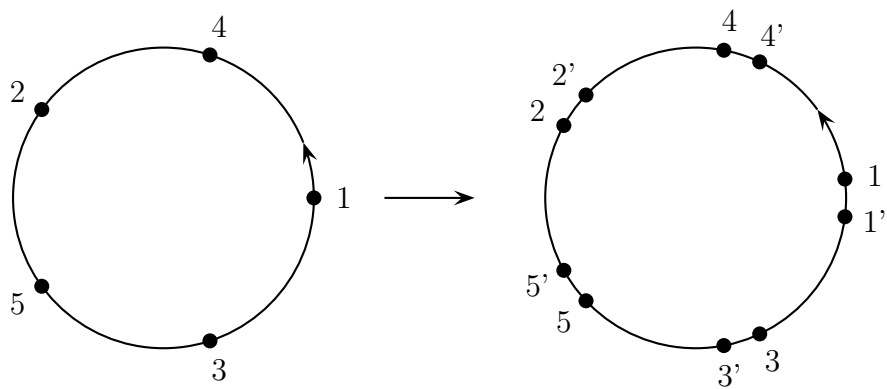}
  \caption{Convention for splitting labels.} \label{fig:split}
\end{figure}
Let surfaced permutations $(S_1,j_1),(S_2,j_2)\in\Sur_n$ be given.
On the boundary of $S_2$ there are marked points labeled by numbers
$1,\dots,n$; let us split every marked point $k$ into a consecutive
pair of points $k$ and $k'$, as it is presented on the example from
Figure \ref{fig:split}. In the second step, for each
$k\in\{1,\dots,n\}$ we glue a small neighborhood of the vertex $k\in
\partial S_1$ to a small neighborhood of the vertex $k'\in\partial
S_2$ in such a way that the orientations of $S_1$ and $S_2$
coincide. In this way we obtain a new surface $S$ which has marked
points on its boundary $\partial S$ and these are exactly the
vertices from $\partial S_2$ labeled $1,\dots,n$; we denote the
resulting surfaced permutation by $(S,j)$ and we call it a product
$(S_1,j_1)(S_2,j_2)$ of the original surfaced permutations. This
choice of gluing surfaces $S_1$ and $S_2$ implies that the support
of $(S_1,j_1)(S_2,j_2)$ is equal to the product of the support of
$(S_1,j_1)$ and the support of $(S_2,j_2)$.

It is not difficult to explain now the definition of the product of
partitioned permutations (Definition \ref{def:mult}): we treat
partitioned permutations as surfaced permutations and compute their
product; if the genus of the resulting surface is zero we can
identify it with another partitioned permutation, otherwise we set
the product to be zero.

It is not difficult to show that for surfaced permutations
the product is associative and the associativity of the product
of partitioned permutations is a simple corollary.

%
%
%
%

\bibliographystyle{alpha}

\begin{thebibliography}{CMSS05}

\bibitem[BS04]{BS} Z.~D.~Bai and J.~Silverstein: CLT for linear
spectral statistics of large-dimensional sample covariance
matrices. {\em Ann. Prob.}, 32: 533--605, 2004.

\bibitem[Bia97]{Biane1997crossings}
P.~Biane:
\newblock Some properties of crossings and partitions.
\newblock {\em Discrete Math.}, 175(1-3):41--53, 1997.

\bibitem[BMS00]{MR1761777}
M.~Bousquet-M\'elou and G.~Schaeffer.
\newblock Enumeration of planar constellations.
\newblock {\em Adv. in Appl. Math.}, 24(4):337--368, 2000.

\bibitem[BZ93]{bz}
\'E. Br\'ezin, A. Zee,
Universality of the correlations between eigenvalues of
large random matrices, {\em Nuclear Phys.} B 402 (1993),
no. 3, 613--627.

\bibitem[Col03]{Collins2002}
B.~Collins:
\newblock Moments and cumulants of polynomial random variables on unitary
  groups, the {I}tzykson-{Z}uber integral, and free probability.
\newblock {\em Int. Math. Res. Not.}, (17):953--982, 2003.

\bibitem[CM\'SS]{HigherOrderFreeness4}
B.~Collins, J.~Mingo, P.~\'Sniady, and R.~Speicher:
\newblock Second order freeness and fluctuations of random matrices: {IV.}
  {A}symptotic freeness of {J}ucys-{M}urphy elements,
\newblock in preparation.

\bibitem[C\'S04]{CollinsSniady2004}
B.~Collins and P.~\'Sniady:
\newblock Integration with respect to the {H}aar measure on unitary, orthogonal
  and symplectic groups.
\newblock {\em Comm. Math. Phy.}, 264, (2006), 773 - 795.

\bibitem[Dia03]{Diaconis}
P.~Diaconis: Patterns in Eigenvalues: The 70th Josiah Willard Gibbs Lecture. {\em Bulletin of
the AMS}, 40: 155-178, 2003.

\bibitem[FMP78]{fmp}
J. B. French, P. A. Mello, A. Pandey, 
Statistical properties of many-particle spectra. II.
Two-point correlations and fluctuations, {\em Ann. of
Physics}, 113 (1978), no. 2, 277--293.

\bibitem[GM05]{GuionnetMaidaRtransform}
A.~Guionnet and M.~Ma\"\i{}da:
\newblock A {F}ourier view on the {$R$}-transform and related asymptotics of
  spherical integrals.
\newblock {\em J. Funct. Anal.}, 222(2):435--490, 2005.

\bibitem[Joh98]{MR1487983}
K.~Johansson:
\newblock On fluctuations of eigenvalues of random {H}ermitian matrices.
\newblock {\em Duke Math. J.}, 91(1):151--204, 1998.

\bibitem[Jon99]{Jonesplanar}
V.~F.~R.~Jones:
\newblock Planar algebras. {I}.
\newblock {\tt math.QA/9909027},
1999, 122 pp.

\bibitem[KKP95]{kkp}
A. Khorunzhy, B. Khoruzhenko, L. Pastur,  On the $1/N$
corrections to the Green functions of random matrices with
independent entries, {\em J. Phys.} A 28 (1995), L31--L35.

\bibitem[Kre72]{Kreweras}
G.~Kreweras:
\newblock Sur les partitions non crois\'ees d'un cycle.
\newblock {\em Discrete Math.}, 1(4):333--350, 1972.

\bibitem[LS59]{LS}
V.~P.~Leonov and A.~N.~Shiryaev:
On a method of semi-invariants.
{\em Theory of Probability and its Applications}, 4, 319--329, 1959.

\bibitem[MN04]{MingoNica2004annular}
J.~Mingo and A.~Nica:
\newblock Annular noncrossing permutations and partitions, and second-order
  asymptotics for random matrices.
\newblock {\em Int. Math. Res. Not.}, (28):1413--1460, 2004.

\bibitem[M\'SS04]{HigherOrderFreeness2}
J.~Mingo, P.~\'Sniady, and R.~Speicher.
\newblock Second order freeness and fluctuations of random matrices: {II.}
  {U}nitary random matrices.
\newblock To appear in {\em Adv. in Math.}.

\bibitem[MS04]{HigherOrderFreeness1}
J.~Mingo and R.~Speicher:
\newblock {Second Order Freeness and Fluctuations of Random Matrices: I.
  Gaussian and Wishart matrices and Cyclic Fock spaces}.
\newblock{\em J. Funct. Anal.}, 235, 2006, pp. 226-270.

\bibitem[NSp97]{NicaSpeicher1997Fourier}
A.~Nica and R.~Speicher:
\newblock A ``{F}ourier transform'' for multiplicative functions on
  non-crossing partitions.
\newblock {\em J. Algebraic Combin.}, 6(2):141--160, 1997.

\bibitem[NSp06]{TheBook}
A.~Nica and R.~Speicher:
Lectures on the Combinatorics of Free Probability, 
London Mathematical Society Lecture Note Series, New York :
Cambridge University Press, to appear.

\bibitem[OZ84]{OBrienZuber84}
K.~H. O'Brien and J.-B. Zuber.
\newblock A note on {${\rm U}(N)$} integrals in the large {$N$} limit.
\newblock {\em Phys. Lett. B}, 144(5-6):407--408, 1984.

\bibitem[Oko00]{Okounkov}
A.~Okounkov: Random matrices and random permutations. {\em
Int. Math. Res. Not.}, 2000, no. 20: 1043--1095, 2000.


\bibitem[Rad04]{Radulescu}
F. Radulescu: Combinatorial aspectes of Connes's embedding conjecture and asymptotic
distribution of traces of products of unitaries. Preprint, 2004. math.0A/0404308.

\bibitem[Rot64]{rota}  G.-C. Rota,
On the foundations of combinatorial
theory. I. Theory of M\"obius functions, {\em Z.
Wahrscheinlichkeitstheorie und Verw. Gebiete}, 2, (1964),
340--368.

\bibitem[\'Sni05]{Sniady2005GaussuanFluctuationsofYoungdiagrams}
P.~\'Sniady.
\newblock {Gaussian fluctuations of characters of symmetric groups and of Young
  diagrams}.
\newblock Preprint {\tt arXiv:math.CO/0501112}, 2005.

\bibitem[Spe94]{Speicher1994}
R.~Speicher:
\newblock Multiplicative functions on the lattice of noncrossing partitions and
  free convolution.
\newblock {\em Math. Ann.}, 298(4):611--628, 1994.

\bibitem[VDN92]{VoiculescuDykemaNica}
D.~V. Voiculescu, K.~J. Dykema, and A.~Nica.
\newblock {\em Free random variables}.
\newblock American Mathematical Society, Providence, RI, 1992.

\bibitem[Voi85]{Voiculescu-first}
D. Voiculescu: Summetries of some reduced free product $C^*$-algebras. In {\it Operator
Algebras and their Connections with Topology and Ergodic Theory}  (Lecture Notes in
Mathematics, vol. 1132, Springer-Verlag): 556-588, 1985.


\bibitem[Voi86]{Voiculescu-R}
D. Voiculescu: Addition of certain non-commuting random variables. {\em J. Funct. Anal.}, 66:
323--346, 1986.

\bibitem[Voi91]{Voiculescu1991}
D.~Voiculescu:
\newblock Limit laws for random matrices and free products.
\newblock {\em Invent. Math.}, 104(1):201--220, 1991.

\bibitem[Wei78]{Weingarten1976}
D.~Weingarten:
\newblock Asymptotic behavior of group integrals in the limit of infinite rank.
\newblock {\em J. Mathematical Phys.}, 19(5):999--1001, 1978.

\bibitem[ZJ99]{MR2000h:82043}
P.~Zinn-Justin.
\newblock Adding and multiplying random matrices: a generalization of
  {V}oiculescu's formulas.
\newblock {\em Phys. Rev. E (3)}, 59(5, part A):4884--4888, 1999.

\end{thebibliography}

\end{document}